\documentclass[final,3p,times]{elsarticle}

\usepackage{amsmath}
\usepackage{amsfonts}
\usepackage{amssymb}
\usepackage{amsthm}
\usepackage[mathscr]{euscript}
\usepackage{graphicx}
\usepackage{enumitem}
\usepackage{framed}
\usepackage{booktabs}
\usepackage{tabularx}
\usepackage{multirow}
\usepackage{caption}
\usepackage[table]{xcolor}
\usepackage{mathtools}
\usepackage{todonotes}
\usepackage[ruled,noend,noline,slide]{algorithm2e}
\usepackage{algpseudocode}

\usepackage{hyperref}
\usepackage{xcolor}
\usepackage{etoolbox}

\definecolor{shadecolor}{RGB}{230,230,230}
\definecolor{shade_0}{RGB}{255,255,255}
\definecolor{shade_1}{RGB}{240,240,240}
\definecolor{shade_2}{RGB}{220,220,220}
\definecolor{shade_3}{RGB}{200,200,200}
\definecolor{shade_4}{RGB}{180,180,180}
\emergencystretch=\maxdimen
\newcommand{\del}{\partial}
\newcommand{\bigO}{\mathcal{O}}

\theoremstyle{definition}
\newtheorem{mylist}{List}
\newtheorem{thm}{Theorem}
\newtheorem{cor}{Corollary}

\journal{Mathematics and Computers in Simulation}

\usepackage{empheq}

\makeatletter
\let\cfbox\fbox
\let\c@frameb@x\@frameb@x
\pretocmd{\cfbox}
{\leavevmode\begingroup\colorlet{currentcolor}{.}\color{red}}
{}{}
\patchcmd\cfbox{\@frameb@x}{\c@frameb@x}{}{}
\patchcmd{\c@frameb@x}{\box\@tempboxa}
{\color{currentcolor}\box\@tempboxa}
{}{}
\apptocmd{\c@frameb@x}{\endgroup}{}{}
\makeatother

\begin{document}

\begin{frontmatter}



\title{On the Spatial and Temporal Order of Convergence of Hyperbolic PDEs}


\author[label1,label2]{Siddhartha Bishnu}
\author[label1]{Mark Petersen}
\author[label2]{Bryan Quaife}

\address[label1]{Computational Physics and Methods Group, Los Alamos National Laboratory}
\address[label2]{Department of Scientific Computing, Florida State University}

\begin{abstract}
In this work, we determine the full expression for the global truncation error of hyperbolic partial differential equations (PDEs). In particular, we use theoretical analysis and symbolic algebra to find exact expressions for the coefficients of the generic global truncation error 
\begin{equation}
\nonumber
\hat{\tau}_G = \bigO\left({\Delta x}^{\alpha}\right) + \Delta t \bigO\left({\Delta x}^{\alpha}\right) + {\Delta t}^2 \bigO\left({\Delta x}^{\alpha}\right) + \cdots + {\Delta t}^{\beta-1} \bigO\left({\Delta x}^{\alpha}\right) + \bigO\left({\Delta t}^{\beta}\right) \approx \bigO\left({\Delta x}^{\alpha}\right) + \bigO\left({\Delta t}^{\beta}\right), \text{ for $\Delta t \ll 1$},
\end{equation}
where $\Delta x$ and $\Delta t$ denote the cell width and the time step size, and $\alpha$ and $\beta$ represent the orders of the spatial and temporal discretizations. Our analysis is valid for any hyperbolic PDE, be it linear or non-linear, and employing finite difference, finite volume, or finite element discretization in space, and advanced in time with a predictor-corrector, multistep, or a deferred correction method, belonging to the Method of Lines.

Furthermore, we discuss the practical implications of this analysis. If we employ a stable numerical scheme and the orders of accuracy of the global solution error and the global truncation error agree, we make the following asymptotic observations: (a) the order of convergence at constant ratio of $\Delta t$ to $\Delta x$ is governed by the minimum of the orders of the spatial and temporal discretizations, and (b) convergence cannot even be guaranteed under only spatial or temporal refinement. An implication of (a) is that it is impractical to invest in a time-stepping method of order higher than the spatial discretization. In addition to (b), we demonstrate that under certain circumstances, the error can even monotonically increase with refinement only in space or only in time, and explain why this phenomenon occurs. To verify our theoretical findings, we conduct convergence studies of linear and non-linear advection equations using finite difference and finite volume spatial discretizations, and predictor-corrector and multistep time-stepping methods. Finally, we study the effect of slope limiters and monotonicity-preserving strategies on the order of accuracy.
\end{abstract}


\end{frontmatter}


\section{Introduction}

A course on numerical methods for applied scientists and engineers is typically taught with standard predictor-corrector and multistep time-stepping methods applied to ODEs in one chapter, followed by spatial discretization operators of PDEs in another. In real-world applications, the discretization of a PDE model of a physical phenomenon consists of both spatial and temporal components. As a result, the order of convergence of the PDE depends on both the cell width $\Delta x$ and the time step $\Delta t$. In past decades, when low-order finite difference and finite volume discretizations were common, the combined effect of the spatial and temporal discretization errors received little attention. Now that high-order discretizations are ubiquitous in PDE applications, the combined effect of $\Delta x$ and $\Delta t$ on the local truncation error has important consequences. After surveying standard numerical analysis textbooks including~\citet{burden1985numerical}, \citet{chapra2010numerical}, \citet{cheney2012numerical}, \citet{iserles2009first}, \citet{strikwerda2004finite}, we found this topic to be missing. Therefore, this paper investigates the simultaneous effects of $\Delta x$ and $\Delta t$ on the local truncation error of hyperbolic PDEs for several standard spatial and temporal discretizations.

There are very few references in the scientific journals as well. One of the best available derivations we found is on a webpage~\citep{Langtangen15}. However, the author did not express the temporal and mixed derivatives of the dependent variable of the PDE as functions of quantities at the current time level, which are known a priori, and this is a key step to arrive at our final form of the local truncation error. Only two published papers on the order of convergence in space and time could be found. \citet{love2013convergence} demonstrate that the order of convergence for finite difference approximations of PDEs may be different for refinement only in time than for refinement in both space and time. They employ linear algebra, especially matrix analysis, to determine the error and convergence rates of linear PDEs with refinement only in time, and substantiate the results with numerical experiments. In contrast, we derive the truncation error of linear and non-linear hyperbolic PDEs from first principles, using theoretical analysis and symbolic algebra, and verify our results numerically using finite difference and finite volume methods, for refinement in both space and time, only in space, and only in time. Our approach is more general and straightforward, and applicable to any discretization. \citet{jeong2019verification} perform convergence tests of three parabolic PDEs---the heat equation, the Allen-Cahn equation, and the Cahn-Hilliard equation---by refining the spatial and temporal steps separately and together. They urge their readers to be cautious of the varying order of convergence, but do not perform any detailed examination of the underlying cause of this discrepancy. The analysis presented in our paper explains the order of the convergence results in \citet{jeong2019verification}. 

There exists numerical methods which treat space and time together such as Lax-Wendroff and Cauchy-Kowalevski. \citet{qiu2005discontinuous} apply a Lax-Wendroff time discretization procedure to the discontinuous Galerkin method for solving hyperbolic conservation laws. Arbitrary derivatives (ADER) time-stepping methods involving space-time basis functions also belong to this category. \citet{normanhigh} develops a high-order WENO-limited finite-volume algorithm for modeling atmospheric flow using the ADER-differential transform time discretization. Models like MITgcm~\citep{marshall1997finite} use direct space-time methods for modeling non-linear advection and other physical phenomena in the ocean. Incremental remapping, cell-integrated, and flux-form semi-Lagrangian approaches are some other instances where space and time considerations are convolved together. Even though our theory holds for Lax-Wendroff methods, we do not test our theory for more sophisticated methods that treat space and time together.

Finally, our theory is applicable only when the global solution error at a time horizon has the same order of accuracy as the global truncation error. If the order of accuracy of the global solution error exceeds that of the global truncation error, superconvergence (or supraconvergence) is observed. Under such circumstances, consistency is not a necessary condition for convergence e.g.~\citep{cockburn1997priori}. \citet{cao2018some} discuss some recent developments in superconvergence of Discontinuous Galerkin methods for time-dependent PDEs, while \citet{peixoto2016accuracy} mentions superconvergence effects in the accuracy analysis of mimetic finite volume operators on geodesic grids. Our theory considers neither superconvergence nor the order reduction arising from the application of (a) dissipative Riemann solvers to compute numerical fluxes, and (b) monotone slope-limiting strategies to ensure oscillation-free profiles. Both techniques are common practices in finite volume methods.

\subsection{Notation}

Consider a differential equation in abstract form
\begin{equation}
L(u) = f, \label{AbstractEquation}
\end{equation}
where $L(u)$ is a function of the dependent variable $u$ and its derivatives with respect to the independent variables, and the forcing term $f$ is only a function of the independent variables. We assume that $L$ includes the boundary condition. Discretizing~\eqref{AbstractEquation} results in the discrete difference equation
\begin{equation}
L_{\Delta}(u_{\Delta}) = f_{\Delta}, \label{DifferentialEquation}
\end{equation}
where the subscript $\Delta$ represents the set of discretization parameters for the spatial and temporal grids. The exact solution $u$ can be distinguished from the numerical solution $u_{\Delta}$, which is a function of $\Delta$. In a finite difference approximation, $u_{\Delta}$ is computed at a set of grid points in space and time, and the abstract function $L_\Delta(u_\Delta)$ at any grid point in space and time typically consists of algebraic equations in $u_\Delta$ at that point and some neighboring ones.

The local truncation error of the difference equation is the residual when substituting the exact solution $u$ into the difference equation~\eqref{DifferentialEquation} 
\begin{equation}
\tau_\Delta = L_{\Delta}(u) - f_{\Delta}. \label{DifferenceEquation}
\end{equation}
A small value of $\tau_\Delta$ indicates that the difference equation closely resembles the differential equation, thereby implying proximity of $u_{\Delta}$ to $u$. A discretization method of a PDE is said to be consistent if the local truncation error $\tau_{\Delta} \to 0$ as $\Delta \to 0$. In the special case when $L_\Delta$ is linear, the error $e_\Delta = u - u_{\Delta}$ can be written in terms of the residual since 
\begin{equation}
L_{\Delta} e_{\Delta} = L_{\Delta}\left(u - u_{\Delta}\right) = L_{\Delta} u - L_{\Delta} u_{\Delta} = L_{\Delta} u - f_{\Delta} = \tau_{\Delta} \hspace{1mm} \text{ i.e. } \hspace{1mm} e_{\Delta} = L_{\Delta}^{-1} \tau_{\Delta},
\end{equation}
which represents the relationship between the actual error and the local truncation error of the difference equation. Taking norms of both sides,
\begin{equation}
\|e_{\Delta}\| = \|L_{\Delta}^{-1} \tau_{\Delta}\| \le \|L_{\Delta}^{-1}\|\hspace{1mm}\|\tau_{\Delta}\|. \label{HalfOfLaxEquivalenceTheoremProof}
\end{equation}
The numerical solution of the PDE converges to its exact analytical counterpart if $\|e_{\Delta}\| \to 0$ as $\Delta \to 0$. 

The numerical scheme is stable if $\|L_{\Delta}^{-1}\|$ is bounded independent of $\Delta$. Therefore, if the numerical solution is stable, then consistency implies convergence using~\eqref{HalfOfLaxEquivalenceTheoremProof}. This in fact constitutes half of the proof of the Lax Equivalence Theorem (\citet{lax1956survey}), which states that a consistent finite difference method for a well-posed linear initial value problem is convergent if and only if it is stable.

\subsection{Outline of the Paper}
Our analysis consists of determining the local truncation error of a set of ordinary differential equations (ODEs) and PDEs from first principles for a variety of spatial and temporal discretizations. Even though the derivation of the local truncation error of the characteristic ODE can be found in textbooks, it is instructive to first introduce our method and notations with the generic ODE, which are present in Section~\ref{sec:odes}. Then, we extend the analysis to the generic hyperbolic PDE in Section~\ref{sec:pdes}, the main subject of this paper. Section~\ref{sec:numerical_results} contains numerical results which demonstrate our theoretical findings. Our conclusions are presented in Section~\ref{sec:conclusion}. The supplementary documents contain two appendices and error expansions for high-order spatial discretizations and non-linear hyperbolic PDEs. Appendix~A outlines the algorithm for implementing Williamson's low-storage third-order Runge-Kutta method~\citep{williamson1980low}, and Carpenter and Kennedy's low-storage fourth-order Runge-Kutta method~\citep{carpenter1994fourth} to advance an ODE or a PDE over one time step, along with the important coefficients. Finally, Appendix B lists the leading order terms in the local truncation error of an ODE and an inhomogeneous variable-coefficient advection equation for a variety of time-stepping methods. 
%
%

\section{Ordinary Differential Equations}
\label{sec:odes}
We start by considering the generic first-order ODE
\begin{equation}
u_t = \mathcal{F}(u,t).
\label{ODE1D}
\end{equation}
The right-hand side $\mathcal{F}$ can be a linear or non-linear function of both $u$ and $t$. Given the exact solution $u^n$ at time level~$t^n$, the {\em local truncation error of the difference equation} at time level $t^{n+1} = t^n + \Delta t$ is 
\begin{equation}
\tau^{n+1} = L_{\Delta} u^n - f_{\Delta}^n. \label{ODE1DTruncationErrorDifferenceEquation_1}
\end{equation}
For an ODE, $\Delta = \Delta t$, and for any time-stepping method, \eqref{ODE1DTruncationErrorDifferenceEquation_1} can be written as 
\begin{equation}
\tau^{n+1} = \frac{1}{\Delta t} \left(u^{n+1} - \hat{u}^{n+1}\right), \label{ODE1DTruncationErrorDifferenceEquation_2}
\end{equation}
where $u^{n+1}$ is the exact solution at time $t^{n+1}$, and $\hat{u}^{n+1}$ is its numerical counterpart as a function of quantities known at time level $t^n$. To maintain clarity of the presentation, we call the {\em local truncation error of the numerical solution} the {\em local truncation error}, which is defined as
\begin{equation}
\hat{\tau}^{n+1} = \Delta t \tau^{n+1} = \Delta t \left(L_{\Delta} u^n - f_{\Delta}^n\right) = u^{n+1} - \hat{u}^{n+1}. \label{ODE1DTruncationErrorNumericalSolution_1}
\end{equation}
Calculating the local truncation error starts by expanding each term in~\eqref{ODE1DTruncationErrorNumericalSolution_1} about a common center. For example, the Taylor expansion of $u^{n+1}$ centered at $t^n$ is
\begin{equation}
u^{n+1} = u^n + \Delta t u_t^n + \frac{\Delta t^2}{2} u_{tt}^n + \frac{\Delta t^3}{6} u_{ttt}^n + \frac{\Delta t^4}{24} u_{tttt}^n + \cdots. \label{ODE1DExactSolutionAtTimeLevelNPlusOne_1}
\end{equation}
Then, the time derivatives in~\eqref{ODE1DExactSolutionAtTimeLevelNPlusOne_1} are written in terms of the $u$ and $t$ derivatives of $\mathcal{F}$ by using the chain rule
\begin{subequations} \label{TimeDerivativesOfDependentVariableODE1D}
\begin{align} 
u_{t} &= \mathcal{F} \equiv \mathcal{F}^{(1)}, \\
u_{tt} &= \mathcal{F}_t + \mathcal{F}_u \mathcal{F} \equiv \mathcal{F}^{(2)}, \\
u_{ttt} &= \mathcal{F}_{tt} + 2 \mathcal{F}_{ut} \mathcal{F} + \mathcal{F}_u \mathcal{F}_t + \mathcal{F}_{uu} \mathcal{F}^2 + \mathcal{F}_u^2 \mathcal{F} \equiv \mathcal{F}^{(3)}, \\
u_{tttt} &= \mathcal{F}_{ttt} + 3 \mathcal{F}_{utt} \mathcal{F} + 3 \mathcal{F}_{ut} \mathcal{F}_t + 5 \mathcal{F}_{ut} \mathcal{F}_u \mathcal{F} + 3 \mathcal{F}_{uu} \mathcal{F} \mathcal{F}_t + \mathcal{F}_u \mathcal{F}_{tt} \nonumber \\
&\hspace{0.325cm}+3 \mathcal{F}_{uut} \mathcal{F}^2 + 4 \mathcal{F}_{uu} \mathcal{F}_u \mathcal{F}^2 + \mathcal{F}_{uuu} \mathcal{F}^3 + \mathcal{F}_u^2 \mathcal{F}_t + \mathcal{F}_u^3 \mathcal{F} \equiv \mathcal{F}^{(4)},
\end{align}
\end{subequations}
and so on. Inserting~\eqref{TimeDerivativesOfDependentVariableODE1D} into~\eqref{ODE1DExactSolutionAtTimeLevelNPlusOne_1},
\begin{equation}
u^{n+1} = u^n + \Delta t \left(\mathcal{F}^{(1)}\right)^n + \frac{\Delta t^2}{2} \left(\mathcal{F}^{(2)}\right)^n + \frac{\Delta t^3}{6} \left(\mathcal{F}^{(3)}\right)^n + \frac{\Delta t^4}{24} \left(\mathcal{F}^{(4)}\right)^n + \cdots = u^n + \sum \limits_{k=1}^{\infty} \frac{\Delta t^k}{k!} \left(\mathcal{F}^{(k)}\right)^n.
\label{ODE1DExactSolutionAtTimeLevelNPlusOne_2}
\end{equation}
We then insert~\eqref{ODE1DExactSolutionAtTimeLevelNPlusOne_2} into~\eqref{ODE1DTruncationErrorNumericalSolution_1}, expand each term in the formula for $\hat{u}^{n+1}$ using a Taylor series with a common center, and the result is the final form of $\hat{\tau}^{n+1}$. 

The numerical solution of any time-stepping method can be expressed as 
\begin{equation}
\hat{u}^{n+1} = u^n + \sum \limits_{k=1}^{\infty} \frac{\Delta t^k}{k!} \left(\widehat{\mathcal{F}}^{(k)}\right)^n, \label{ODE1DNumericalSolutionAtTimeLevelNPlusOne}
\end{equation}
as we will illustrate with a number of examples. Here $\left(\widehat{\mathcal{F}}^{(k)}\right)^n$ is the discrete equivalent of $\left(\mathcal{F}^{(k)}\right)^n$ for $k = 1,2,\ldots$, and is defined by the time-stepping method. More specifically, $\left(\widehat{\mathcal{F}}^{(k)}\right)^n$ represents the coefficient of $\frac{\Delta t^k}{k!}$ in the numerical solution $\hat{u}^{n+1}$. For a time-stepping method of order $\beta$, $\left(\widehat{\mathcal{F}}^{(k)}\right)^n = \left(\mathcal{F}^{(k)}\right)^n$ for $k = 1,2,\ldots,\beta$. Inserting~\eqref{ODE1DExactSolutionAtTimeLevelNPlusOne_2} and~\eqref{ODE1DNumericalSolutionAtTimeLevelNPlusOne} into~\eqref{ODE1DTruncationErrorNumericalSolution_1}, 
\begin{equation}
\hat{\tau}^{n+1} = \sum \limits_{k=\beta+1}^{\infty} \frac{\Delta t^k}{k!} \left\{\left(\mathcal{F}^{(k)}\right)^n - \left(\widehat{\mathcal{F}}^{(k)}\right)^n\right\} = \frac{c_{\beta+1}^n}{(\beta+1)!} \Delta t^{\beta+1} + \bigO\left(\Delta t^{\beta+2}\right) = \bigO\left(\Delta t^{\beta+1}\right), \label{ODE1DTruncationErrorNumericalSolution_2}
\end{equation}
where $c_{\beta+1} = \mathcal{F}^{(\beta+1)} - \widehat{\mathcal{F}}^{(\beta+1)} \ne 0$. At any time horizon $T = N \Delta t$, the global truncation error of the numerical solution results from the accumulation of these local truncation errors over $N = T/\Delta t$ time steps, and is one order of $\Delta t$ less than the local truncation error. 

As pointed out by \citet{leveque2002finite}, the order of accuracy of a numerical method is not the only important attribute worth considering. The magnitude of the numerical error also depends on the coefficients of the leading order terms of the truncation error, which in turn depends on the problem being solved, the spatial and temporal discretizations, and the time horizon at which the error is being computed. If this coefficient is a few orders of magnitude larger for a high-order method than for a low-order method, the latter may result in a lower magnitude of the numerical error and turn out to be the better option. Moreover, it is only in the asymptotic regime, where the discretization parameters tend to zero, that the higher order terms are negligible with respect to the leading order terms. However, in practice, one may not employ such small values of the discretization parameters, in which case the coefficients of the higher terms cannot be neglected, and the leading order terms may not reign supreme in terms of the error magnitude. Finally, in some applications like ocean modeling, numerical stability which guarantees that small errors are not amplified by the numerical method, and the conservation of physical quantities are assigned higher priority over accuracy. In his paper, however, we focus mostly on the order of accuracy and not as much on numerical stability.

We now derive the final form of the local truncation error for the following sets of explicit and implicit time-stepping methods, belonging to the Method of Lines, and demonstrate that the numerical solution after one time step assumes the form \eqref{ODE1DNumericalSolutionAtTimeLevelNPlusOne} and the local truncation error assumes the form \eqref{ODE1DTruncationErrorNumericalSolution_2}.

\begin{mylist}\mbox{Explicit time-stepping methods for local truncation error analysis:}
\begin{enumerate}[label=(\alph*),noitemsep] 
\label{myListOfExplicitTimeSteppingMethodsForAnalysis}
\item first-order Forward Euler method;
\item explicit midpoint method, belonging to the second-order Runge-Kutta family;
\item low-storage third-order Runge-Kutta method of~\citet{williamson1980low};
\item second-order Adams-Bashforth method;
\item third-order Adams-Bashforth method.
\end{enumerate}
\end{mylist}

\begin{mylist}\mbox{Implicit time-stepping methods for local truncation error analysis:}
\begin{enumerate}[label=(\alph*),noitemsep] 
\label{myListOfImplicitTimeSteppingMethodsForAnalysis}
\item first-order Backward Euler method;
\item second-order implicit midpoint method;
\item second-order trapezoidal rule (Crank-Nicolson).
\end{enumerate}
\end{mylist}

\subsection{Forward Euler Time-Stepping Method}

The first-order Forward Euler method is
\begin{equation}
\hat{u}^{n+1} = u^n + \Delta t \mathcal{F}^n,
\end{equation}
and the local truncation error is
\begin{equation}
\hat{\tau}^{n+1} = u^{n+1} - \hat{u}^{n+1} = \frac{\Delta t^2}{2!} \left(\mathcal{F}^{(2)}\right)^n
+ O \left(\Delta t^3\right) = \bigO \left(\Delta t^2\right).
\end{equation}

\subsection{Runge-Kutta Time-Stepping Methods}
We next consider two explicit Runge-Kutta methods: the second-order explicit midpoint method and the low-storage third-order method of \citet{williamson1980low}. The explicit midpoint method is
\begin{align}
\hat{u}^{n+1} = u^n + \Delta t \mathcal{F} \left(\hat{u}^{n+\frac{1}{2}}, t^{n+\frac{1}{2}}\right) &= u^n + \Delta t \mathcal{F} \left(u^n + \frac{\Delta t}{2} \mathcal{F}^n, t^n + \frac{\Delta t}{2}\right) \nonumber \\
&= u^n + \frac{\Delta t}{1!} \left(\mathcal{F}^{(1)}\right)^n + \frac{\Delta t^2}{2!} \left(\mathcal{F}^{(2)}\right)^n + \frac{\Delta t^3}{3!} \left(\widehat{\mathcal{F}}^{(3)}\right)^n + \bigO\left(\Delta t^4\right), \label{ODE1DDepedentVariableNPlusOneEMM}
\end{align}
where we have expanded $\mathcal{F} \left(u^n + \frac{\Delta t}{2} \mathcal{F}^n, t^n + \frac{\Delta t}{2}\right)$ in a Taylor series about $u^n$ and $t^n$, where
\begin{equation}
\left(\widehat{\mathcal{F}}^{(3)}\right)^n = \frac{3}{4} \left(\mathcal{F}_{uu} \mathcal{F}^2 + 2 \mathcal{F}_{ut} \mathcal{F} + \mathcal{F}_{tt}\right)^n \ne \left(\mathcal{F}^{(3)}\right)^n,
\end{equation}
and $\hat{u}^{n+\frac{1}{2}} = u^n + \frac{\Delta t}{2} \mathcal{F}^n$ is the predicted solution at time level $t^{n+\frac{1}{2}} = t^n + \frac{\Delta t}{2}$. The local truncation error is
\begin{align}
\hat{\tau}^{n+1} = u^{n+1} - \hat{u}^{n+1} = \frac{\Delta t^3}{3!} c_3^n + \bigO\left(\Delta t^4\right) = \bigO\left(\Delta t^3\right),
\end{align}
where
\begin{equation}
c_3^n = \left(\mathcal{F}^{(3)}\right)^n - \left(\widehat{\mathcal{F}}^{(3)}\right)^n = \left(\mathcal{F}_u \mathcal{F}_t + \mathcal{F}_u^2 \mathcal{F}\right)^n + \frac{1}{4} \left(\mathcal{F}_{uu} \mathcal{F}^2 + 2 \mathcal{F}_{ut} \mathcal{F} + \mathcal{F}_{tt}\right)^n \ne 0.
\end{equation}

The third-order Runge-Kutta method of~\citet{williamson1980low} is outlined in Algorithm~1. This method has the advantage that it requires only two levels of storage. The stages are \vspace{1.5mm}\\
\textbf{Stage 1}
\begin{align}
\hat{u}^{n+\frac{1}{3}} &= u^n + \frac{1}{3} \Delta t \mathcal{F} \left(u^n,t^n\right) \equiv u^n + \Delta \hat{u}^{n + \frac{1}{3}}, \\
\mathcal{F}_{\text{mean}} \left(\hat{u}^{n+\frac{1}{3}},t^{n+\frac{1}{3}}\right) &= 
-\frac{5}{9} \mathcal{F} \left(u^n,t^n\right) + \mathcal{F} \left(\hat{u}^{n+\frac{1}{3}},t^{n+\frac{1}{3}}\right) \equiv -\frac{5}{9} \mathcal{F}^n + \mathcal{F} \left(u^n + \Delta \hat{u}^{n + \frac{1}{3}}, t^n + \frac{1}{3} \Delta t\right).
\end{align}
\textbf{Stage 2}
\begin{align}
\hat{u}^{n+\frac{3}{4}} &= u^{n+\frac{1}{3}} + \frac{15}{16} \Delta t \mathcal{F}_{\text{mean}} \left(\hat{u}^{n+\frac{1}{3}},t^{n+\frac{1}{3}}\right) \equiv u^n + \Delta \hat{u}^{n + \frac{3}{4}}, \\
\mathcal{F}_{\text{mean}} \left(\hat{u}^{n+\frac{3}{4}},t^{n+\frac{3}{4}}\right) &= 
-\frac{153}{128} \mathcal{F}_{\text{mean}} \left(\hat{u}^{n+\frac{1}{3}},t^{n+\frac{1}{3}}\right) + \mathcal{F} \left(\hat{u}^{n+\frac{3}{4}},t^{n+\frac{3}{4}}\right) \equiv -\frac{153}{128} \mathcal{F}_{\text{mean}}^{n+\frac{1}{3}} + \mathcal{F} \left(u^n + \Delta \hat{u}^{n + \frac{3}{4}},t^n + \frac{3}{4} \Delta t\right).
\end{align}
\textbf{Stage 3}
\begin{align}
\hat{u}^{n+1} &= u^{n+\frac{3}{4}} + \frac{8}{15} \Delta t \mathcal{F}_{\text{mean}} \left(\hat{u}^{n+\frac{3}{4}},t^{n+\frac{3}{4}}\right) = u^n + \frac{\Delta t}{1!} \left(\mathcal{F}^{(1)}\right)^n + \frac{\Delta t^2}{2!} \left(\mathcal{F}^{(2)}\right)^n + \frac{\Delta t^3}{3!} \left(\mathcal{F}^{(3)}\right)^n + \frac{\Delta t^4}{4!} \left(\widehat{\mathcal{F}}^{(4)}\right)^n + \bigO\left(\Delta t^5\right).
\end{align}
The full expressions for $\mathcal{F}_{\text{mean}} \left(\hat{u}^{n+\theta_k},t^{n+\theta_k}\right)$ are obtained by Taylor expanding
\begin{equation}
\mathcal{F} \left(\hat{u}^{n+\theta_k},t^{n+\theta_k}\right) \equiv \mathcal{F} \left(u_j^n + \Delta \hat{u}^{n+\theta_k}, t^n+\theta_k \Delta t\right),
\end{equation}
about $u_j^n$ and $t^n$ for $k=1$, $2$ and $\theta_1 = \frac{1}{3}$, $\theta_2 = \frac{3}{4}$. The final expression of the numerical solution at time level $t^{n+1}$ is
\begin{align}
\left(\widehat{\mathcal{F}}^{(4)}\right)^n &= \frac{1}{18} \left(17 \mathcal{F}^3 \mathcal{F}_{uuu} + 66 \mathcal{F}^2 \mathcal{F}_u \mathcal{F}_{uu} + 51 \mathcal{F}^2 \mathcal{F}_{uut} + 54 \mathcal{F} \mathcal{F}_t \mathcal{F}_{uu}\right. \nonumber \\
&\hspace{0.975cm}+ 78 \mathcal{F} \mathcal{F}_u \mathcal{F}_{ut} + 51 \mathcal{F} \mathcal{F}_{utt} + 54 \mathcal{F}_t \mathcal{F}_{ut} + 12 \mathcal{F}_{tt} \mathcal{F}_u + 17 \mathcal{F}_{ttt}\Big)^n \ne \left(\mathcal{F}^{(4)}\right)^n.
\end{align}
Therefore, the local truncation error is 
\begin{align}
\hat{\tau}^{n+1} = u^{n+1} - \hat{u}^{n+1} &= \frac{\Delta t^4}{4!} c_4^n + \bigO\left(\Delta t^5\right) = \bigO\left(\Delta t^4\right),
\end{align}
where
\begin{align}
c_4^n &= \left(\mathcal{F}^{(4)}\right)^n - \left(\widehat{\mathcal{F}}^{(4)}\right)^n \nonumber \\
&= \frac{1}{18} \left(\mathcal{F}^3 \mathcal{F}_{uuu} + 6 \mathcal{F}^2 \mathcal{F}_u \mathcal{F}_{uu} + 3 \mathcal{F}^2 \mathcal{F}_{uut} + 18 \mathcal{F} \mathcal{F}_u^3 + 12 \mathcal{F} \mathcal{F}_u \mathcal{F}_{ut} + 3 \mathcal{F} \mathcal{F}_{utt} + 18 \mathcal{F}_t \mathcal{F}_u^2 + 6 \mathcal{F}_{tt} \mathcal{F}_u + \mathcal{F}_{ttt}\right)^n \ne 0.
\end{align}

Summarizing, for a predictor-corrector Runge-Kutta method of order $\beta$, $\left(\widehat{\mathcal{F}}^{(k)}\right)^n = \left(\mathcal{F}^{(k)}\right)^n$ for $k = 1,2,\ldots,\beta$, and $\left(\widehat{\mathcal{F}}^{(\beta+1)}\right)^n$ only consists of terms in $\left(\mathcal{F}^{(\beta+1)}\right)^n$, but not necessarily with the correct multiplicative factor. As a result, $\left(\widehat{\mathcal{F}}^{(\beta+1)}\right)^n \ne \left(\mathcal{F}^{(\beta+1)}\right)^n$ and the local truncation error assumes the form~\eqref{ODE1DTruncationErrorNumericalSolution_2}.

\subsection{Adams-Bashforth Time-Stepping Methods}
We now consider multistep Adams-Bashforth methods. These methods involve the solution at time levels $t^{n-m}$ for $m = 1, 2, \ldots$ that is given by
\begin{equation}
u^{n-m} = u^n - \frac{m \Delta t}{1!} \left(\mathcal{F}^{(1)}\right)^n + \frac{(m \Delta t)^2}{2!} \left(\mathcal{F}^{(2)}\right)^n - \frac{(m \Delta t)^3}{3!} \left(\mathcal{F}^{(3)}\right)^n + \bigO\left(\Delta t^4\right) \equiv u^n + \Delta u^{n-m},
\end{equation}
where $\Delta u^{n-m} = \sum \limits_{k=1}^{\infty} \frac{(-m \Delta t)^k}{k!} \left(\mathcal{F}^{(k)}\right)^n$. The second-order Adams-Bashforth method leads to the numerical solution
\begin{align}
\hat{u}^{n+1} &= u^n + \Delta t \left\{\frac{3}{2} \mathcal{F}\left(u^n,t^n\right) - \frac{1}{2} \mathcal{F}\left(u^{n-1},t^{n-1}\right)\right\} \equiv u^n + \Delta t \left\{\frac{3}{2} \mathcal{F}^n - \frac{1}{2} \mathcal{F}\left(u^n + \Delta u^{n-1}, t^n - \Delta t\right)\right\} \nonumber \\
&= u^n + \frac{\Delta t}{1!} \left(\mathcal{F}^{(1)}\right)^n + \frac{\Delta t^2}{2!} \left(\mathcal{F}^{(2)}\right)^n + \frac{\Delta t^3}{3!} \left(\widehat{\mathcal{F}}^{(3)}\right)^n + \bigO\left(\Delta t^4\right),
\end{align}
where $\left(\widehat{\mathcal{F}}^{(3)}\right)^n = -\frac{3}{2} \left(\mathcal{F}^{(3)}\right)^n \ne \left(\mathcal{F}^{(3)}\right)^n$. Therefore, the local truncation error is
\begin{equation}
\hat{\tau}^{n+1} = u^{n+1} - \hat{u}^{n+1} = \frac{\Delta t^3}{3!} c_3^n + \bigO\left(\Delta t^4\right) = \bigO\left(\Delta t^3\right),
\end{equation}
where $c_3^n = \left(\mathcal{F}^{(3)}\right)^n - \left(\widehat{\mathcal{F}}^{(3)}\right)^n = \frac{5}{2} \left(\widehat{\mathcal{F}}^{(3)}\right)^n \ne 0$.

Repeating the same calculation for the third-order Adams-Bashforth method, we obtain the numerical solution
\begin{align}
\hat{u}^{n+1} &= u^n + \Delta t \left\{\frac{23}{12} \mathcal{F}\left(u^n,t^n\right) - \frac{16}{12} \mathcal{F}\left(u^{n-1},t^{n-1}\right) + \frac{5}{12} \mathcal{F}\left(u^{n-2},t^{n-2}\right)\right\} \nonumber \\
&\equiv u^n + \Delta t \left\{\frac{23}{12} \mathcal{F}\left(u^n,t^n\right) - \frac{16}{12} \mathcal{F}\left(u^n + \Delta u^{n-1}, t - \Delta t\right) + \frac{5}{12} \mathcal{F}\left(u^n + \Delta u^{n-2}, t - 2\Delta t\right)\right\} \nonumber \\
&= u^n + \frac{\Delta t}{1!} \left(\mathcal{F}^{(1)}\right)^n + \frac{\Delta t^2}{2!} \left(\mathcal{F}^{(2)}\right)^n + \frac{\Delta t^3}{3!} \left(\mathcal{F}^{(3)}\right)^n + \frac{\Delta t^4}{4!}\left(\widehat{\mathcal{F}}^{(4)}\right)^n + \bigO\left(\Delta t^5\right),
\end{align}
and the local truncation error
\begin{equation}
\hat{\tau}^{n+1} = u^{n+1} - \hat{u}^{n+1} = \frac{\Delta t^4}{4!} c_4^n + \bigO\left(\Delta t^5\right) = \bigO\left(\Delta t^4\right),
\end{equation}
where $\left(\widehat{\mathcal{F}}^{(4)}\right)^n = -8 \left(\mathcal{F}^{(4)}\right)^n \ne \left(\mathcal{F}^{(4)}\right)^n$, and $c_4^n = \left(\mathcal{F}^{(4)}\right)^n - \left(\widehat{\mathcal{F}}^{(4)}\right)^n = 9 \left(\widehat{\mathcal{F}}^{(4)}\right)^n \ne 0$.

Summarizing, a multistep Adams-Bashforth method of order $\beta$ results in $\left(\widehat{\mathcal{F}}^{(k)}\right)^n = \left(\mathcal{F}^{(k)}\right)^n$ for $k = 1,2,\ldots,\beta$, and $\left(\widehat{\mathcal{F}}^{(\beta+1)}\right)^n = \gamma \left(\mathcal{F}^{(\beta+1)}\right)^n \ne \left(\mathcal{F}^{(\beta+1)}\right)^n$ for some $\gamma \ne 1$, which in turn produces the same form of the local truncation error as~\eqref{ODE1DTruncationErrorNumericalSolution_2}.

\subsection{Implicit Time-Stepping Methods}

We consider the three implicit time-stepping methods of List~\ref{myListOfImplicitTimeSteppingMethodsForAnalysis}. The Backward Euler method is
\begin{equation}
\hat{u}^{n+1} = u^n + \Delta t \mathcal{F} \left(u^{n+1},t^{n+1}\right) = u^n + \Delta t \mathcal{F} \left(u^n + \Delta u^{n+1},t^n + \Delta t\right) = u^n + \frac{\Delta t}{1!} \left(\mathcal{F}^{(1)}\right)^n + \frac{\Delta t^2}{2!} \left(\widehat{\mathcal{F}}^{(2)}\right)^n + \bigO\left(\Delta t^3\right),
\end{equation}
where $\Delta u^{n+1} = \sum \limits_{k=1}^{\infty} \frac{\Delta t^k}{k!} \left(\mathcal{F}^{(k)}\right)^n$ from~\eqref{ODE1DExactSolutionAtTimeLevelNPlusOne_2} and $\left(\widehat{\mathcal{F}}^{(2)}\right)^n = 2 \left(\mathcal{F}^{(2)}\right)^n \ne \left(\mathcal{F}^{(2)}\right)^n$. The local truncation error is
\begin{equation}
\hat{\tau}^{n+1} = u^{n+1} - \hat{u}^{n+1} = \frac{\Delta t^2}{2!} c_2^n + \bigO\left(\Delta t^3\right) = \bigO\left(\Delta t^2\right),
\end{equation}
where $c_2^n = \left(\mathcal{F}^{(2)}\right)^n - \left(\widehat{\mathcal{F}}^{(2)}\right)^n = -\left(\mathcal{F}^{(2)}\right)^n \ne 0$.

The second-order implicit midpoint method is
\begin{align}
\hat{u}^{n+1} &= u^n + \Delta t \mathcal{F} \left(\frac{1}{2} \left(u^n + u^{n+1}\right),t^{n+\frac{1}{2}}\right) \equiv u^n + \Delta t \mathcal{F} \left(u^n + \frac{1}{2} \Delta u^{n+1},t^n + \frac{\Delta t}{2}\right) \nonumber \\
&= u^n + \frac{\Delta t}{1!} \left(\mathcal{F}^{(1)}\right)^n + \frac{\Delta t^2}{2!} \left(\mathcal{F}^{(2)}\right)^n + \frac{\Delta t^3}{3!} \left(\widehat{\mathcal{F}}^{(3)}\right)^n + \bigO\left(\Delta t^4\right),
\end{align}
where
\begin{equation}
\left(\widehat{\mathcal{F}}^{(3)}\right)^n = \frac{3}{4} \left(\mathcal{F}_{uu} \mathcal{F}^2 + 2 \mathcal{F}_u^2 \mathcal{F} + \mathcal{F}_{ut} \mathcal{F} + 2 \mathcal{F}_u \mathcal{F}_t + \mathcal{F}_{tt}\right)^n \ne \left(\mathcal{F}^{(3)}\right)^n.
\end{equation}
The local truncation error is
\begin{equation}
\hat{\tau}^{n+1} = u^{n+1} - \hat{u}^{n+1} = \frac{\Delta t^3}{3!} c_3^n + \bigO\left(\Delta t^4\right) = \bigO\left(\Delta t^3\right),
\end{equation}
where
\begin{equation}
c_3^n = \left(\mathcal{F}^{(3)}\right)^n - \left(\widehat{\mathcal{F}}^{(3)}\right)^n = \frac{1}{24} \left(\mathcal{F}_{uu} \mathcal{F}^2 - 2 \mathcal{F}_u^2 \mathcal{F} + 2 \mathcal{F}_{ut} \mathcal{F} - 2 \mathcal{F}_u \mathcal{F}_t + \mathcal{F}_{tt}\right)^n \ne 0.
\end{equation}

Finally, the trapezoidal rule (Crank-Nicolson) is
\begin{align}
\hat{u}^{n+1} &= u^n + \frac{\Delta t}{2} \left\{\mathcal{F} \left(u^n,t^n\right) + \mathcal{F} \left(u^{n+1},t^{n+1}\right)\right\} \equiv u^n + \frac{\Delta t}{2} \left\{\mathcal{F} \left(u^n,t^n\right) + \mathcal{F} \left(u^n + \Delta u^{n+1},t^n + \Delta t\right)\right\} \nonumber \\
&= u^n + \frac{\Delta t}{1!} \left(\mathcal{F}^{(1)}\right)^n + \frac{\Delta t^2}{2!} \left(\mathcal{F}^{(2)}\right)^n + \frac{\Delta t^3}{3!} \left(\widehat{\mathcal{F}}^{(3)}\right)^n + \bigO\left(\Delta t^4\right),
\end{align}
and the local truncation error is
\begin{equation}
\hat{\tau}^{n+1} = u^{n+1} - \hat{u}^{n+1} = \frac{\Delta t^3}{3!} c_3^n + \bigO\left(\Delta t^4\right) = \bigO\left(\Delta t^3\right),
\end{equation}
where $\left(\widehat{\mathcal{F}}^{(3)}\right)^n = \frac{3}{2} \left(\mathcal{F}^{(3)}\right)^n \ne \left(\mathcal{F}^{(3)}\right)^n$, and $c_3^n = \left(\mathcal{F}^{(3)}\right)^n - \left(\widehat{\mathcal{F}}^{(3)}\right)^n = -\frac{1}{2} \left(\mathcal{F}^{(3)}\right)^n \ne 0$.

Based on our examples, we observe that an implicit time-stepping method of order $\beta$ results in $\left(\widehat{\mathcal{F}}^{(k)}\right)^n = \left(\mathcal{F}^{(k)}\right)^n$ for $k = 1,2,\ldots,\beta$. When considering a predictor-corrector method, $\left(\widehat{\mathcal{F}}^{(\beta+1)}\right)^n$ consists of all the terms as in $\left(\mathcal{F}^{(\beta+1)}\right)^n$, but mostly with different pre-factors, and if considering a multistep method, $\left(\widehat{\mathcal{F}}^{(\beta+1)}\right)^n$ is a scalar multiple of $\left(\mathcal{F}^{(\beta+1)}\right)^n$, with the scalar factor not equal to one. In all cases, though, $\left(\widehat{\mathcal{F}}^{(\beta+1)}\right)^n \ne \left(\mathcal{F}^{(\beta+1)}\right)^n$, resulting in the the local truncation error assuming the form~\eqref{ODE1DTruncationErrorNumericalSolution_2}.

\subsection{Local Truncation Errors of a First-Order Linear ODE} \label{LTE_First_Order_Linear_ODE}
Motivated by the form of the linear inhomogeneous variable-coefficient advection equation appearing in Section \ref{LTE_Linear_Inhomogeneous_Variable_Coefficient_Advection_Equation}, we consider the ODE
\begin{equation}
u_t + \left(p_0 + q_1\right) u = f(t). \label{ODE1DParticularChoice}
\end{equation}
Using the notation in~\eqref{ODE1D}, we have
\begin{equation}
\mathcal{F}(u,t) = -\left(p_0 + q_1\right) u + f(t), \label{mathcal_F_Particular_Choice}
\end{equation}
where $p_0$, $q_1$ are constants, and $f(t)$ is a function of the independent variable $t$. The last two rows of Table~B.1 in Appendix~B express the analytical second- and third-order derivatives of the dependent variable $u$ as functions of $u$, $p_0$, $q_1$, $f(t)$ and its derivatives i.e.~$f_t$, $f_{tt}$, $f_{ttt}$, $\ldots$ at time level $t^n$. Therefore, these rows of Table~B.1 are equivalent to~\eqref{TimeDerivativesOfDependentVariableODE1D} with this specific form of $\mathcal{F}(u,t)$. Table~B.2 lists the local truncation error of~\eqref{ODE1DParticularChoice} advanced with the five explicit time-stepping methods of List~\ref{myListOfExplicitTimeSteppingMethodsForAnalysis}. Unfortunately, we cannot employ implicit time-stepping methods to advance~\eqref{ODE1DParticularChoice} since it would require knowledge of the functional form of $f(t)$. We were able to do so with the generic ODE~\eqref{ODE1D}, since at any later time level $t^{n+k} = t^n + k \Delta t$, for $k > 0$, we could substitute the Taylor expansion of the exact solution $u^{n+k}$ about $u^n$ in $\mathcal{F}\left(u^{n+k},t^{n+k}\right)$, and expand it as a Taylor series about $u^n$ and $t^n$. It is noteworthy that the local truncation error of \eqref{ODE1DParticularChoice} satisfies the form~\eqref{ODE1DTruncationErrorNumericalSolution_2} for any choice of $\mathcal{F}(u,t)$. The motivation behind the particular choice \eqref{mathcal_F_Particular_Choice} of $\mathcal{F}(u,t)$ in this example will become apparent in Section \ref{LTE_Linear_Inhomogeneous_Variable_Coefficient_Advection_Equation}.

\section{Partial Differential Equations} \label{sec:pdes}
In this paper we consider first-order one-dimensional hyperbolic PDEs of the form $u_t = \mathcal{F}(u,u_x,x,t)$. For notational convenience, we replace $u_x$ with $v$ so that the generic hyperbolic PDE we investigate is
\begin{equation}
u_t = \mathcal{F}(u,v,x,t). \label{AdvectionEquation1DFunctionalForm}
\end{equation}
The exact solution of~\eqref{AdvectionEquation1DFunctionalForm} at spatial location $x_j = j \Delta x$ and time level $t^{n+1} = t^n + \Delta t$ is
\begin{equation}
u_j^{n+1} = u_j^n + \Delta t \left(u_t\right)_j^n + \frac{\Delta t^2}{2} \left(u_{tt}\right)_j^n + \frac{\Delta t^3}{6} \left(u_{ttt}\right)_j^n + \frac{\Delta t^4}{24} \left(u_{tttt}\right)_j^n + \cdots. \label{PDE1DDependentVariableAtTimeLevelNPlusOne_1}
\end{equation}
Repeatedly differentiating~\eqref{AdvectionEquation1DFunctionalForm} with respect to time, and expressing the right-hand side in terms of known quantities at the current time level, as we did in~\eqref{TimeDerivativesOfDependentVariableODE1D}, 
\begin{subequations} \label{TemporalDerivativesOfDependentVariablePDE1D}
\begin{align} 
\frac{\del^k u}{\del t^k} &= \mathcal{F}^{(k)}, \label{TemporalDerivativesOfDependentVariablePDE1D_1} \\
\frac{\del^k v}{\del t^k} \equiv \frac{\del^k}{\del t^k} \left(\frac{\del u}{\del x}\right) &= \mathcal{G}^{(k)}, \label{TemporalDerivativesOfDependentVariablePDE1D_2}
\end{align} 
\end{subequations}
for $k=1,2,3,\ldots$ with $\mathcal{F}^{(1)} = \mathcal{F}$, $\mathcal{G}^{(1)} = \mathcal{G}$. Additional expressions of $\mathcal{F}^{(k)}$ and $\mathcal{G}^{(k)}$ are in Tables~\ref{TemporalDerivativesOfDependentVariableAsFunctionsOfKnownQuantities_GenericPDE} and~\ref{TemporalDerivativesOfSpatialGradientOfDependentVariableAsFunctionsOfKnownQuantities_GenericPDE}. For a generic PDE, all derivatives of $\mathcal{F}$ need to be considered, but since we are interested in linear and non-linear advection equations, we assume that
\begin{subequations} \label{myApproximationsForPDEsUsedToModelPhysicalPhenomena}
\begin{align}
\frac{\del^{k} \mathcal{F}}{\del u^{k}} &= 0, \text{ for $k = 2, 3, \ldots$}, \label{myApproximationsForPDEsUsedToModelPhysicalPhenomena_1} \\
\frac{\del^{k} \mathcal{F}}{\del v^{k}} &= 0, \text{ for $k = 2, 3, \ldots$}, \label{myApproximationsForPDEsUsedToModelPhysicalPhenomena_2} \\
\frac{\del^{l}}{\del t^{l}} \left(\frac{\del^{k} \mathcal{F}}{\del u^{k}}\right) &= 0, \text{ for $k = 1, 2, \ldots$, $l = 1, 2, \ldots$}, \label{myApproximationsForPDEsUsedToModelPhysicalPhenomena_3} \\
\frac{\del^{l}}{\del t^{l}} \left(\frac{\del^{k} \mathcal{F}}{\del v^{k}}\right) &= 0, \text{ for $k = 1, 2, \ldots$, $l = 1, 2, \ldots$}. \label{myApproximationsForPDEsUsedToModelPhysicalPhenomena_4} 
\end{align}
\end{subequations} 
Assumptions \eqref{myApproximationsForPDEsUsedToModelPhysicalPhenomena_1} and \eqref{myApproximationsForPDEsUsedToModelPhysicalPhenomena_2} imply that $\mathcal{F}$ only consists of linear functions of $u$ or linear functions of $v$ or products of linear functions of $u$ and $v$. For example, the inviscid Burgers' equation $u_t + u u_x \equiv u_t + uv = 0$ can be expressed as $u_t = \mathcal{F} \equiv -uv$, so it satisfies assumptions \eqref{myApproximationsForPDEsUsedToModelPhysicalPhenomena_1} and \eqref{myApproximationsForPDEsUsedToModelPhysicalPhenomena_2}. Assumptions \eqref{myApproximationsForPDEsUsedToModelPhysicalPhenomena_3} and \eqref{myApproximationsForPDEsUsedToModelPhysicalPhenomena_4} eliminate the existence of terms of the form $f(t) g(u)$ and $f(t) h(v)$, where $f(t)$, $g(u)$, and $h(v)$ are non-constant functions of $t$, $u$, and $v$ respectively. Even though assumptions \eqref{myApproximationsForPDEsUsedToModelPhysicalPhenomena_1}--\eqref{myApproximationsForPDEsUsedToModelPhysicalPhenomena_4} help us reduce the extent of the leading order terms of the local truncation error from hundreds of pages to a few pages, they are not necessary to arrive at the final results. Inserting~\eqref{TemporalDerivativesOfDependentVariablePDE1D_1} into~\eqref{PDE1DDependentVariableAtTimeLevelNPlusOne_1}, we rewrite $u_j^{n+1}$ as a function of known quantities at time level $t^n$
\begin{equation}
u_j^{n+1} = u_j^n + \Delta t \left(\mathcal{F}^{(1)}\right)_j^n + \frac{\Delta t^2}{2} \left(\mathcal{F}^{(2)}\right)_j^n + \frac{\Delta t^3}{6} \left(\mathcal{F}^{(3)}\right)_j^n + \frac{\Delta t^4}{24} \left(\mathcal{F}^{(4)}\right)_j^n + \cdots = u_j^n + \sum \limits_{k=1}^{\infty} \frac{\Delta t^k}{k!} \left(\mathcal{F}^{(k)}\right)_j^n. \label{PDE1DDependentVariableAtTimeLevelNPlusOne_2}
\end{equation}

\begin{table}[!htp]
\centering
\caption{Temporal derivatives of the dependent variable of the generic hyperbolic PDE \eqref{AdvectionEquation1DFunctionalForm} with approximations~\eqref{myApproximationsForPDEsUsedToModelPhysicalPhenomena}, up to fourth order, expressed as functions of quantities known at the current time level, and $v=u_x$, $w_1 = u_{xx}$, $w_2 = u_{xxx}$, and $w_3 = u_{xxxx}$ for notational convenience.}
\vspace{3mm}
\setlength{\tabcolsep}{0.35em}
\begin{tabular}{cc} 
\toprule
{\colorbox{shade_1}
{\parbox{1.3cm}{\centering $u_t \equiv \mathcal{F}^{(1)}$}}} & {\colorbox{shade_1}
{\parbox{0.3cm}{$\mathcal{F}$}}} \vspace{1mm} \\
{\colorbox{shade_2}
{\parbox{1.375cm}{\centering $u_{tt} \equiv \mathcal{F}^{(2)}$}}} & {\colorbox{shade_2}
{\parbox{4.75cm}{\centering $\mathcal{F} \mathcal{F}_u + \mathcal{F}_t + \mathcal{F}_u \mathcal{F}_v v + \mathcal{F}_v^2 w_1 + \mathcal{F}_v \mathcal{F}_x$}}} \vspace{1mm} \\ 
{\colorbox{shade_3}
{\parbox{1.45cm}{\centering $u_{ttt} \equiv \mathcal{F}^{(3)}$}}} & {\colorbox{shade_3}
{\parbox{10cm}{$\mathcal{F} \mathcal{F}_u^2 + 2 \mathcal{F} \mathcal{F}_u \mathcal{F}_{uv} v + 3 \mathcal{F} \mathcal{F}_{uv} \mathcal{F}_v w_1 + 2 \mathcal{F} \mathcal{F}_{uv} \mathcal{F}_x + \mathcal{F} \mathcal{F}_{ux} \mathcal{F}_v + \mathcal{F}_t \mathcal{F}_u + \mathcal{F}_{tt} \vspace{1mm} \\ + 2 \mathcal{F}_u^2 \mathcal{F}_v v + \mathcal{F}_u \mathcal{F}_{uv} \mathcal{F}_v v^2 + 3 \mathcal{F}_u \mathcal{F}_v^2 w_1 + \mathcal{F}_u \mathcal{F}_v \mathcal{F}_{vx} v + 2 \mathcal{F}_u \mathcal{F}_v \mathcal{F}_x + 3 \mathcal{F}_{uv} \mathcal{F}_v^2 v w_1 \vspace{1mm} \\ + \mathcal{F}_{uv} \mathcal{F}_v \mathcal{F}_x v + 2 \mathcal{F}_{ux} \mathcal{F}_v^2 v + \mathcal{F}_v^3 w_2 + 3 \mathcal{F}_v^2 \mathcal{F}_{vx} w_1 + \mathcal{F}_v^2 \mathcal{F}_{xx} + \mathcal{F}_v \mathcal{F}_{vx} \mathcal{F}_x + \mathcal{F}_v \mathcal{F}_{xt}$}}} \vspace{1mm} \\
{\colorbox{shade_4}
{\parbox{1.525cm}{\centering $u_{tttt} \equiv \mathcal{F}^{(4)}$}}} & {\colorbox{shade_4}
{\parbox{13.75cm}{$3 \mathcal{F}^2 \mathcal{F}_{uv}^2 w_1 + 3 \mathcal{F}^2 \mathcal{F}_{uv} \mathcal{F}_{ux} + \mathcal{F} \mathcal{F}_u^3 + 8 \mathcal{F} \mathcal{F}_u^2 \mathcal{F}_{uv} v + 3 \mathcal{F} \mathcal{F}_u \mathcal{F}_{uv}^2 v^2 + 16 \mathcal{F} \mathcal{F}_u \mathcal{F}_{uv} \mathcal{F}_v w_1 + 3 \mathcal{F} \mathcal{F}_u \mathcal{F}_{uv} \mathcal{F}_{vx} v \vspace{1mm} \\ + 8 \mathcal{F} \mathcal{F}_u \mathcal{F}_{uv} \mathcal{F}_x + 2 \mathcal{F} \mathcal{F}_u \mathcal{F}_{uvx} \mathcal{F}_v v + 3 \mathcal{F} \mathcal{F}_u \mathcal{F}_{ux} \mathcal{F}_v + 14 \mathcal{F} \mathcal{F}_{uv}^2 \mathcal{F}_v v w_1 + 3 \mathcal{F} \mathcal{F}_{uv}^2 \mathcal{F}_x v + 11 \mathcal{F} \mathcal{F}_{uv} \mathcal{F}_{ux} \mathcal{F}_v v \vspace{1mm} \\ + 6 \mathcal{F} \mathcal{F}_{uv} \mathcal{F}_v^2 w_2 + 14 \mathcal{F} \mathcal{F}_{uv} \mathcal{F}_v \mathcal{F}_{vx} w_1 + 5 \mathcal{F} \mathcal{F}_{uv} \mathcal{F}_v \mathcal{F}_{xx} + 3 \mathcal{F} \mathcal{F}_{uv} \mathcal{F}_{vx} \mathcal{F}_x + 3 \mathcal{F} \mathcal{F}_{uv} \mathcal{F}_{xt} + 4 \mathcal{F} \mathcal{F}_{uvx} \mathcal{F}_v^2 w_1 \vspace{1mm} \\ + 2 \mathcal{F} \mathcal{F}_{uvx} \mathcal{F}_v \mathcal{F}_x + \mathcal{F} \mathcal{F}_{ux} \mathcal{F}_v \mathcal{F}_{vx} + \mathcal{F} \mathcal{F}_{uxx} \mathcal{F}_v^2 + \mathcal{F}_t \mathcal{F}_u^2 + 3 \mathcal{F}_t \mathcal{F}_u \mathcal{F}_{uv} v + 4 \mathcal{F}_t \mathcal{F}_{uv} \mathcal{F}_v w_1 + 3 \mathcal{F}_t \mathcal{F}_{uv} \mathcal{F}_x \vspace{1mm} \\ + \mathcal{F}_t \mathcal{F}_{ux} \mathcal{F}_v + \mathcal{F}_{tt} \mathcal{F}_u + \mathcal{F}_{ttt} + 3 \mathcal{F}_u^3 \mathcal{F}_v v + 8 \mathcal{F}_u^2 \mathcal{F}_{uv} \mathcal{F}_v v^2 + 6 \mathcal{F}_u^2 \mathcal{F}_v^2 w_1 + 3 \mathcal{F}_u^2 \mathcal{F}_v \mathcal{F}_{vx} v + 3 \mathcal{F}_u^2 \mathcal{F}_v \mathcal{F}_x \vspace{1mm} \\ + \mathcal{F}_u \mathcal{F}_{uv}^2 \mathcal{F}_v v^3 + 26 \mathcal{F}_u \mathcal{F}_{uv} \mathcal{F}_v^2 v w_1 + 2 \mathcal{F}_u \mathcal{F}_{uv} \mathcal{F}_v \mathcal{F}_{vx} v^2 + 13 \mathcal{F}_u \mathcal{F}_{uv} \mathcal{F}_v \mathcal{F}_x v + 2 \mathcal{F}_u \mathcal{F}_{uvx} \mathcal{F}_v^2 v^2 + 9 \mathcal{F}_u \mathcal{F}_{ux} \mathcal{F}_v^2 v \vspace{1mm} \\ + 4 \mathcal{F}_u \mathcal{F}_v^3 w_2 + 12 \mathcal{F}_u \mathcal{F}_v^2 \mathcal{F}_{vx} w_1 + \mathcal{F}_u \mathcal{F}_v^2 \mathcal{F}_{vxx} v + 3 \mathcal{F}_u \mathcal{F}_v^2 \mathcal{F}_{xx} + \mathcal{F}_u \mathcal{F}_v \mathcal{F}_{vx}^2 v + 3 \mathcal{F}_u \mathcal{F}_v \mathcal{F}_{vx} \mathcal{F}_x + 2 \mathcal{F}_u \mathcal{F}_v \mathcal{F}_{xt} \vspace{1mm} \\ + 7 \mathcal{F}_{uv}^2 \mathcal{F}_v^2 v^2 w_1 + \mathcal{F}_{uv}^2 \mathcal{F}_v \mathcal{F}_x v^2 + 6 \mathcal{F}_{uv} \mathcal{F}_{ux} \mathcal{F}_v^2 v^2 + 6 \mathcal{F}_{uv} \mathcal{F}_v^3 v w_2 + 12 \mathcal{F}_{uv} \mathcal{F}_v^3 w_1^2 + 14 \mathcal{F}_{uv} \mathcal{F}_v^2 \mathcal{F}_{vx} v w_1 \vspace{1mm} \\ + 14 \mathcal{F}_{uv} \mathcal{F}_v^2 \mathcal{F}_x w_1 + 3 \mathcal{F}_{uv} \mathcal{F}_v^2 \mathcal{F}_{xx} v + 2 \mathcal{F}_{uv} \mathcal{F}_v \mathcal{F}_{vx} \mathcal{F}_x v + 5 \mathcal{F}_{uv} \mathcal{F}_v \mathcal{F}_x^2 + \mathcal{F}_{uv} \mathcal{F}_v \mathcal{F}_{xt} v + 8 \mathcal{F}_{uvx} \mathcal{F}_v^3 v w_1 \vspace{1mm} \\ + 2 \mathcal{F}_{uvx} \mathcal{F}_v^2 \mathcal{F}_x v + 6 \mathcal{F}_{ux} \mathcal{F}_v^3 w_1 + 6 \mathcal{F}_{ux} \mathcal{F}_v^2 \mathcal{F}_{vx} v + 3 \mathcal{F}_{ux} \mathcal{F}_v^2 \mathcal{F}_x + 3 \mathcal{F}_{uxx} \mathcal{F}_v^3 v + \mathcal{F}_v^4 w_3 + 6 \mathcal{F}_v^3 \mathcal{F}_{vx} w_2 \vspace{1mm} \\ + 4 \mathcal{F}_v^3 \mathcal{F}_{vxx} w_1 + \mathcal{F}_v^3 \mathcal{F}_{xxx} + 7 \mathcal{F}_v^2 \mathcal{F}_{vx}^2 w_1 + 3 \mathcal{F}_v^2 \mathcal{F}_{vx} \mathcal{F}_{xx} + \mathcal{F}_v^2 \mathcal{F}_{vxx} \mathcal{F}_x + \mathcal{F}_v^2 \mathcal{F}_{xxt} + \mathcal{F}_v \mathcal{F}_{vx}^2 \mathcal{F}_x + \mathcal{F}_v \mathcal{F}_{vx} \mathcal{F}_{xt} + \mathcal{F}_v \mathcal{F}_{xtt}$}}} \vspace{1mm} \\
\bottomrule \\
\end{tabular} \label{TemporalDerivativesOfDependentVariableAsFunctionsOfKnownQuantities_GenericPDE}
\end{table}

\begin{table}[!htp]
\centering
\caption{Temporal derivatives of the spatial gradient of the dependent variable of the generic hyperbolic PDE \eqref{AdvectionEquation1DFunctionalForm} with approximations~\eqref{myApproximationsForPDEsUsedToModelPhysicalPhenomena}, up to third order, expressed as functions of quantities known at the current time level, and $v=u_x$, $w_1 = u_{xx}$, $w_2 = u_{xxx}$, and $w_3 = u_{xxxx}$ for notational convenience.}
\vspace{3mm}
\setlength{\tabcolsep}{0.35em}
\begin{tabular}{cc} 
\toprule
{\colorbox{shade_1}
{\parbox{1.3cm}{\centering $v_t \equiv \mathcal{G}^{(1)}$}}} & {\colorbox{shade_1}
{\parbox{2.25cm}{\centering $\mathcal{F}_u v + \mathcal{F}_v w_1 + \mathcal{F}_x$}}} \vspace{1mm} \\
{\colorbox{shade_2}
{\parbox{1.375cm}{\centering $v_{tt} \equiv \mathcal{G}^{(2)}$}}} & {\colorbox{shade_2}
{\parbox{10.5cm}{$\mathcal{F} \mathcal{F}_{uv} w_1 + \mathcal{F} \mathcal{F}_{ux} + \mathcal{F}_u^2 v + \mathcal{F}_u \mathcal{F}_{uv} v^2 + 2 \mathcal{F}_u \mathcal{F}_v w_1 + \mathcal{F}_u \mathcal{F}_{vx} v + \mathcal{F}_u \mathcal{F}_x + 3 \mathcal{F}_{uv} \mathcal{F}_v v w_1 \vspace{1mm} \\ + \mathcal{F}_{uv} \mathcal{F}_x v + 2 \mathcal{F}_{ux} \mathcal{F}_v v + \mathcal{F}_v^2 w_2 + 3 \mathcal{F}_v \mathcal{F}_{vx} w_1 + \mathcal{F}_v \mathcal{F}_{xx} + \mathcal{F}_{vx} \mathcal{F}_x + \mathcal{F}_{xt}$}}} \vspace{1mm} \\ 
{\colorbox{shade_3}
{\parbox{1.45cm}{\centering $v_{ttt} \equiv \mathcal{G}^{(3)}$}}} & {\colorbox{shade_3}
{\parbox{13.75cm}{$4 \mathcal{F} \mathcal{F}_u \mathcal{F}_{uv} w_1 + 2 \mathcal{F} \mathcal{F}_u \mathcal{F}_{uvx} v + 2 \mathcal{F} \mathcal{F}_u \mathcal{F}_{ux} + 5 \mathcal{F} \mathcal{F}_{uv}^2 v w_1 + 5 \mathcal{F} \mathcal{F}_{uv} \mathcal{F}_{ux} v + 3 \mathcal{F} \mathcal{F}_{uv} \mathcal{F}_v w_2 + 5 \mathcal{F} \mathcal{F}_{uv} \mathcal{F}_{vx} w_1 \vspace{1mm} \\ + 2 \mathcal{F} \mathcal{F}_{uv} \mathcal{F}_{xx} + 4 \mathcal{F} \mathcal{F}_{uvx} \mathcal{F}_v w_1 + 2 \mathcal{F} \mathcal{F}_{uvx} \mathcal{F}_x + \mathcal{F} \mathcal{F}_{ux} \mathcal{F}_{vx} + \mathcal{F} \mathcal{F}_{uxx} \mathcal{F}_v + \mathcal{F}_t \mathcal{F}_{uv} w_1 + \mathcal{F}_t \mathcal{F}_{ux} + \mathcal{F}_u^3 v \vspace{1mm} \\ + 4 \mathcal{F}_u^2 \mathcal{F}_{uv} v^2 + 3 \mathcal{F}_u^2 \mathcal{F}_v w_1 + 2 \mathcal{F}_u^2 \mathcal{F}_{vx} v + \mathcal{F}_u^2 \mathcal{F}_x + \mathcal{F}_u \mathcal{F}_{uv}^2 v^3 + 17 \mathcal{F}_u \mathcal{F}_{uv} \mathcal{F}_v v w_1 + 2 \mathcal{F}_u \mathcal{F}_{uv} \mathcal{F}_{vx} v^2 \vspace{1mm} \\ + 6 \mathcal{F}_u \mathcal{F}_{uv} \mathcal{F}_x v + 2 \mathcal{F}_u \mathcal{F}_{uvx} \mathcal{F}_v v^2 + 7 \mathcal{F}_u \mathcal{F}_{ux} \mathcal{F}_v v + 3 \mathcal{F}_u \mathcal{F}_v^2 w_2 + 9 \mathcal{F}_u \mathcal{F}_v \mathcal{F}_{vx} w_1 + \mathcal{F}_u \mathcal{F}_v \mathcal{F}_{vxx} v + 2 \mathcal{F}_u \mathcal{F}_v \mathcal{F}_{xx} \vspace{1mm} \\ + \mathcal{F}_u \mathcal{F}_{vx}^2 v + 2 \mathcal{F}_u \mathcal{F}_{vx} \mathcal{F}_x + \mathcal{F}_u \mathcal{F}_{xt} + 7 \mathcal{F}_{uv}^2 \mathcal{F}_v v^2 w_1 + \mathcal{F}_{uv}^2 \mathcal{F}_x v^2 + 6 \mathcal{F}_{uv} \mathcal{F}_{ux} \mathcal{F}_v v^2 + 6 \mathcal{F}_{uv} \mathcal{F}_v^2 v w_2 + 9 \mathcal{F}_{uv} \mathcal{F}_v^2 w_1^2 \vspace{1mm} \\ + 14 \mathcal{F}_{uv} \mathcal{F}_v \mathcal{F}_{vx} v w_1 + 8 \mathcal{F}_{uv} \mathcal{F}_v \mathcal{F}_x w_1 + 3 \mathcal{F}_{uv} \mathcal{F}_v \mathcal{F}_{xx} v + 2 \mathcal{F}_{uv} \mathcal{F}_{vx} \mathcal{F}_x v + 2 \mathcal{F}_{uv} \mathcal{F}_x^2 + \mathcal{F}_{uv} \mathcal{F}_{xt} v + 8 \mathcal{F}_{uvx} \mathcal{F}_v^2 v w_1 \vspace{1mm} \\ + 2 \mathcal{F}_{uvx} \mathcal{F}_v \mathcal{F}_x v + 6 \mathcal{F}_{ux} \mathcal{F}_v^2 w_1 + 6 \mathcal{F}_{ux} \mathcal{F}_v \mathcal{F}_{vx} v + 3 \mathcal{F}_{ux} \mathcal{F}_v \mathcal{F}_x + 3 \mathcal{F}_{uxx} \mathcal{F}_v^2 v + \mathcal{F}_v^3 w_3 + 6 \mathcal{F}_v^2 \mathcal{F}_{vx} w_2 \vspace{1mm} \\ + 4 \mathcal{F}_v^2 \mathcal{F}_{vxx} w_1 + \mathcal{F}_v^2 \mathcal{F}_{xxx} + 7 \mathcal{F}_v \mathcal{F}_{vx}^2 w_1 + 3 \mathcal{F}_v \mathcal{F}_{vx} \mathcal{F}_{xx} + \mathcal{F}_v \mathcal{F}_{vxx} \mathcal{F}_x + \mathcal{F}_v \mathcal{F}_{xxt} + \mathcal{F}_{vx}^2 \mathcal{F}_x + \mathcal{F}_{vx} \mathcal{F}_{xt} + \mathcal{F}_{xtt}$}}} \vspace{1mm} \\
\bottomrule \\
\end{tabular} \label{TemporalDerivativesOfSpatialGradientOfDependentVariableAsFunctionsOfKnownQuantities_GenericPDE}
\end{table}

We now derive the local truncation error of the generic hyperbolic PDE~\eqref{AdvectionEquation1DFunctionalForm} discretized in space using a finite difference method, and advanced in time with the time-stepping methods in Lists~\ref{myListOfExplicitTimeSteppingMethodsForAnalysis} and~\ref{myListOfImplicitTimeSteppingMethodsForAnalysis}. The final form of this local truncation error, however, remains the same for any hyperbolic PDE, employing any type of spatial discretization, including finite element and finite volume methods, and any explicit or implicit time-stepping method, including predictor-corrector and multistep.

\subsection{Forward Euler Time-Stepping Method}

With the first-order Forward Euler time-stepping method, the numerical solution of the advection equation~\eqref{AdvectionEquation1DFunctionalForm} at spatial location $x_j$ and time level $t^{n+1}$ is 
\begin{equation}
\hat{u}_j^{n+1} = u_j^n + \Delta t \widetilde{\mathcal{F}} \left(u^n_{j-j_{\text{lower}}}, u^n_{j-j_{\text{lower}}+1}, \ldots, u^n_j, \ldots, u^n_{j-j_{\text{lower}}+\alpha},x_{j-j_{\text{lower}}}, x_{j-j_{\text{lower}}+1}, \ldots, x_j, \ldots, x_{j-j_{\text{lower}}+\alpha},t^n\right), \label{PDE1DSolutionAtTimelevelNPlusOneUsingForwardEuler_1}
\end{equation}
where $\widetilde{\mathcal{F}}$ is a spatially discretized version of $\mathcal{F}$ with $\alpha$ being the order of the spatial discretization, and the index $j_{\text{lower}}$ depends on the finite difference scheme. If we apply the first-order upwind finite difference scheme to an advection problem with positive advection velocity, then $\alpha = 1$, $j_{\text{lower}} = 1$, and~\eqref{PDE1DSolutionAtTimelevelNPlusOneUsingForwardEuler_1} reduces to 
\begin{equation}
\hat{u}_j^{n+1} = u_j^n + \Delta t \widetilde{\mathcal{F}} \left(u_{j-1}^n,u_j^n,x_{j-1},x_j,t^n\right). \label{PDE1DSolutionAtTimelevelNPlusOneUsingForwardEuler_1_FirstOrderFiniteVolumeUpwind}
\end{equation}
Using the definition of $\widetilde{\mathcal{F}}$ corresponding to Forward Euler, we rewrite~\eqref{PDE1DSolutionAtTimelevelNPlusOneUsingForwardEuler_1} as
\begin{equation}
\hat{u}_j^{n+1} = u_j^n + \Delta t \mathcal{F} \left(u_j^n, v_j^n + \bigO\left(\Delta x^{\alpha}\right), x_j, t^n\right) = u_j^n + \Delta t \left\{\mathcal{F} \left(u_j^n, v_j^n,x_j, t^n\right) + \bigO\left(\Delta x^{\alpha}\right)\right\} = u_j^n + \Delta t \left(\mathcal{F}^{(1)} + \bigO\left(\Delta x^{\alpha}\right)\right)_j^n, \label{PDE1DSolutionAtTimelevelNPlusOneUsingForwardEuler_2}
\end{equation}
and the local truncation error is
\begin{equation}
\hat{\tau}_j^{n+1} = u_j^{n+1} - \hat{u}_j^{n+1} = \Delta t \bigO\left(\Delta x^{\alpha}\right) + \frac{\Delta t^2}{2} \left(\mathcal{F}^{(2)}\right)_j^n + \bigO\left(\Delta t^3\right) = \Delta t \bigO\left(\Delta x^{\alpha}\right) + \bigO\left(\Delta t^2\right).
\end{equation}

\subsection{Runge-Kutta Time-Stepping Methods}
We next derive the local truncation error of the numerical solution resulting from the explicit second-order midpoint method and Williamson's low-storage third-order Runge-Kutta~\citep{williamson1980low}. The explicit midpoint method is
\begin{equation}
\hat{u}_j^{n+1} = u_j^n + \Delta t \mathcal{F} \left(\hat{u}_j^{n+\frac{1}{2}},\hat{v}_j^{n+\frac{1}{2}}+O\left(\Delta x^{\alpha}\right),x_j,t^{n+\frac{1}{2}}\right), \label{PDE1DSolutionAtTimelevelNPlusOneUsingEMM_1}
\end{equation}
where
\begin{equation}
\hat{u}_j^{n+\frac{1}{2}} = u_j^n + \frac{\Delta t}{2} \mathcal{F} \left(u_j^n,v_j^n + \bigO\left(\Delta x^{\alpha}\right),x_j,t^n\right) = u_j^n + \frac{\Delta t}{2} \left(\mathcal{F} + \bigO\left(\Delta x^{\alpha}\right)\right)_j^n \equiv u_j^n + \Delta \hat{u}^{n+\frac{1}{2}}_j. \label{PDE1DSolutionAtTimelevelNPlusHalfUsingEMM}
\end{equation}
Equation~\eqref{PDE1DSolutionAtTimelevelNPlusHalfUsingEMM} is the predicted value of $u$ at spatial location $x_j$ and time level $t^{n+\frac{1}{2}} = t^n + \frac{\Delta t}{2}$. Its exact spatial derivative is
\begin{equation}
\hat{v}_j^{n+\frac{1}{2}} = v_j^n + \frac{\Delta t}{2} \left(\mathcal{F}_x + \mathcal{F}_u v + \mathcal{F}_v w_1 + \bigO\left(\Delta x^{\alpha}\right)\right)_j^n \equiv v_j^n + \Delta \hat{v}_j^{n+\frac{1}{2}}, \label{PDE1DSolutionSpatialDerivativeAtTimelevelNPlusHalfUsingEMM}
\end{equation}
where $w_1 = v_x = u_{xx}$. Inserting~\eqref{PDE1DSolutionAtTimelevelNPlusHalfUsingEMM} and~\eqref{PDE1DSolutionSpatialDerivativeAtTimelevelNPlusHalfUsingEMM} into~\eqref{PDE1DSolutionAtTimelevelNPlusOneUsingEMM_1},
\begin{align}
\hat{u}_j^{n+1} &= u_j^n + \Delta t \mathcal{F} \left(u_j^n + \Delta \hat{u}^{n+\frac{1}{2}}_k,v_j^n + \Delta \hat{v}^{n+\frac{1}{2}}_k + \bigO\left(\Delta x^{\alpha}\right),x_j,t^n + \frac{\Delta t}{2}\right) \nonumber \\
&= u_j^n + \frac{\Delta t}{1!} \left(\mathcal{F}^{(1)} + \bigO\left(\Delta x^{\alpha}\right)\right)_j^n + \frac{\Delta t^2}{2!} \left(\mathcal{F}^{(2)} + \bigO\left(\Delta x^{\alpha}\right)\right)_j^n + \frac{\Delta t^3}{3!} \left(\widehat{\mathcal{F}}^{(3)} + \bigO\left(\Delta x^{\alpha}\right)\right)_j^n + \bigO\left(\Delta t^4\right),
\end{align}
where $\left(\widehat{\mathcal{F}}^{(3)}\right)_j^n \ne \left(\mathcal{F}^{(3)}\right)_j^n$. The local truncation error is
\begin{align}
\hat{\tau}_j^{n+1} &= u_j^{n+1} - \hat{u}_j^{n+1} \nonumber \\
&= \frac{\Delta t}{1!} \bigO\left(\Delta x^{\alpha}\right) + \frac{\Delta t^2}{2!} \bigO\left(\Delta x^{\alpha}\right) + \frac{\Delta t^3}{3!} \left(c_3 + \bigO\left(\Delta x^{\alpha}\right)\right)_j^n + \bigO\left(\Delta t^4\right) = \Delta t \bigO\left(\Delta x^{\alpha}\right) + \Delta t^2 \bigO\left(\Delta x^{\alpha}\right) + \bigO\left(\Delta t^3\right),
\end{align}
where $\left(c_3\right)_j^n = \left(\mathcal{F}^{(3)}\right)_j^n - \left(\mathcal{F}^{\hat{(3)}}\right)_j^n \ne 0$. The full expressions for $\frac{1}{3!} \widehat{\mathcal{F}}^{(3)}$ and $\frac{1}{3!} c_3$ are in Table~\ref{FHat3c3ExplicitMidpointMethod}.

The derivation of the local truncation error of the numerical solution resulting from Williamson's low-storage third-order Runge-Kutta time-stepping method~\citep{williamson1980low} has the following stages: \vspace{1.5mm} \\
\textbf{Stage 1}
\begin{align}
\hat{u}_j^{n+\frac{1}{3}} &= u_j^n + \frac{\Delta t}{3} \mathcal{F} \left(u_j^n,v_j^n + \bigO\left(\Delta x^{\alpha}\right),x_j,t^n\right) = u_j^n + \frac{\Delta t}{3} \left(\mathcal{F} + \bigO\left(\Delta x^{\alpha}\right)\right)_j^n \equiv u_j^n + \Delta \hat{u}_j^{n + \frac{1}{3}}, \\
\hat{v}_j^{n+\frac{1}{3}} &= v_j^n + \frac{\Delta t}{3} \left(\mathcal{F}_x + \mathcal{F}_u v + \mathcal{F}_v w_1 + \bigO\left(\Delta x^{\alpha}\right)\right)_j^n \equiv v_j^n + \Delta \hat{v}_j^{n + \frac{1}{3}}, \\
\widetilde{\mathcal{F}}_{\text{mean}} \left(\hat{u}_j^{n+\frac{1}{3}},\hat{v}_j^{n+\frac{1}{3}},x_j,t^{n+\frac{1}{3}}\right) &= 
-\frac{5}{9} \left(\mathcal{F} + \bigO\left(\Delta x^{\alpha}\right)\right)_j^n + \mathcal{F} \left(\hat{u}_j^{n+\frac{1}{3}},\hat{v}_j^{n+\frac{1}{3}} + \bigO\left(\Delta x^{\alpha}\right),x_j,t^{n+\frac{1}{3}}\right) \nonumber \\
&= -\frac{5}{9} \left(\mathcal{F} + \bigO\left(\Delta x^{\alpha}\right)\right)_j^n + \mathcal{F} \left(u_j^n + \Delta \hat{u}_j^{n+\frac{1}{3}}, v_j^n + \Delta \hat{v}_j^{n+\frac{1}{3}} + \bigO\left(\Delta x^{\alpha}\right), x_j, t^n+\frac{\Delta t}{3}\right).
\end{align}
\textbf{Stage 2}
\begin{align}
\hat{u}_j^{n+\frac{3}{4}} &= \hat{u}_j^{n+\frac{1}{3}} + \frac{15}{16} \Delta t \left(\widetilde{\mathcal{F}}_{\text{mean}}\right)_j^{n+\frac{1}{3}} \equiv u_j^n + \Delta \hat{u}_j^{n + \frac{3}{4}}, \\
\hat{v}_j^{n+\frac{3}{4}} &= \hat{v}_j^{n+\frac{1}{3}} + \frac{15}{16} \Delta t \left(\widetilde{\mathcal{F}}_{\text{mean},x} + \widetilde{\mathcal{F}}_{\text{mean},u} v + \widetilde{\mathcal{F}}_{\text{mean},v} w_1\right)_j^{n+\frac{1}{3}} \equiv v_j^n + \Delta \hat{v}_j^{n + \frac{3}{4}}, \\
\widetilde{\mathcal{F}}_{\text{mean}} \left(\hat{u}_j^{n+\frac{3}{4}},\hat{v}_j^{n+\frac{3}{4}},x_j,t^{n+\frac{3}{4}}\right) &= 
-\frac{153}{128} \left(\widetilde{\mathcal{F}}_{\text{mean}}\right)_j^{n+\frac{1}{3}} + \mathcal{F} \left(\hat{u}_j^{n+\frac{3}{4}},\hat{v}_j^{n+\frac{3}{4}} + \bigO\left(\Delta x^{\alpha}\right),x_j,t^{n+\frac{3}{4}}\right) \nonumber \\
&\equiv -\frac{153}{128} \left(\widetilde{\mathcal{F}}_{\text{mean}}\right)_j^{n+\frac{1}{3}} + \mathcal{F} \left(u_j^n + \Delta \hat{u}_j^{n+\frac{3}{4}}, v_j^n + \Delta \hat{v}_j^{n+\frac{3}{4}} + \bigO\left(\Delta x^{\alpha}\right), x_j, t^n+\frac{3}{4} \Delta t\right).
\end{align}

\begin{table}[!htp]
\centering
\caption{The term $\frac{1}{3!} \widehat{\mathcal{F}}^{(3)}$ in the numerical solution $\hat{u}_j^{n+1}$ given by \eqref{Exact_Solution_Time_Level_nP1_Theorem}, and the term $\frac{1}{3!} c_3$ in the local truncation error $\hat{\tau}_j^{n+1}$ given by \eqref{LocalTruncationErrorNumericalSolutionFinalForm_1}, of the generic hyperbolic PDE~\eqref{AdvectionEquation1DFunctionalForm} with approximations~\eqref{myApproximationsForPDEsUsedToModelPhysicalPhenomena}, advanced in time with the explicit midpoint method, at spatial location~$x_j$ and time level $t^{n+1} = t^n + \Delta t$, expressed as functions of quantities known at the current time level~$t^n$, and $v=u_x$, $w_1 = u_{xx}$, and $w_2 = u_{xxx}$ for notational convenience.}
\vspace{3mm}
\setlength{\tabcolsep}{0.35em}
\begin{tabular}{cc} 
\toprule
{\colorbox{shade_1}
{\parbox{0.85cm}{\centering $\frac{1}{3!} \widehat{\mathcal{F}}^{(3)}$}}} & {\colorbox{shade_1}
{\parbox{6.125cm}{$\frac{1}{4} \mathcal{F} \mathcal{F}_u \mathcal{F}_{uv} v + \frac{1}{4} \mathcal{F} \mathcal{F}_{uv} \mathcal{F}_v w_1 + \frac{1}{4} \mathcal{F} \mathcal{F}_{uv} \mathcal{F}_x + \frac{1}{8} \mathcal{F}_{tt}$}}} \vspace{1mm} \\ 
{\colorbox{shade_2}
{\parbox{0.55cm}{\centering $\frac{1}{3!} c_3$}}} & {\colorbox{shade_2}
{\parbox{11.25cm}{$\frac{1}{6} \mathcal{F} \mathcal{F}_u^2 + \frac{1}{12} \mathcal{F} \mathcal{F}_u \mathcal{F}_{uv} v + \frac{1}{4} \mathcal{F} \mathcal{F}_{uv} \mathcal{F}_v w_1 + \frac{1}{12} \mathcal{F} \mathcal{F}_{uv} \mathcal{F}_x + \frac{1}{6} \mathcal{F} \mathcal{F}_{ux} \mathcal{F}_v + \frac{1}{6} \mathcal{F}_t \mathcal{F}_u + \frac{1}{24} \mathcal{F}_{tt} \vspace{1mm} \\ + \frac{1}{3} \mathcal{F}_u^2 \mathcal{F}_v v + \frac{1}{6} \mathcal{F}_u \mathcal{F}_{uv} \mathcal{F}_v v^2 + \frac{1}{2} \mathcal{F}_u \mathcal{F}_v^2 w_1 + \frac{1}{6} \mathcal{F}_u \mathcal{F}_v \mathcal{F}_{vx} v + \frac{1}{3} \mathcal{F}_u \mathcal{F}_v \mathcal{F}_x + \frac{1}{2} \mathcal{F}_{uv} \mathcal{F}_v^2 v w_1 \vspace{1mm} \\ + \frac{1}{6} \mathcal{F}_{uv} \mathcal{F}_v \mathcal{F}_x v + \frac{1}{3} \mathcal{F}_{ux} \mathcal{F}_v^2 v + \frac{1}{6} \mathcal{F}_v^3 w_2 + \frac{1}{2} \mathcal{F}_v^2 \mathcal{F}_{vx} w_1 + \frac{1}{6} \mathcal{F}_v^2 \mathcal{F}_{xx} + \frac{1}{6} \mathcal{F}_v \mathcal{F}_{vx} \mathcal{F}_x + \frac{1}{6} \mathcal{F}_v \mathcal{F}_{xt}$}}} \vspace{1mm} \\
\bottomrule \\
\end{tabular} \label{FHat3c3ExplicitMidpointMethod}
\end{table}

\begin{table}[!htp]
\centering
\caption{The term $\frac{1}{4!} \widehat{\mathcal{F}}^{(4)}$ in the numerical solution $\hat{u}_j^{n+1}$ given by \eqref{Exact_Solution_Time_Level_nP1_Theorem}, and the term $\frac{1}{4!} c_4$ in the local truncation error $\hat{\tau}_j^{n+1}$ given by \eqref{LocalTruncationErrorNumericalSolutionFinalForm_1}, of the generic hyperbolic PDE~\eqref{AdvectionEquation1DFunctionalForm} with approximations~\eqref{myApproximationsForPDEsUsedToModelPhysicalPhenomena}, advanced in time with the low-storage third-order Runge-Kutta method of \citet{williamson1980low}, at spatial location~$x_j$ and time level $t^{n+1} = t^n + \Delta t$, expressed as functions of quantities known at the current time level~$t^n$, and $v=u_x$, $w_1 = u_{xx}$, and $w_2 = u_{xxx}$ for notational convenience.}
\vspace{3mm}
\setlength{\tabcolsep}{0.35em}
\begin{tabular}{cc} 
\toprule
{\colorbox{shade_1}
{\parbox{0.85cm}{\centering $\frac{1}{4!} \widehat{\mathcal{F}}^{(4)}$}}} & {\colorbox{shade_1}
{\parbox{13.75cm}{$\frac{1}{8} \mathcal{F}^2 \mathcal{F}_{uv}^2 w_1 + \frac{1}{8} \mathcal{F}^2 \mathcal{F}_{uv} \mathcal{F}_{ux} + \frac{11}{36} \mathcal{F} \mathcal{F}_u^2 \mathcal{F}_{uv} v + \frac{1}{8} \mathcal{F} \mathcal{F}_u \mathcal{F}_{uv}^2 v^2 + \frac{35}{72} \mathcal{F} \mathcal{F}_u \mathcal{F}_{uv} \mathcal{F}_v w_1 + \frac{1}{8} \mathcal{F} \mathcal{F}_u \mathcal{F}_{uv} \mathcal{F}_{vx} v \vspace{1mm} \\ + \frac{11}{36} \mathcal{F} \mathcal{F}_u \mathcal{F}_{uv} \mathcal{F}_x + \frac{1}{18} \mathcal{F} \mathcal{F}_u \mathcal{F}_{uvx} \mathcal{F}_v v + \frac{35}{72} \mathcal{F} \mathcal{F}_{uv}^2 \mathcal{F}_v v w_1 + \frac{1}{8} \mathcal{F} \mathcal{F}_{uv}^2 \mathcal{F}_x v + \frac{13}{36} \mathcal{F} \mathcal{F}_{uv} \mathcal{F}_{ux} \mathcal{F}_v v \vspace{1mm} \\ + \frac{13}{72} \mathcal{F} \mathcal{F}_{uv} \mathcal{F}_v^2 w_2 + \frac{35}{72} \mathcal{F} \mathcal{F}_{uv} \mathcal{F}_v \mathcal{F}_{vx} w_1 + \frac{13}{72} \mathcal{F} \mathcal{F}_{uv} \mathcal{F}_v \mathcal{F}_{xx} + \frac{1}{8} \mathcal{F} \mathcal{F}_{uv} \mathcal{F}_{vx} \mathcal{F}_x + \frac{1}{8} \mathcal{F} \mathcal{F}_{uv} \mathcal{F}_{xt} + \frac{1}{18} \mathcal{F} \mathcal{F}_{uvx} \mathcal{F}_v^2 w_1 \vspace{1mm} \\ + \frac{1}{18} \mathcal{F} \mathcal{F}_{uvx} \mathcal{F}_v \mathcal{F}_x + \frac{1}{8} \mathcal{F}_t \mathcal{F}_u \mathcal{F}_{uv} v + \frac{1}{8} \mathcal{F}_t \mathcal{F}_{uv} \mathcal{F}_v w_1 + \frac{1}{8} \mathcal{F}_t \mathcal{F}_{uv} \mathcal{F}_x + \frac{1}{36} \mathcal{F}_{tt} \mathcal{F}_u + \frac{17}{432} \mathcal{F}_{ttt} + \frac{13}{72} \mathcal{F}_u^2 \mathcal{F}_{uv} \mathcal{F}_v v^2 \vspace{1mm} \\ + \frac{13}{36} \mathcal{F}_u \mathcal{F}_{uv} \mathcal{F}_v^2 v w_1 + \frac{13}{36} \mathcal{F}_u \mathcal{F}_{uv} \mathcal{F}_v \mathcal{F}_x v + \frac{13}{72} \mathcal{F}_{uv} \mathcal{F}_v^3 w_1^2 + \frac{13}{36} \mathcal{F}_{uv} \mathcal{F}_v^2 \mathcal{F}_x w_1 + \frac{13}{72} \mathcal{F}_{uv} \mathcal{F}_v \mathcal{F}_x^2 + \frac{1}{36} \mathcal{F}_v \mathcal{F}_{xtt}$}}} \vspace{1mm} \\ 
{\colorbox{shade_2}
{\parbox{0.55cm}{\centering $\frac{1}{4!} c_4$}}} & {\colorbox{shade_2}
{\parbox{14.25cm}{$\frac{1}{24} \mathcal{F} \mathcal{F}_u^3 + \frac{1}{36} \mathcal{F} \mathcal{F}_u^2 \mathcal{F}_{uv} v + \frac{13}{72} \mathcal{F} \mathcal{F}_u \mathcal{F}_{uv} \mathcal{F}_v w_1 + \frac{1}{36} \mathcal{F} \mathcal{F}_u \mathcal{F}_{uv} \mathcal{F}_x + \frac{1}{36} \mathcal{F} \mathcal{F}_u \mathcal{F}_{uvx} \mathcal{F}_v v + \frac{1}{8} \mathcal{F} \mathcal{F}_u \mathcal{F}_{ux} \mathcal{F}_v \vspace{1mm} \\ + \frac{7}{72} \mathcal{F} \mathcal{F}_{uv}^2 \mathcal{F}_v v w_1 + \frac{7}{72} \mathcal{F} \mathcal{F}_{uv} \mathcal{F}_{ux} \mathcal{F}_v v + \frac{5}{72} \mathcal{F} \mathcal{F}_{uv} \mathcal{F}_v^2 w_2 + \frac{7}{72} \mathcal{F} \mathcal{F}_{uv} \mathcal{F}_v \mathcal{F}_{vx} w_1 + \frac{1}{36} \mathcal{F} \mathcal{F}_{uv} \mathcal{F}_v \mathcal{F}_{xx} + \frac{1}{9} \mathcal{F} \mathcal{F}_{uvx} \mathcal{F}_v^2 w_1 \vspace{1mm} \\ + \frac{1}{36} \mathcal{F} \mathcal{F}_{uvx} \mathcal{F}_v \mathcal{F}_x + \frac{1}{24} \mathcal{F} \mathcal{F}_{ux} \mathcal{F}_v \mathcal{F}_{vx} + \frac{1}{24} \mathcal{F} \mathcal{F}_{uxx} \mathcal{F}_v^2 + \frac{1}{24} \mathcal{F}_t \mathcal{F}_u^2 + \frac{1}{24} \mathcal{F}_t \mathcal{F}_{uv} \mathcal{F}_v w_1 + \frac{1}{24} \mathcal{F}_t \mathcal{F}_{ux} \mathcal{F}_v + \frac{1}{72} \mathcal{F}_{tt} \mathcal{F}_u \vspace{1mm} \\ + \frac{1}{432} \mathcal{F}_{ttt} + \frac{1}{8} \mathcal{F}_u^3 \mathcal{F}_v v + \frac{11}{72} \mathcal{F}_u^2 \mathcal{F}_{uv} \mathcal{F}_v v^2 + \frac{1}{4} \mathcal{F}_u^2 \mathcal{F}_v^2 w_1 + \frac{1}{8} \mathcal{F}_u^2 \mathcal{F}_v \mathcal{F}_{vx} v + \frac{1}{8} \mathcal{F}_u^2 \mathcal{F}_v \mathcal{F}_x + \frac{1}{24} \mathcal{F}_u \mathcal{F}_{uv}^2 \mathcal{F}_v v^3 \vspace{1mm} \\ + \frac{13}{18} \mathcal{F}_u \mathcal{F}_{uv} \mathcal{F}_v^2 v w_1 + \frac{1}{12} \mathcal{F}_u \mathcal{F}_{uv} \mathcal{F}_v \mathcal{F}_{vx} v^2 + \frac{13}{72} \mathcal{F}_u \mathcal{F}_{uv} \mathcal{F}_v \mathcal{F}_x v + \frac{1}{12} \mathcal{F}_u \mathcal{F}_{uvx} \mathcal{F}_v^2 v^2 + \frac{3}{8} \mathcal{F}_u \mathcal{F}_{ux} \mathcal{F}_v^2 v + \frac{1}{8} \mathcal{F}_u \mathcal{F}_v^3 w_2 \vspace{1mm} \\ + \frac{1}{2} \mathcal{F}_u \mathcal{F}_v^2 \mathcal{F}_{vx} w_1 + \frac{1}{24} \mathcal{F}_u \mathcal{F}_v^2 \mathcal{F}_{vxx} v + \frac{1}{8} \mathcal{F}_u \mathcal{F}_v^2 \mathcal{F}_{xx} + \frac{1}{24} \mathcal{F}_u \mathcal{F}_v \mathcal{F}_{vx}^2 v + \frac{1}{8} \mathcal{F}_u \mathcal{F}_v \mathcal{F}_{vx} \mathcal{F}_x + \frac{1}{12} \mathcal{F}_u \mathcal{F}_v \mathcal{F}_{xt} + \frac{7}{24} \mathcal{F}_{uv}^2 \mathcal{F}_v^2 v^2 w_1 \vspace{1mm} \\ + \frac{1}{24} \mathcal{F}_{uv}^2 \mathcal{F}_v \mathcal{F}_x v^2 + \frac{1}{4} \mathcal{F}_{uv} \mathcal{F}_{ux} \mathcal{F}_v^2 v^2 + \frac{1}{4} \mathcal{F}_{uv} \mathcal{F}_v^3 v w_2 + \frac{23}{72} \mathcal{F}_{uv} \mathcal{F}_v^3 w_1^2 + \frac{7}{12} \mathcal{F}_{uv} \mathcal{F}_v^2 \mathcal{F}_{vx} v w_1 + \frac{2}{9} \mathcal{F}_{uv} \mathcal{F}_v^2 \mathcal{F}_x w_1 \vspace{1mm} \\ + \frac{1}{8} \mathcal{F}_{uv} \mathcal{F}_v^2 \mathcal{F}_{xx} v + \frac{1}{12} \mathcal{F}_{uv} \mathcal{F}_v \mathcal{F}_{vx} \mathcal{F}_x v + \frac{1}{36} \mathcal{F}_{uv} \mathcal{F}_v \mathcal{F}_x^2 + \frac{1}{24} \mathcal{F}_{uv} \mathcal{F}_v \mathcal{F}_{xt} v + \frac{1}{3} \mathcal{F}_{uvx} \mathcal{F}_v^3 v w_1 + \frac{1}{12} \mathcal{F}_{uvx} \mathcal{F}_v^2 \mathcal{F}_x v \vspace{1mm} \\ + \frac{1}{4} \mathcal{F}_{ux} \mathcal{F}_v^3 w_1 + \frac{1}{24} \mathcal{F}_{ux} \mathcal{F}_v^3 w_2 + \frac{1}{4} \mathcal{F}_{ux} \mathcal{F}_v^2 \mathcal{F}_{vx} v + \frac{1}{8} \mathcal{F}_{ux} \mathcal{F}_v^2 \mathcal{F}_x + \frac{1}{8} \mathcal{F}_{uxx} \mathcal{F}_v^3 v + \frac{1}{24} \mathcal{F}_v^4 w_2 + \frac{1}{4} \mathcal{F}_v^3 \mathcal{F}_{vx} w_2 + \frac{1}{6} \mathcal{F}_v^3 \mathcal{F}_{vxx} w_1 \vspace{1mm} \\ + \frac{1}{24} \mathcal{F}_v^3 \mathcal{F}_{xxx} + \frac{7}{24} \mathcal{F}_v^2 \mathcal{F}_{vx}^2 w_1 + \frac{1}{8} \mathcal{F}_v^2 \mathcal{F}_{vx} \mathcal{F}_{xx} + \frac{1}{24} \mathcal{F}_v^2 \mathcal{F}_{vxx} \mathcal{F}_x + \frac{1}{24} \mathcal{F}_v^2 \mathcal{F}_{xxt} + \frac{1}{24} \mathcal{F}_v \mathcal{F}_{vx}^2 \mathcal{F}_x + \frac{1}{24} \mathcal{F}_v \mathcal{F}_{vx} \mathcal{F}_{xt} + \frac{1}{72} \mathcal{F}_v \mathcal{F}_{xtt}$}}} \vspace{1mm} \\ 
\bottomrule \\
\end{tabular} \label{FHat4c4WilliamsonLowStorageThirdOrderRungeKuttaMethod}
\end{table}

\noindent \textbf{Stage 3}
\begin{align}
\hat{u}_j^{n+1} &= \hat{u}^{n+\frac{3}{4}} + \frac{8}{15} \Delta t \left(\widetilde{\mathcal{F}}_{\text{mean}}\right)_j^{n+\frac{3}{4}} \nonumber \\
&= u_j^n + \frac{\Delta t}{1!} \left(\mathcal{F}^{(1)} + \bigO\left(\Delta x^{\alpha}\right)\right)_j^n + \frac{\Delta t^2}{2!} \left(\mathcal{F}^{(2)} + \bigO\left(\Delta x^{\alpha}\right)\right)_j^n \nonumber \\
&\hspace{0.35cm}+ \frac{\Delta t^3}{3!} \left(\mathcal{F}^{(3)} + \bigO\left(\Delta x^{\alpha}\right)\right)_j^n + \frac{\Delta t^4}{4!} \left(\widehat{\mathcal{F}}^{(4)} + \bigO\left(\Delta x^{\alpha}\right)\right)_j^n + \bigO\left(\Delta t^5\right),
\end{align}
where $\left(\widehat{\mathcal{F}}^{(4)}\right)_j^n \ne \left(\mathcal{F}^{(4)}\right)_j^n$. The full expressions for $\mathcal{F}_{\text{mean}} \left(\hat{u}_j^{n+\theta_k},\hat{v}_j^{n+\theta_k},x_j,t^{n+\theta_k}\right)$ are found by Taylor expanding $\mathcal{F} \left(\hat{u}_j^{n+\theta_k},\hat{v}_j^{n+\theta_k},x_j,t^{n+\theta_k}\right) \equiv \mathcal{F} \left(u_j^n + \Delta \hat{u}_j^{n+\theta_k}, v_j^n + \Delta \hat{v}_j^{n+\theta_k}, x_j, t^n+\theta_k \Delta t\right)$
about $u_j^n$, $v_j^n$, and $t^n$ for $k=1$, $2$ and $\theta_1 = \frac{1}{3}$, $\theta_2 = \frac{3}{4}$. After Stage 3, we can determine the local truncation error
\begin{align}
\hat{\tau}_j^{n+1} = u_j^{n+1} - \hat{u}_j^{n+1} &= \frac{\Delta t}{1!} \bigO\left(\Delta x^{\alpha}\right) + \frac{\Delta t^2}{2!} \bigO\left(\Delta x^{\alpha}\right) + \frac{\Delta t^3}{3!} \bigO\left(\Delta x^{\alpha}\right) + \frac{\Delta t^4}{4!} \left(c_4 + \bigO\left(\Delta x^{\alpha}\right)\right)_j^n + \bigO\left(\Delta t^5\right) \nonumber \\
&=\Delta t \bigO\left(\Delta x^{\alpha}\right) + \Delta t^2 \bigO\left(\Delta x^{\alpha}\right) + \Delta t^3 \bigO\left(\Delta x^{\alpha}\right) + \bigO\left(\Delta t^4\right),
\end{align}
where $\left(c_4\right)_j^n = \left(\mathcal{F}^{(4)}\right)_j^n - \left(\mathcal{F}^{\hat{(4)}}\right)_j^n \ne 0$. The full expressions for $\frac{1}{4!} \widehat{\mathcal{F}}^{(4)}$ and $\frac{1}{4!} c_4$ are in Table~\ref{FHat4c4WilliamsonLowStorageThirdOrderRungeKuttaMethod}. 

Summarizing, a predictor-corrector Runge-Kutta method of order $\beta$ results in $\left(\widehat{\mathcal{F}}^{(k)}\right)_j^n = \left(\mathcal{F}^{(k)}\right)_j^n$ for $k = 1,2,\ldots,\beta$, and $\left(\widehat{\mathcal{F}}^{(\beta+1)}\right)_j^n$ only consists of some of the terms in $\left(\mathcal{F}^{(\beta+1)}\right)_j^n$ but not necessarily with the same multiplicative factors, so $\left(\widehat{\mathcal{F}}^{(\beta+1)}\right)_j^n \ne \left(\mathcal{F}^{(\beta+1)}\right)_j^n$.

\subsection{Adams-Bashforth Time-Stepping Methods}

We now consider Adams-Bashforth methods, which involve the solution $u_j^{n-m} = u_j^n + \Delta u_j^{n-m}$ and its spatial derivative $v_j^{n-m} = v_j^n + \Delta v_j^{n-m}$, at spatial location~\ $x_j$ and time level $t^{n-m}$, where 
\begin{align}
\Delta u_j^{n-m} = \sum_{k=1}^{\infty} \frac{(-m \Delta t)^k}{k!} \left(\mathcal{F}^{(k)}\right)_j^n \quad \text{and} \quad
\Delta v_j^{n-m} = \sum_{k=1}^{\infty} \frac{(-m \Delta t)^k}{k!} \left(\mathcal{G}^{(k)}\right)_j^n.
\end{align}

Second-order Adams-Bashforth results in the numerical solution
\begin{align}
\hat{u}_j^{n+1} &= u_j^n + \Delta t \left\{\frac{3}{2} \mathcal{F}\left(u_j^n,v_j^n + \bigO\left(\Delta x^{\alpha}\right),x_j,t^n\right) - \frac{1}{2} \mathcal{F}\left(u_j^{n-1},v_j^{n-1} + \bigO\left(\Delta x^{\alpha}\right),x_j,t^{n-1}\right)\right\} \nonumber \\
&= u_j^n + \Delta t \left\{\frac{3}{2} \mathcal{F}\left(u_j^n,v_j^n + \bigO\left(\Delta x^{\alpha}\right),x_j,t^n\right) - \frac{1}{2} \mathcal{F}\left(u_j^n + \Delta u_j^{n-1},v_j^n + \Delta v_j^{n-1} + \bigO\left(\Delta x^{\alpha}\right),x_j,t^n - \Delta t\right)\right\} \nonumber \\
&= u_j^n + \frac{\Delta t}{1!} \left(\mathcal{F}^{(1)} + \bigO\left(\Delta x^{\alpha}\right)\right)_j^n + \frac{\Delta t^2}{2!} \left(\mathcal{F}^{(2)} + \bigO\left(\Delta x^{\alpha}\right)\right)_j^n + \frac{\Delta t^3}{3!}\left(\widehat{\mathcal{F}}^{(3)} + \bigO\left(\Delta x^{\alpha}\right)\right)_j^n + \bigO\left(\Delta t^4\right),
\end{align}
and the local truncation error
\begin{align}
\hat{\tau}_j^{n+1} = u_j^{n+1} - \hat{u}_j^{n+1} &= \frac{\Delta t}{1!} \bigO\left(\Delta x^{\alpha}\right) + \frac{\Delta t^2}{2!} \bigO\left(\Delta x^{\alpha}\right) + \frac{\Delta t^3}{3!} \left(c_3 + \bigO\left(\Delta x^{\alpha}\right)\right)_j^n + \bigO\left(\Delta t^4\right) \nonumber \\
&= \Delta t \bigO\left(\Delta x^{\alpha}\right) + \Delta t^2 \bigO\left(\Delta x^{\alpha}\right) + \bigO\left(\Delta t^3\right),
\end{align}
where $\left(\widehat{\mathcal{F}}^{(3)}\right)_j^n = -\frac{3}{2} \left(\mathcal{F}^{(3)}\right)_j^n \ne \left(\mathcal{F}^{(3)}\right)_j^n$, and $\left(c_3\right)_j^n = \left(\mathcal{F}^{(3)}\right)_j^n - \left(\widehat{\mathcal{F}}^{(3)}\right)_j^n = \frac{5}{2} \left(\mathcal{F}^{(3)}\right)_j^n \ne 0$. 

Third-order Adams-Bashforth method results in the numerical solution
\begin{align}
\hat{u}_j^{n+1} &= u_j^n + \Delta t \left\{\frac{23}{12} \mathcal{F}\left(u_j^n,v_j^n + \bigO\left(\Delta x^{\alpha}\right),x_j,t^n\right) - \frac{16}{12} \mathcal{F}\left(u_j^{n-1},v_j^{n-1} + \bigO\left(\Delta x^{\alpha}\right),x_j,t^{n-1}\right)\right. \nonumber \\
&\hspace{1.625cm} \left.+ \frac{5}{12} \mathcal{F}\left(u_j^{n-2},v_j^{n-2} + \bigO\left(\Delta x^{\alpha}\right),x_j,t^{n-2}\right)\right\} \nonumber \\
&= u_j^n + \Delta t \left\{\frac{23}{12} \mathcal{F}\left(u_j^n,v_j^n,t^n\right) - \frac{16}{12} \mathcal{F}\left(u_j^n + \Delta u_j^{n-1},v_j^n + \Delta v_j^{n-1} + \bigO\left(\Delta x^{\alpha}\right),t^n - \Delta t\right)\right. \nonumber \\
&\hspace{1.625cm} \left.+ \frac{5}{12} \mathcal{F}\left(u_j^n + \Delta u_j^{n-2},v_j^n + \Delta v_j^{n-2} + \bigO\left(\Delta x^{\alpha}\right),t^n - 2 \Delta t\right) + \bigO\left(\Delta x^{\alpha}\right)\right\} \nonumber \\
&= u_j^n + \frac{\Delta t}{1!} \left(\mathcal{F}^{(1)} + \bigO\left(\Delta x^{\alpha}\right)\right)_j^n + \frac{\Delta t^2}{2!} \left(\mathcal{F}^{(2)} + \bigO\left(\Delta x^{\alpha}\right)\right)_j^n \nonumber \\ 
&\hspace{0.35cm}+ \frac{\Delta t^3}{3!}\left(\mathcal{F}^{(3)} + \bigO\left(\Delta x^{\alpha}\right)\right)_j^n + \frac{\Delta t^4}{4!}\left(\widehat{\mathcal{F}}^{(4)} + \bigO\left(\Delta x^{\alpha}\right)\right)_j^n + \bigO\left(\Delta t^5\right),
\end{align}
and the local truncation error 
\begin{align}
\hat{\tau}_j^{n+1} = u_j^{n+1} - \hat{u}_j^{n+1} &= \frac{\Delta t}{1!} \bigO\left(\Delta x^{\alpha}\right) + \frac{\Delta t^2}{2!} \bigO\left(\Delta x^{\alpha}\right) + \frac{\Delta t^3}{3!} \bigO\left(\Delta x^{\alpha}\right) + \frac{\Delta t^4}{4!} \left(c_4 + \bigO\left(\Delta x^{\alpha}\right)\right)_j^n + \bigO\left(\Delta t^5\right) \nonumber \\
&= \Delta t \bigO\left(\Delta x^{\alpha}\right) + \Delta t^2 \bigO\left(\Delta x^{\alpha}\right) + \Delta t^3 \bigO\left(\Delta x^{\alpha}\right) + \bigO\left(\Delta t^4\right),
\end{align}
where $\left(\widehat{\mathcal{F}}^{(4)}\right)_j^n = -8 \left(\mathcal{F}^{(4)}\right)_j^n \ne \left(\mathcal{F}^{(4)}\right)_j^n$, and $\left(c_4\right)_j^n = \left(\mathcal{F}^{(4)}\right)_j^n - \left(\widehat{\mathcal{F}}^{(4)}\right)_j^n = 9 \left(\mathcal{F}^{(4)}\right)_j^n \ne 0$.

Summarizing, when applying a $\beta$-order Adams-Bashforth method, we obtain $\left(\widehat{\mathcal{F}}^{(k)}\right)_j^n = \left(\mathcal{F}^{(k)}\right)_j^n$ for $k = 1,2,\ldots,\beta$, and $\left(\widehat{\mathcal{F}}^{(\beta+1)}\right)_j^n = \gamma \left(\mathcal{F}^{(\beta+1)}\right)_j^n \ne \left(\mathcal{F}^{(\beta+1)}\right)_j^n$ where $\gamma \ne 1$.

\subsection{Implicit Time-Stepping Methods}

We consider the three implicit time-stepping methods in List~\ref{myListOfImplicitTimeSteppingMethodsForAnalysis}. Each method's local truncation error involves the solution $u_j^{n+1} = u_j^n + \Delta u_j^{n+1}$ and its spatial derivative $v_j^{n+1} = v_j^n + \Delta v_j^{n+1}$ at spatial location $x_j$ and time level $t^{n+1}$, where 
\begin{align}
\Delta u_j^{n+1} = \sum_{k=1}^{\infty} \frac{(\Delta t)^k}{k!} \left(\mathcal{F}^{(k)}\right)_j^n \quad \text{and} \quad
\Delta v_j^{n+1} = \sum_{k=1}^{\infty} \frac{(\Delta t)^k}{k!} \left(\mathcal{G}^{(k)}\right)_j^n.
\end{align}
Applying the first-order Backward Euler method,
\begin{align}
\hat{u}_j^{n+1} &= u_j^n + \Delta t \mathcal{F} \left(u_j^{n+1},v_j^{n+1}+O\left(\Delta x^{\alpha}\right),x_j,t^{n+1}\right) \equiv u_j^n + \Delta t \mathcal{F} \left(u_j^n + \Delta u_j^{n+1},v_j^n + \Delta v_j^{n+1} + \bigO\left(\Delta x^{\alpha}\right),x_j,t^n + \Delta t\right) \nonumber \\
&= u_j^n + \frac{\Delta t}{1!} \left(\mathcal{F}^{(1)} + \bigO\left(\Delta x^{\alpha}\right)\right)_j^n + \frac{\Delta t^2}{2!} \left(\widehat{\mathcal{F}}^{(2)} + \bigO\left(\Delta x^{\alpha}\right)\right)_j^n + \bigO\left(\Delta t^3\right),
\end{align}
where $\left(\widehat{\mathcal{F}}^{(2)}\right)_j^n = 2 \left(\mathcal{F}^{(2)}\right)_j^n \ne \left(\mathcal{F}^{(2)}\right)_j^n$. The local truncation error is
\begin{equation}
\hat{\tau}_j^{n+1} = u_j^{n+1} - \hat{u}_j^{n+1} = \frac{\Delta t}{1!} \bigO\left(\Delta x^{\alpha}\right) + \frac{\Delta t^2}{2!} \left(c_2 + \bigO\left(\Delta x^{\alpha}\right)\right)_j^n + \bigO\left(\Delta t^3\right) = \Delta t \bigO\left(\Delta x^{\alpha}\right) + \bigO\left(\Delta t^2\right),
\end{equation}
where $\left(c_2\right)^n = \left(\mathcal{F}^{(2)}\right)_j^n - \left(\widehat{\mathcal{F}}^{(2)}\right)_j^n = -\left(\mathcal{F}^{(2)}\right)_j^n \ne 0$.

\begin{table}[!htp]
\centering
\caption{The term $\frac{1}{3!} \widehat{\mathcal{F}}^{(3)}$ in the numerical solution $\hat{u}_j^{n+1}$ given by \eqref{Exact_Solution_Time_Level_nP1_Theorem}, and the term $\frac{1}{3!} c_3$ in the local truncation error $\hat{\tau}_j^{n+1}$ given by \eqref{LocalTruncationErrorNumericalSolutionFinalForm_1}, of the generic hyperbolic PDE~\eqref{AdvectionEquation1DFunctionalForm} with approximations~\eqref{myApproximationsForPDEsUsedToModelPhysicalPhenomena}, advanced in time with the implicit midpoint method, at spatial location~$x_j$ and time level $t^{n+1} = t^n + \Delta t$, expressed as functions of quantities known at the current time level~$t^n$, and $v=u_x$, $w_1 = u_{xx}$, and $w_2 = u_{xxx}$ for notational convenience.}
\vspace{3mm}
\setlength{\tabcolsep}{0.35em}
\begin{tabular}{cc} 
\toprule
{\colorbox{shade_1}
{\parbox{0.85cm}{\centering $\frac{1}{3!} \widehat{\mathcal{F}}^{(3)}$}}} & {\colorbox{shade_1}
{\parbox{11.25cm}{$\frac{1}{4} \mathcal{F} \mathcal{F}_u^2 + \frac{1}{4} \mathcal{F} \mathcal{F}_u \mathcal{F}_{uv} v + \frac{1}{2} \mathcal{F} \mathcal{F}_{uv} \mathcal{F}_v w_1 + \frac{1}{4} \mathcal{F} \mathcal{F}_{uv} \mathcal{F}_x + \frac{1}{4} \mathcal{F} \mathcal{F}_{ux} \mathcal{F}_v + \frac{1}{4} \mathcal{F}_t \mathcal{F}_u + \frac{1}{8} \mathcal{F}_{tt} \vspace{1mm} \\ + \frac{1}{2} \mathcal{F}_u^2 \mathcal{F}_v v + \frac{1}{4} \mathcal{F}_u \mathcal{F}_{uv} \mathcal{F}_v v^2 + \frac{3}{4} \mathcal{F}_u \mathcal{F}_v^2 w_1 + \frac{1}{4} \mathcal{F}_u \mathcal{F}_v \mathcal{F}_{vx} v + \frac{1}{2} \mathcal{F}_u \mathcal{F}_v \mathcal{F}_x + \frac{3}{4} \mathcal{F}_{uv} \mathcal{F}_v^2 v w_1 \vspace{1mm} \\ + \frac{1}{4} \mathcal{F}_{uv} \mathcal{F}_v \mathcal{F}_x v + \frac{1}{2} \mathcal{F}_{ux} \mathcal{F}_v^2 v + \frac{1}{4} \mathcal{F}_v^3 w_2 + \frac{3}{4} \mathcal{F}_v^2 \mathcal{F}_{vx} w_1 + \frac{1}{4} \mathcal{F}_v^2 \mathcal{F}_{xx} + \frac{1}{4} \mathcal{F}_v \mathcal{F}_{vx} \mathcal{F}_x + \frac{1}{4} \mathcal{F}_v \mathcal{F}_{xt}$}}} \vspace{1mm} \\ 
{\colorbox{shade_2}
{\parbox{0.55cm}{\centering $\frac{1}{3!} c_3$}}} & {\colorbox{shade_2}
{\parbox{11.25cm}{$-\frac{1}{12} \mathcal{F} \mathcal{F}_u^2 + \frac{1}{12} \mathcal{F} \mathcal{F}_u \mathcal{F}_{uv} v + \frac{1}{12} \mathcal{F} \mathcal{F}_{uv} \mathcal{F}_x - \frac{1}{12} \mathcal{F} \mathcal{F}_{ux} \mathcal{F}_v - \frac{1}{12} \mathcal{F}_t \mathcal{F}_u + \frac{1}{24} \mathcal{F}_{tt} - \frac{1}{6} \mathcal{F}_u^2 \mathcal{F}_v v \vspace{1mm} \\ - \frac{1}{12} \mathcal{F}_u \mathcal{F}_{uv} \mathcal{F}_v v^2 - \frac{1}{4} \mathcal{F}_u \mathcal{F}_v^2 w_1 - \frac{1}{12} \mathcal{F}_u \mathcal{F}_v \mathcal{F}_{vx} v - \frac{1}{6} \mathcal{F}_u \mathcal{F}_v \mathcal{F}_x - \frac{1}{4} \mathcal{F}_{uv} \mathcal{F}_v^2 v w_1 - \frac{1}{12} \mathcal{F}_{uv} \mathcal{F}_v \mathcal{F}_x v \vspace{1mm} \\ - \frac{1}{6} \mathcal{F}_{ux} \mathcal{F}_v^2 v - \frac{1}{12} \mathcal{F}_v^3 w_2 - \frac{1}{4} \mathcal{F}_v^2 \mathcal{F}_{vx} w_1 - \frac{1}{12} \mathcal{F}_v^2 \mathcal{F}_{xx} - \frac{1}{12} \mathcal{F}_v \mathcal{F}_{vx} \mathcal{F}_x - \frac{1}{12} \mathcal{F}_v \mathcal{F}_{xt}$}}} \vspace{1mm} \\
\bottomrule \\
\end{tabular} \label{FHat3c3ImplicitMidpointMethod}
\end{table}

The second-order predictor-corrector implicit midpoint method is
\begin{align}
\hat{u}_j^{n+1} &= u_j^n + \Delta t \mathcal{F} \left(\frac{1}{2} \left(u_j^n + u_j^{n+1}\right),\frac{1}{2} \left(v_j^n + v_j^{n+1}\right) + \bigO\left(\Delta x^{\alpha}\right),t^{n+\frac{1}{2}}\right) \nonumber \\
&\equiv u_j^n + \Delta t \mathcal{F} \left(u_j^n + \frac{1}{2} \Delta u_j^{n+1},v_j^n + \frac{1}{2} \Delta v_j^{n+1} + \bigO\left(\Delta x^{\alpha}\right),t^n + \frac{\Delta t}{2}\right) \nonumber \\
&= u_j^n + \frac{\Delta t}{1!} \left(\mathcal{F}^{(1)} + \bigO\left(\Delta x^{\alpha}\right)\right)_j^n + \frac{\Delta t^2}{2!} \left(\mathcal{F}^{(2)} + \bigO\left(\Delta x^{\alpha}\right)\right)_j^n + \frac{\Delta t^3}{3!} \left(\widehat{\mathcal{F}}^{(3)} + \bigO\left(\Delta x^{\alpha}\right)\right)_j^n + \bigO\left(\Delta t^4\right),
\end{align}
with the local truncation error 
\begin{align}
\hat{\tau}_j^{n+1} = u_j^{n+1} - \hat{u}_j^{n+1} &= \frac{\Delta t}{1!} \bigO\left(\Delta x^{\alpha}\right) + \frac{\Delta t^2}{2!} \bigO\left(\Delta x^{\alpha}\right) + \frac{\Delta t^3}{3!} \left(c_3 + \bigO\left(\Delta x^{\alpha}\right)\right)_j^n + \bigO\left(\Delta t^4\right) \nonumber \\
&= \Delta t \bigO\left(\Delta x^{\alpha}\right) + \Delta t^2 \bigO\left(\Delta x^{\alpha}\right) + \bigO\left(\Delta t^3\right),
\end{align}
where $\left(\widehat{\mathcal{F}}^{(3)}\right)_j^n \ne \left(\mathcal{F}^{(3)}\right)_j^n$ and $\left(c_3\right)_j^n = \left(\mathcal{F}^{(3)}\right)_j^n - \left(\widehat{\mathcal{F}}^{(3)}\right)_j^n \ne 0$. The full expressions for $\left(\mathcal{F}^{(3)}\right)_j^n$ and $\left(c_3\right)_j^n$ are in Table~\ref{FHat3c3ImplicitMidpointMethod}.

Finally, the trapezoidal rule results in the numerical solution
\begin{align}
\hat{u}_j^{n+1} &= u_j^n + \frac{\Delta t}{2} \left\{\mathcal{F} \left(u_j^n,v_j^n + \bigO\left(\Delta x^{\alpha}\right),x_j,t^n\right) + \mathcal{F} \left(u_j^{n+1},v_j^{n+1} + \bigO\left(\Delta x^{\alpha}\right),x_j,t^{n+1}\right)\right\} \nonumber \\
&\equiv u_j^n + \frac{\Delta t}{2} \left\{\mathcal{F} \left(u_j^n,v_j^n + \bigO\left(\Delta x^{\alpha}\right),x_j,t^n\right) + \mathcal{F} \left(u_j^n + \Delta u_j^{n+1},v_j^n + \Delta v_j^{n+1} + \bigO\left(\Delta x^{\alpha}\right),x_j,t^n + \Delta t\right)\right\} \nonumber \\
&= u_j^n + \frac{\Delta t}{1!} \left(\mathcal{F}^{(1)} + \bigO\left(\Delta x^{\alpha}\right)\right)_j^n + \frac{\Delta t^2}{2!} \left(\mathcal{F}^{(2)} + \bigO\left(\Delta x^{\alpha}\right)\right)_j^n + \frac{\Delta t^3}{3!} \left(\widehat{\mathcal{F}}^{(3)} + \bigO\left(\Delta x^{\alpha}\right)\right)_j^n + \bigO\left(\Delta t^4\right),
\end{align}
and the local truncation error
\begin{align}
\hat{\tau}_j^{n+1} = u_j^{n+1} - \hat{u}_j^{n+1} &= \frac{\Delta t}{1!} \bigO\left(\Delta x^{\alpha}\right) + \frac{\Delta t^2}{2!} \bigO\left(\Delta x^{\alpha}\right) + \frac{\Delta t^3}{3!} \left(c_3 + \bigO\left(\Delta x^{\alpha}\right)\right)_j^n + \bigO\left(\Delta t^4\right) \nonumber \\
&= \Delta t \bigO\left(\Delta x^{\alpha}\right) + \Delta t^2 \bigO\left(\Delta x^{\alpha}\right) + \bigO\left(\Delta t^3\right),
\end{align}
where $\left(\widehat{\mathcal{F}}^{(3)}\right)_j^n = \frac{3}{2} \left(\mathcal{F}^{(3)}\right)_j^n \ne \left(\mathcal{F}^{(3)}\right)_j^n$, and $c_3^n = \left(\mathcal{F}^{(3)}\right)_j^n - \left(\widehat{\mathcal{F}}^{(3)}\right)_j^n = -\frac{1}{2} \left(\mathcal{F}^{(3)}\right)_j^n \ne 0$.

Based on our examples, an implicit time-stepping method of order $\beta$ results in $\left(\widehat{\mathcal{F}}^{(k)}\right)_j^n = \left(\mathcal{F}^{(k)}\right)_j^n$ for $k = 1,2,\ldots,\beta$. For predictor-corrector methods, the term $\left(\widehat{\mathcal{F}}^{(\beta+1)}\right)_j^n$ consists of all the terms as in $\left(\mathcal{F}^{(\beta+1)}\right)_j^n$, but mostly with different pre-factors, and for multistep methods, the term $\left(\widehat{\mathcal{F}}^{(\beta+1)}\right)_j^n$ is a scalar multiple of $\left(\mathcal{F}^{(\beta+1)}\right)_j^n$, with the scalar factor not equal to one. In either case, $\left(\widehat{\mathcal{F}}^{(\beta+1)}\right)_j^n \ne \left(\mathcal{F}^{(\beta+1)}\right)_j^n$.

\subsection{Local and Global Truncation Errors of a Hyperbolic PDE}

We have observed that the local truncation error of both the predictor-corrector and the multistep time-stepping methods when applied to the first-order hyperbolic PDE have the same generic form. The result is the following Theorem.
\begin{thm} \label{Theorem_1}
The global truncation error of a hyperbolic PDE $u_t = \mathcal{F}(u,u_x,x,t)$ on a uniform mesh with spacing $\Delta x$ after an integral number of time steps of magnitude $\Delta t$ is
\begin{equation}
\hat{\tau}_G = \bigO\left({\Delta x}^{\alpha}\right) + \Delta t \bigO\left({\Delta x}^{\alpha}\right) + {\Delta t}^2 \bigO\left({\Delta x}^{\alpha}\right) + \cdots + {\Delta t}^{\beta-1} \bigO\left({\Delta x}^{\alpha}\right) + \bigO\left({\Delta t}^{\beta}\right), \label{GlobalTruncationErrorNumericalSolutionFinalForm}
\end{equation}
which reduces to
\begin{equation}
\hat{\tau}_G \approx \bigO\left({\Delta x}^{\alpha}\right) + \bigO\left({\Delta t}^{\beta}\right) \text{ for $\Delta t \ll 1$}. \label{GlobalTruncationErrorNumericalSolutionFinalForm_ApproximationForDeltaTMuchLessThanOne}
\end{equation}
\begin{proof}
Given the exact solution $u_j^n$ of a hyperbolic PDE $u_t = \mathcal{F}(u,u_x,x,t)$ on a uniform mesh with spacing $\Delta x$, at spatial locations $x_j$ for $j=1,2,\ldots$, and at time level $t^n$, the exact solution at time level $t^{n+1} = t^n + \Delta t$ may be obtained by Taylor expanding $u_j^n$ about time level $t^n$ as
\begin{equation}
u_j^{n+1} = u_j^n + \sum \limits_{k=1}^{\infty} \frac{\Delta t^k}{k!} \left(\frac{\del^k u}{\del t^k}\right)_j^n \equiv u_j^n + \sum \limits_{k=1}^{\infty} \frac{\Delta t^k}{k!} \left(\mathcal{F}^{(k)}\right)_j^n, \label{Exact_Solution_Time_Level_nP1_Theorem}
\end{equation}
where $\left(\mathcal{F}^{(k)}\right)_j^n = \left(\frac{\del^k u}{\del t^k}\right)_j^n$ is the $k^{\text{th}}$-order spatial derivative at $x_j$ and $t^n$. The numerical solution at time level $t^{n+1}$, obtained with a time-stepping method belonging to the Method of Lines, may be written in the general form
\begin{equation}
\hat{u}_j^{n+1} = u_j^n + \sum \limits_{k=1}^{\infty} \frac{\Delta t^k}{k!} \left(\widehat{\mathcal{F}}^{(k)} + \bigO\left(\Delta x^{\alpha}\right)\right)_j^n, \label{Numerical_Solution_Time_Level_nP1_Theorem}
\end{equation}
where $\alpha$ is the order of the spatial discretization and $\widehat{\mathcal{F}}^{(k)}$ is specified by the time-stepping method. If $\beta$ represents the order of the time-stepping method,
\begin{equation}
\left(\widehat{\mathcal{F}}^{(k)}\right)_j^n = \left(\mathcal{F}^{(k)}\right)_j^n \equiv \left(\frac{\del^k u}{\del t^k}\right)_j^n \mbox{ for } k = 1,2,\ldots,\beta. \label{FHat_F_Theorem}
\end{equation}
Subtracting \eqref{Numerical_Solution_Time_Level_nP1_Theorem} from \eqref{Exact_Solution_Time_Level_nP1_Theorem}, we obtain the local truncation error
\begin{align}
\hat{\tau}_j^{n+1} &= u_j^{n+1} - \hat{u}_j^{n+1} \nonumber \\
&= \sum \limits_{k=1}^{\infty} \frac{\Delta t^k}{k!} \left\{\left(\mathcal{F}^{(k)}\right)_j^n - \left(\widehat{\mathcal{F}}^{(k)}\right)_j^n + \bigO\left(\Delta x^{\alpha}\right)\right\}. \label{Local_Truncation_Error_Time_Level_nP1_Theorem}
\end{align}
Combining~\eqref{Local_Truncation_Error_Time_Level_nP1_Theorem} and \eqref{FHat_F_Theorem},
\begin{align}
\hat{\tau}_j^{n+1} &= \sum \limits_{k=1}^{\beta} \frac{\Delta t^k}{k!} \bigO\left(\Delta x^{\alpha}\right) + \sum \limits_{k=\beta+1}^{\infty} \frac{\Delta t^k}{k!} \left\{\left(\mathcal{F}^{(k)}\right)_j^n - \left(\widehat{\mathcal{F}}^{(k)}\right)_j^n + \bigO\left(\Delta x^{\alpha}\right)\right\} \nonumber \\
&= \frac{\Delta t}{1!} \bigO\left(\Delta x^{\alpha}\right) + \frac{\Delta t^2}{2!} \bigO\left(\Delta x^{\alpha}\right) + \frac{\Delta t^3}{3!} \bigO\left(\Delta x^{\alpha}\right) + \cdots + \frac{\Delta t^{\beta}}{{\beta}!} \bigO\left(\Delta x^{\alpha}\right) + \frac{\Delta t^{\beta+1}}{(\beta+1)!} \left(c_{\beta+1} + \bigO\left(\Delta x^{\alpha}\right)\right)_j^n + \bigO\left(\Delta t^{\beta+2}\right) \label{LocalTruncationErrorNumericalSolutionFinalForm_1} \\
&= \Delta t \bigO\left({\Delta x}^{\alpha}\right) + {\Delta t}^2 \bigO\left({\Delta x}^{\alpha}\right) + {\Delta t}^3 \bigO\left({\Delta x}^{\alpha}\right) + \cdots + {\Delta t}^{\beta} \bigO\left({\Delta x}^{\alpha}\right) + \bigO\left({\Delta t}^{\beta+1}\right), \label{LocalTruncationErrorNumericalSolutionFinalForm}
\end{align}
where $\left(c_{\beta+1}\right)_j^n = \left(\mathcal{F}^{(\beta+1)}\right)_j^n - \left(\widehat{\mathcal{F}}^{(\beta+1)}\right)_j^n \ne 0$. The global truncation error at a time horizon, after an integral number of time steps, is one order of $\Delta t$ less than its local counterpart, and can be expressed as
\begin{equation}
\left(\hat{\tau}_G\right)_j = \bigO\left({\Delta x}^{\alpha}\right) + \Delta t \bigO\left({\Delta x}^{\alpha}\right) + {\Delta t}^2 \bigO\left({\Delta x}^{\alpha}\right) + \cdots + {\Delta t}^{\beta-1} \bigO\left({\Delta x}^{\alpha}\right) + \bigO\left({\Delta t}^{\beta}\right), \label{GlobalTruncationErrorNumericalSolutionFinalForm_GivenSpatialLocation}
\end{equation}
which reduces to 
\begin{equation}
\left(\hat{\tau}_G\right)_j = \bigO\left({\Delta x}^{\alpha}\right) + \bigO\left({\Delta t}^{\beta}\right) \text{ for $\Delta t \ll 1$}. \label{GlobalTruncationErrorNumericalSolutionFinalForm_ApproximationForDeltaTMuchLessThanOne_GivenSpatialLocation}
\end{equation}
Replacing $\left(\hat{\tau}_G\right)_j$ by its norm over all spatial locations $x_j$, \eqref{GlobalTruncationErrorNumericalSolutionFinalForm_GivenSpatialLocation} and \eqref{GlobalTruncationErrorNumericalSolutionFinalForm_ApproximationForDeltaTMuchLessThanOne_GivenSpatialLocation} become \eqref{GlobalTruncationErrorNumericalSolutionFinalForm} and \eqref{GlobalTruncationErrorNumericalSolutionFinalForm_ApproximationForDeltaTMuchLessThanOne} respectively.
\end{proof}
If we employ a stable numerical scheme and the global solution error is the same order of accuracy as the global truncation error, we arrive at the following Corollaries.
\end{thm}
\begin{cor} \label{Corollary_1}
The order of convergence of a hyperbolic PDE in the asymptotic regime at constant ratio of time step to cell width is specified by the minimum of the orders of the spatial and temporal discretizations.
\end{cor}
\begin{cor} \label{Corollary_2}
To achieve the maximum possible order of convergence in the asymptotic regime at constant ratio of time step to cell width, the time-stepping method used to advance a hyperbolic PDE should at least have the same order of accuracy as the spatial discretization. 
\end{cor}
\begin{cor} \label{Corollary_3}
The discretization of a hyperbolic PDE under only spatial or only temporal refinement in the asymptotic regime is not guaranteed to converge.
\end{cor}

We can compare the behavior of the local truncation error of the generic hyperbolic PDE with that of the generic ODE. We know that an order $\beta$ time-stepping method is constructed so that if the solution of an ODE were exact at time $t^n$, then the error at time step $t^{n+1} = t^n + \Delta t$ will be $\hat{\tau}^{n+1} = \bigO\left(\Delta t^{\beta + 1}\right)$. If we express the exact solution $u^{n+1}$ and the numerical solution $\hat{u}^{n+1}$ as polynomials in $\Delta t$, the coefficients of $\Delta t^k$ for $k = 0, 1, 2, \ldots, \beta$ in $\hat{u}^{n+1}$ match those of $u^{n+1}$. When the ODE is expressed in the generic form~\eqref{ODE1D}, these coefficients of $\Delta t^k$ for $k = 1, 2, \ldots$ consist of partial and mixed derivatives of the right-hand side term $\mathcal{F}(u,t)$, to be referred to as the tendency term from here onward, with respect to $u$ and $t$. To pinpoint the source of these derivatives, we recapitulate that $\hat{u}^{n+1}$ consists of tendency terms, either at intermediate time levels between (and including) $t^n$ and $t^{n+1}$ for a predictor-corrector time-stepping method, or at current and previous time levels $t^{n-k}$ for $k=0,1,2,\ldots$ for a multistep time-stepping method. When the tendency terms are Taylor expanded about $(u^n,t^n)$, the result is the above-mentioned partial and mixed $u$- and $t$-derivatives of $\mathcal{F}(u,t)$ in the polynomial expression for $\hat{u}^{n+1}$. After expressing the mixed and $t$-derivatives of $\mathcal{F}(u,t)$ as functions of known quantities at time level $t^n$, we observe that the coefficients of $\Delta t^k$ for $k=0,1,\ldots,\beta$ in this polynomial are equal to those of $u^{n+1}$, which in turn results in $\hat{\tau}^{n+1} = u^{n+1} - \hat{u}^{n+1} = \bigO\left(\Delta t^{\beta+1}\right)$. A fundamental reason for this result is that for given values of $u$ and $t$, the tendency term $\mathcal{F}(u,t)$ is exact for an ODE. 

The derivation of the local truncation error of a generic hyperbolic PDE involves the same operations as the generic ODE, but with one fundamental difference---the above-mentioned tendency terms are replaced with their spatially discretized versions. Since the time derivative of the dependent variable of a hyperbolic PDE is a function of the dependent and independent variables, and also the spatial derivatives of the dependent variable, we need to perform a discretization in space while computing the tendency term at any instant of time. Whenever we perform this operation, we introduce an $\bigO\left(\Delta x^{\alpha}\right)$ term, where $\alpha$ is the order of the spatial discretization. Denoting $v=u_x$, $w_1=u_{xx}$, $w_2=u_{xxx}$, $\ldots$, the Taylor expansion of the spatially discretized tendency terms centered at $u_j^n$, $v_j^n$, $\left(w_1\right)_j^n$, $\left(w_2\right)_j^n$, $\ldots$, $x_j$ and $t^n$, contain
\begin{enumerate}[label=(\alph*),noitemsep] 
\item the dependent variables $u$, $v$, $w_1$, $w_2$, $\ldots$ defined at $x_j$ and $t^n$;

\item the partial and mixed derivatives of $\mathcal{F}(u,v,x,t)$ with respect to the dependent and independent variables $u$, $v$, $x$, and $t$ all defined at $x_j$ and $t^n$;

\item the $\bigO\left(\Delta x^{\alpha}\right)$ terms.
\end{enumerate}
If the spatial discretization operator were exact, the $\bigO\left({\Delta x}^{\alpha}\right)$ terms would be absent. In other words, the discretization error would only consist of its temporal component, and the local truncation error would assume the form $\bigO\left({\Delta t}^{\beta+1}\right)$, identical to the result when applying an order $\beta$ time-stepping method to an ODE. However, since we cannot make this assumption for a general PDE, terms involving $\bigO\left({\Delta x}^{\alpha}\right)$ are expected to appear when we replace a tendency term with its spatially discretized version. This introduces an $\bigO\left({\Delta x}^{\alpha}\right)$ term in the coefficient of ${\Delta t}^k$ for $k = 1, 2, 3, \ldots$ in the final expression for $\hat{u}_j^{n+1}$ that is not present in the corresponding coefficient of ${\Delta t}^k$ in $u_j^{n+1}$. This $\bigO\left({\Delta x}^{\alpha}\right)$ term in the coefficient of ${\Delta t}^k$ for $k = 1, 2, 3, \ldots$ cannot be ignored in the final expression for $\hat{\tau}_j^{n+1}$ defined as the difference between $u_j^{n+1}$ and $\hat{u}_j^{n+1}$.

\subsection{Convergence at Constant Ratio of Time Step to Cell Width}
Here we assume that our numerical scheme is stable, and the global solution error is the same order of accuracy as the global truncation error. Then, in the asymptotic regime, where the magnitude of the truncation error is dominated by the powers of $\Delta t$ and $\Delta x$ rather than their coefficients, the order of convergence of the global solution error norm is the minimum of the $\bigO\left(\Delta x^{\alpha}\right)$ and the $\bigO\left(\Delta t^{\beta}\right)$ terms in~\eqref{GlobalTruncationErrorNumericalSolutionFinalForm}. If $\Delta t$ is proportional to $\Delta x$, meaning that the ratio of the time step to cell width is held fixed (or the Courant number is kept constant for a one-dimensional linear constant-coefficient advection problem), the order of the global solution error becomes
\begin{equation}
\hat{\tau}_G = \bigO\left({\Delta x}^{\alpha}\right) + \bigO\left({\Delta t}^{\beta}\right) =  \bigO\left({\Delta x}^{\alpha}\right) + \bigO\left({\gamma^{\beta} \Delta x}^{\beta}\right) = \bigO\left({\Delta x}^{\alpha}\right) + \bigO\left({\Delta x}^{\beta}\right) \approx \bigO\left({\Delta x}^{\text{min}({\alpha},{\beta})}\right),
\end{equation}
where $\gamma = \Delta t/\Delta x$. Therefore, the order of convergence can not exceed the order of the spatial discretization~$\alpha$, and to achieve this order of convergence, we need to apply a time-stepping method of order~$\beta$ with $\beta \ge \alpha$. Within a specific family of time-stepping methods, the most computationally efficient choice is $\beta = \alpha$. Table~\ref{Table_GobalSolutionErrorNorm} lists the order of convergence of the error norm in the asymptotic regime at constant ratio of $\Delta t$ to $\Delta x$ for varying orders of spatial and temporal discretizations.
\begin{table}[!tp]
\centering
\caption{Order of convergence of the error norm in the asymptotic regime at constant ratio of time step to cell width with spatial and temporal discretizations up to order four. FE denotes first-order Forward Euler, RK2, RK3, RK4 denote the second-, third-, fourth-order (predictor-corrector) Runge-Kutta methods, respectively, and AB2, AB3, AB4 denote the second-, third-, fourth-order (multistep) Adams-Bashforth methods, respectively.}
\vspace{3mm}
\begin{tabular}{cccc}
\toprule
{Order of} & {Time-Stepping} & {Order of} & {Order of Convergence of Error Norm in} \\ 
{Spatial} & {Method} & {Time-Stepping} & {Asymptotic Regime at Constant Ratio} \\ 
{Discretization} & {Employed} & {Method} & {of Time Step to Cell Width} \\ 
\midrule
1 & FE & 1 & 1 \\
1 & RK2 or AB2 & 2 & 1 \\
1 & RK3 or AB3 & 3 & 1 \\
1 & RK4 or AB4 & 4 & 1 \\
2 & FE & 1 & 1 \\
2 & RK2 or AB2 & 2 & 2 \\
2 & RK3 or AB3 & 3 & 2 \\
2 & RK4 or AB4 & 4 & 2 \\
3 & FE & 1 & 1 \\
3 & RK2 or AB2 & 2 & 2 \\
3 & RK3 or AB3 & 3 & 3 \\
3 & RK4 or AB4 & 4 & 3 \\
4 & FE & 1 & 1 \\
4 & RK2 or AB2 & 2 & 2 \\
4 & RK3 or AB3 & 3 & 3 \\
4 & RK4 or AB4 & 4 & 4 \\
\bottomrule \\
\end{tabular} \label{Table_GobalSolutionErrorNorm}
\end{table}

\subsection{Refinement Only in Space or Only in Time}
If only spatial or temporal refinement is performed with a stable numerical scheme, the leading order terms of the global solution error are the $\bigO\left(\Delta x^{\alpha}\right)$ and the $\bigO\left(\Delta t^{\beta}\right)$ terms. In this case, convergence cannot be guaranteed due to the $\bigO\left(\Delta t^{\beta}\right)$ term for refinement only in space, and due to the $\bigO\left(\Delta x^{\alpha}\right)$ term for refinement only in time. More specifically, under only spatial or temporal refinement, the global truncation error does not necessarily converge to zero. As a result, our numerical solution may not even be consistent, and convergence may be impossible.

Under certain circumstances, the magnitude of the global solution error norm can even increase with only temporal refinement. The simplest example is the one-dimensional linear homogeneous constant-coefficient advection equation
\begin{equation}
u_t + a u_x = 0, \label{LinearAdvecton1D_SimplestExample}
\end{equation}
discretized in space with the first-order upwind finite difference scheme and advanced in time with the first-order Forward Euler method. The local truncation error for this problem at spatial location $x_j$ and time level $t^{n+1}$ is
\begin{equation}
\hat{\tau}_j^{n+1} = \left(-\frac{1}{2} |a| \Delta t \Delta x + \frac{1}{2} a^2 \Delta t^2\right) \left(u_{xx}\right)_j^n + \cdots,
\end{equation}
where $|a|$ is the magnitude of the constant wave speed $a$. The leading order term of the local truncation error are diffusive in nature, and can be expressed as
\begin{equation}
\left[\hat{\tau}_j^{n+1}\right]_{\text{leading order}} = \left(-\frac{1}{2} |a| \Delta t \Delta x + \frac{1}{2} a^2 \Delta t^2\right) \left(u_{xx}\right)_j^n = -\frac{1}{2} |a| \Delta x \Delta t \left(1 - \frac{|a| \Delta t}{\Delta x}\right) \left(u_{xx}\right)_j^n \equiv -\frac{1}{2} |a| \Delta x \Delta t \left(1 - C\right) \left(u_{xx}\right)_j^n,
\end{equation}
where $C = |a| \Delta t/\Delta x$ is the Courant number, which is positive and must be less than one to ensure numerical stability. The global truncation error is one order of $\Delta t$ less, and can be approximated as
\begin{equation}
\left[\left(\hat{\tau}_G\right)_j\right]_{\text{leading order}} = -\frac{1}{2} |a| \Delta x \left(1 - \frac{|a| \Delta t}{\Delta x}\right) \left(u_{xx}\right)_j^n = -\frac{1}{2} |a| \Delta x \left(1 - C\right) \left(u_{xx}\right)_j^n.
\end{equation}
Maintaining $C < 1$, if $\Delta x$ is held constant and $\Delta t$ is refined, then $(1-C)$ increases towards $1$, and the magnitude of the global truncation error increases. Moreover, the error will be diffusive in nature.

\begin{figure}[!htp]
\centering
\includegraphics[scale=.35]{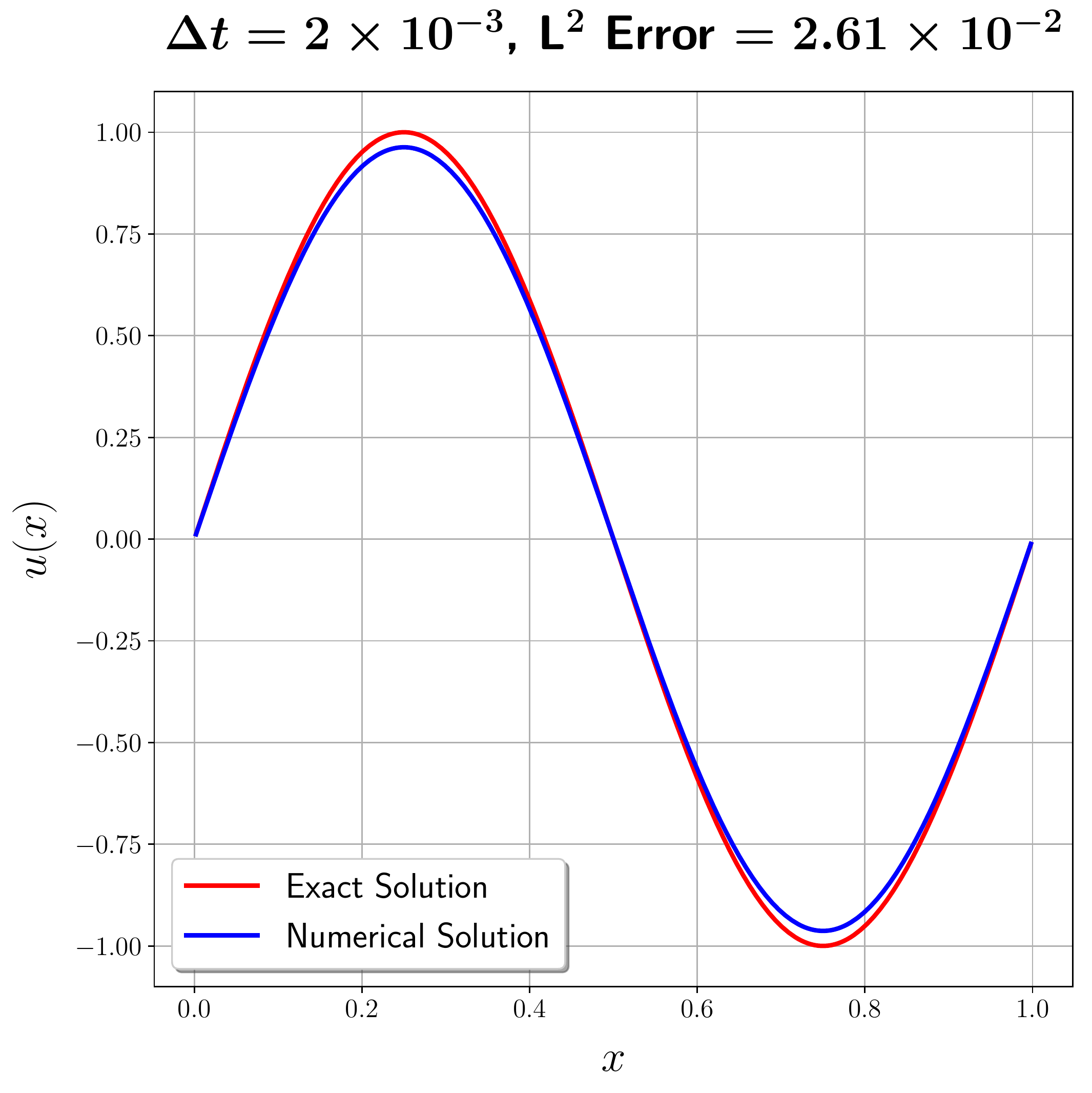} \hspace{0.5cm}
\includegraphics[scale=.35]{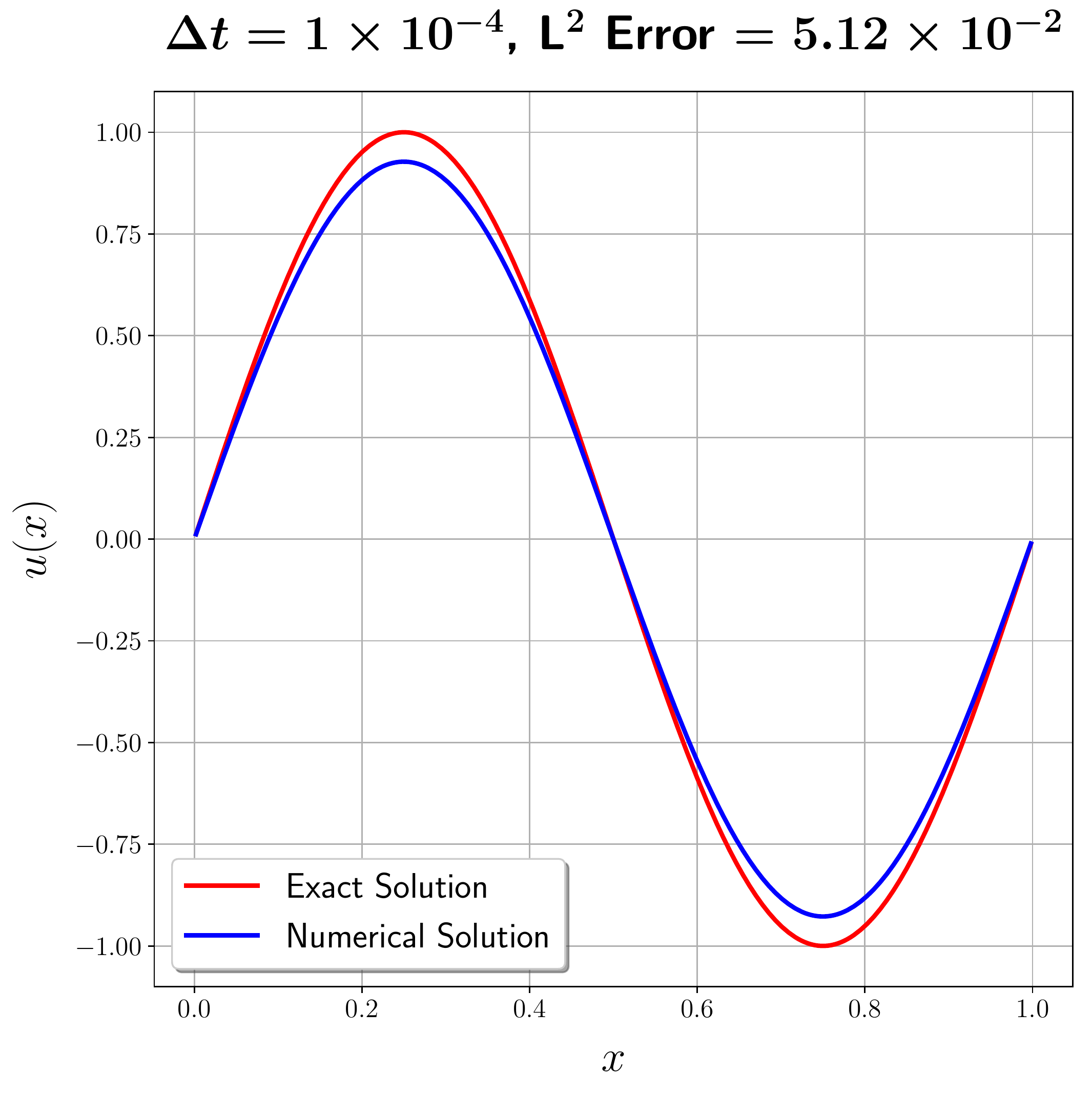}
\caption{The numerical solution of $u_t + u_x = 0$ with periodic boundary conditions at $t=1.0$ with $\Delta x = 1/2^8$ and two different time step sizes. As predicted by the theory, larger errors are incurred when a smaller time step size is used.} \label{LinearConstantCoefficientAdvection1D_ReductionInTimeStep}
\end{figure}

Figure~\ref{LinearConstantCoefficientAdvection1D_ReductionInTimeStep} shows the numerical solution of the linear advection equation~\eqref{LinearAdvecton1D_SimplestExample} on the domain $[0,1]$ with wave speed $a=1$, periodic boundary conditions, initial condition $u(x,0) = u_0(x) = \sin (2\pi x)$, and spatial resolution $\Delta x = 1/2^8$. The exact solution is $u(x,t) = u_0 \sin(2\pi(x-t))$. At $t=1.0$, we see that the error is larger with a time step $\Delta t = 10^{-4}$ when compared to the error with a 20 times larger time step. This is because the numerical diffusion, contributing to the error, is larger for the numerical solution using a smaller value of $\Delta t$, as evidenced by the higher reduction in the solution amplitude. 

To consider the effect of refinement only in space, we write the leading order term of the global truncation error as
\begin{equation}
\left[\left(\hat{\tau}_G\right)_j\right]_{\text{leading order}} = -\frac{1}{2} a^2 \Delta t \left(\frac{\Delta x}{|a| \Delta t} - 1\right) \left(u_{xx}\right)_j^n \equiv -\frac{1}{2} a^2 \Delta t \left(\frac{1}{C} - 1\right) \left(u_{xx}\right)_j^n.
\end{equation}
If $\Delta x$ is refined and $\Delta t$ held fixed so that $C < 1$ at all spatial resolutions, then $\left(\frac{1}{C} - 1\right)$ decreases towards 1, and the magnitude of the global truncation error, approximating the global solution error, decreases.

The unexpected behavior of the error norm with only temporal refinement can be attributed to the interaction of the leading order terms in the global truncation error. Since these terms have opposite signs, the magnitude of their difference increases with the reduction in $\Delta t$ at constant $\Delta x$ in the regime of $C \in (0,1)$. If, however, refinement were performed in both space and time by keeping $\Delta t$ proportional to $\Delta x$, the global truncation error would be dominated by the term (or the sum of the terms) with the lowest power of $\Delta x$ (or $\Delta t$) and only its magnitude, and not its sign, will play the pivotal role in the error. 

\subsection{Verification of the Spatial or Temporal Order of Accuracy}

We have established that asymptotic convergence may not be achieved with only spatial or temporal refinement. However, we can apply a technique to capture the order of the spatial and temporal discretizations. By considering only the leading order terms, the global solution error at a spatial location $x_j$ and a time horizon can be approximated as
\begin{equation}
\left(\hat{\tau}_G\right)_j \approx \left[\left(\hat{\tau}_G\right)_j\right]_{\text{leading order}} = \bigO\left({\Delta x}^{\alpha}\right) + \bigO\left({\Delta t}^{\beta}\right) = \zeta \Delta x^{\alpha} + {\zeta}_{\beta+1} \Delta t^{\beta}, \label{Global_Solution_Error_with_Leading_Order_Terms}
\end{equation}
where the coefficients $\zeta$ and ${\zeta}_{\beta+1}$ are independent of $\Delta x$ and $\Delta t$, and ${\zeta}_{\beta+1} = c_{\beta+1}/(\beta+1)!$ from~\eqref{LocalTruncationErrorNumericalSolutionFinalForm_1}. If ${\zeta} \Delta x^{\alpha} \gg {\zeta}_{\beta+1} \Delta t^{\beta}$, then the spatial order of convergence can be calculated by refining $\Delta x$ while keeping $\Delta t$ fixed. However, a general setting requires an alternative method to find the spatial order of convergence. This can be done by considering two uniform meshes with cell widths $\Delta x_i$ and $\Delta x_{i+1}$ with $\Delta x_{i+1} < \Delta x_i$. Then we can write
\begin{subequations}
\begin{align}
\left(\hat{\tau}_{G^x_i}\right)_j &\approx {\zeta} \Delta x_i^{\alpha} + {\zeta}_{\beta+1} \Delta t^{\beta}, \\
\left(\hat{\tau}_{G^x_{i+1}}\right)_j &\approx {\zeta} \Delta x_{i+1}^{\alpha} + {\zeta}_{\beta+1} \Delta t^{\beta}. 
\end{align}
\end{subequations}
Assuming $\left(\hat{\tau}_{G^x_{i+1}}\right)_j < \left(\hat{\tau}_{G^x_i}\right)_j$, we define
\begin{equation}
\Delta \left\{\left(\hat{\tau}_{G^x_{i,i+1}}\right)_j\right\} = \left(\hat{\tau}_{G^x_i}\right)_j - \left(\hat{\tau}_{G^x_{i+1}}\right)_j = {\zeta} \left(\Delta x_i^{\alpha} - \Delta x_{i+1}^{\alpha}\right) = {\zeta} \Delta x_{i+1}^{\alpha} \left\{\left(\frac{\Delta x_i}{\Delta x_{i+1}}\right)^{\alpha} - 1\right\} > 0.
\end{equation}
Defining $p = \Delta x_{i+1}/\Delta x_i < 1$ to be the ratio between the two mesh sizes, we can write
\begin{equation}
\Delta \left\{\left(\hat{\tau}_{G^x_{i,i+1}}\right)_j\right\} = {\zeta} \Delta x_{i+1}^{\alpha} \left(p^{-\alpha} - 1\right).
\end{equation}
Taking the logarithm of both sides,
\begin{equation}
\log \left[\Delta \left\{\left(\hat{\tau}_{G^x_{i,i+1}}\right)_j\right\}\right] = 
\theta + \alpha \log \left(\Delta x_{i+1}\right),
\end{equation}
where $\theta = \log \left\{{\zeta} \left(p^{-\alpha} - 1\right)\right\}$ is constant. So, we can compute the spatial order of accuracy by first choosing a sequence of grids with ${\Delta x}_{i+1}/{\Delta x}_i = p$ for $i=1,2,\ldots,M$, and all satisfying any CFL condition. Then, after interpolating the error to the coarsest mesh with spacing $\Delta x_1$, we can find the line of best fit of the norm of the difference between successive global solution errors $\Delta \left\{\left(\hat{\tau}_{G^x_{i,i+1}}\right)_{\text{norm}}\right\}$ vs.~the cell width~$\Delta x_{i+1}$ on a log-log scale for $i=1,2,\ldots,M-1$, and determine its slope which is the spatial order of accuracy.

Proceeding in a similar fashion, by refining only the time step by a constant ratio at a fixed spatial resolution, and plotting the norm of the difference between successive global solution errors, we can obtain the temporal order of accuracy. Since the spatial resolutions remain the same, we skip the interpolation step for refinement only in time.

Now, the exact solution is independent of the spatial resolution and the time step. So, for refinement only in time, if $u_j^{n+1}$ represents the solutions at time level $t^{n+1}$ and $\left(\hat{u}_{i}\right)_j^{n+1}$ represents its numerical counterpart obtained with time step~$\Delta t_i$, we can write
\begin{equation}
\left(\Delta \hat{\tau}_{G^t_{i,i+1}}\right)_j^{n+1} \equiv \left(\hat{\tau}_{G^t_i}\right)_j^{n+1} - \left(\hat{\tau}_{G^t_{i+1}}\right)_j^{n+1} = \left\{u_j^{n+1} - \left(\hat{u}_i\right)_j^{n+1}\right\} - \left\{u_j^{n+1} - \left(\hat{u}_{i+1}\right)_j^{n+1}\right\} = \left(\hat{u}_{i+1}\right)_j^{n+1} - \left(\hat{u}_i\right)_j^{n+1} \equiv \left(\Delta \hat{u}_{i,i+1}\right)_j^{n+1}.
\end{equation}
So, if we take the difference $\left(\Delta \hat{u}_{i,i+1}\right)_j^{n+1}$ between the numerical solutions $\left(\hat{u}_i\right)_j^{n+1}$ and $\left(\hat{u}_{i+1}\right)_j^{n+1}$ obtained with time steps $\Delta t_i$ and $\Delta t_{i+1}$ for $i = 1,2,\ldots,M-1$ at every mesh point $x_j$ and time level $t^{n+1}$, compute its norm $\left(\Delta \hat{u}_{G^t_{i,i+1}}\right)_{\text{norm}}^{n+1}$ and plot it against $\Delta t_{i+1}$, we will attain convergence with order equal to that of the time-stepping method. If we are performing only a spatial refinement, we first need to interpolate the numerical solution to the coarsest mesh and then follow the same steps to obtain convergence with the same order as that of the spatial discretization. From practical considerations, this approach has the clear advantage of not having to deal with an exact or manufactured solution.

It is worth keeping in mind that the sole purpose of these atypical convergence exercises is to verify the correct implementation of the spatial or temporal discretizations. The solution error norm under only spatial or temporal refinement is not expected to converge in the asymptotic regime. It is only when the time step and the cell width are refined simultaneously while keeping their ratio constant that we can expect convergence.

As an alternative, one can perform refinement in both $\Delta x$ and $\Delta t$, while maintaining $\Delta x^{\alpha} \propto \Delta t^{\beta}$, and plot the error norm (a) against $\Delta x$ to capture the spatial order of accuracy, and (b) against $\Delta t$ to capture the temporal order of accuracy. In this paper, we have not performed convergence studies with this refinement strategy, but we want to mention it for the sake of completeness. The first limitation of this refinement strategy is that for a high-order spatial and a low-order temporal discretization, keeping $\Delta t$ proportional to $\Delta x^{\alpha/\beta}$ can refine the time step to such an extent that the machine precision error dominates the discretization error. The second limitation is that one needs to know the order of the spatial and temporal discretizations i.e. the values of $\alpha$ and $\beta$ apriori. This knowledge is not necessary for plotting the differences in the numerical solution (or the error) for successive resolutions with refinement only in space (or only in time) to capture the spatial (or temporal) order of accuracy. As a result, we can even apply this technique to obtain the order of accuracy of complex spatial or temporal discretizations, even when it is difficult to extract the orders of accuracy analytically.

\subsection{Reduction in the Observed Order of Convergence}

Under certain circumstances, order reduction in the global solution is observed. Tables~\ref{Table_NatureOfConvergence_SpaceAndTime} and~\ref{Table_NatureOfConvergence_OnlySpace_OnlyTime} discuss the nature of convergence of the global solution error approximated as~\eqref{Global_Solution_Error_with_Leading_Order_Terms} before and after reaching the asymptotic regime. We consider the behavior when refining in space and time, and when refining only in space or time. As we reach the asymptotic regime, we can observe reduction in the order of convergence if
\begin{enumerate}[label=(\alph*),noitemsep] 
\item $\alpha > \beta$ and ${\zeta} \Delta x^{\alpha} \gg {\zeta}_{\beta+1} \Delta t^{\beta}$, or $\alpha < \beta$ and ${\zeta} \Delta x^{\alpha} \ll {\zeta}_{\beta+1} \Delta t^{\beta}$, for refinement in both space and time while keeping $\Delta t$ proportional to $\Delta x$;
\item ${\zeta} \Delta x^{\alpha} \gg {\zeta}_{\beta+1} \Delta t^{\beta}$ for refinement only in space;
\item ${\zeta} \Delta x^{\alpha} \ll {\zeta}_{\beta+1} \Delta t^{\beta}$ for refinement only in time.
\end{enumerate}
In Section \ref{sec:convergence_plots}, we will encounter order reduction with the convergence plots of a linear variable-coefficient advection equation and a non-linear advection equation, with refinement in both space and time, when they are (a) discretized in space with a non-monotone finite volume method, and advanced in time with the explicit midpoint and the second-order Adams-Bashforth methods, and (b) discretized in space with a monotone finite volume method, and advanced in time with Forward Euler, the explicit midpoint and the second-order Adams-Bashforth methods.

\begin{table}[!t]
\centering
\caption{Convergence of the global solution error approximated as \eqref{Global_Solution_Error_with_Leading_Order_Terms} before and after reaching the asymptotic regime, while maintaining $\Delta t/\Delta x = r$.}
\vspace{3mm}
\setlength{\tabcolsep}{0.25em}
\begin{tabular}{ccc}
\toprule
& $\alpha > \beta$ & $\alpha < \beta$ \\
{\colorbox{shade_0}
{\parbox{2.8cm}{${\zeta} \Delta x^{\alpha} \gg {\zeta}_{\beta+1} \Delta t^{\beta} \\ \text{i.e. } {\zeta} \Delta x^{\alpha-\beta} \gg r {\zeta}_{\beta+1}$}}} &
{\colorbox{shade_1}
{\parbox{5.9cm}{Before asymptotic regime:
$\tau_G \approx {\zeta} \Delta x^{\alpha}$ \\
Convergence: Attained with slope $\alpha$ \\
After asymptotic regime:
$\tau_G \approx r {\zeta}_{\beta+1} \Delta x^{\beta}$ \\
Convergence: Attained with slope $\beta$ \\
Order reduction: Observed}}} &
{\colorbox{shade_1}
{\parbox{5.9cm}{Before asymptotic regime:
$\tau_G \approx {\zeta} \Delta x^{\alpha}$ \\
Convergence: Attained with slope $\alpha$ \\
After asymptotic regime:
$\tau_G \approx {\zeta} \Delta x^{\alpha}$ \\
Convergence: Attained with slope $\alpha$ \\
Order reduction: Not observed}}} \vspace{1mm} \\
{\colorbox{shade_0}
{\parbox{2.8cm}{${\zeta} \Delta x^{\alpha} \ll {\zeta}_{\beta+1} \Delta t^{\beta} \\ \text{i.e. } {\zeta} \Delta x^{\alpha-\beta} \ll r {\zeta}_{\beta+1}$}}} &
{\colorbox{shade_1}
{\parbox{5.9cm}{Before asymptotic regime:
$\tau_G \approx r {\zeta}_{\beta+1} \Delta x^{\beta}$ \\
Convergence: Attained with slope $\beta$ \\
After asymptotic regime:
$\tau_G \approx r {\zeta}_{\beta+1} \Delta x^{\beta}$ \\
Convergence: Attained with slope $\beta$ \\
Order reduction: Not observed}}} &
{\colorbox{shade_1}
{\parbox{5.9cm}{Before asymptotic regime:
$\tau_G \approx r {\zeta}_{\beta+1} \Delta x^{\beta}$ \\
Convergence: Attained with slope $\beta$ \\
After asymptotic regime:
$\tau_G \approx {\zeta} \Delta x^{\alpha}$ \\
Convergence: Attained with slope $\alpha$ \\
Order reduction: Observed}}} \vspace{1mm} \\
\bottomrule \\
\end{tabular} \label{Table_NatureOfConvergence_SpaceAndTime}
\end{table}

\begin{table}[!t]
\centering
\caption{Convergence of the global solution error approximated as \eqref{Global_Solution_Error_with_Leading_Order_Terms} before and after reaching the asymptotic regime with refinement only in space or only in time.}
\vspace{3mm}
\setlength{\tabcolsep}{0.25em}
\begin{tabular}{ccc}
\toprule
& {Refinement only in space} & {Refinement only in time} \\
{\colorbox{shade_0}
{\parbox{2.3cm}{\centering ${\zeta} \Delta x^{\alpha} \gg {\zeta}_{\beta+1} \Delta t^{\beta}$}}} &
{\colorbox{shade_1}
{\parbox{5.65cm}{Before asymptotic regime:
$\tau_G \approx {\zeta} \Delta x^{\alpha}$ \\
Convergence: Attained with slope $\alpha$ \\
After asymptotic regime:
$\tau_G \approx {\zeta}_{\beta+1} \Delta t^{\beta}$ \\
Convergence: Not attained \\
Order reduction: Observed}}} &
{\colorbox{shade_1}
{\parbox{5.65cm}{Before asymptotic regime:
$\tau_G \approx {\zeta} \Delta x^{\alpha}$ \\
Convergence: Not attained \\
After asymptotic regime:
$\tau_G \approx {\zeta} \Delta x^{\alpha}$ \\
Convergence: Not attained \\
Order reduction: Not applicable}}} \vspace{1mm} \\
{\colorbox{shade_0}
{\parbox{2.3cm}{\centering ${\zeta} \Delta x^{\alpha} \ll {\zeta}_{\beta+1} \Delta t^{\beta}$}}} &
{\colorbox{shade_1}
{\parbox{5.65cm}{Before asymptotic regime:
$\tau_G \approx {\zeta}_{\beta+1} \Delta t^{\beta}$ \\
Convergence: Not attained \\
After asymptotic regime:
$\tau_G \approx {\zeta}_{\beta+1} \Delta t^{\beta}$ \\
Convergence: Not attained \\
Order reduction: Not applicable}}} &
{\colorbox{shade_1}
{\parbox{5.65cm}{Before asymptotic regime:
$\tau_G \approx {\zeta}_{\beta+1} \Delta t^{\beta}$ \\
Convergence: Attained with slope $\beta$ \\
After asymptotic regime:
$\tau_G \approx {\zeta} \Delta x^{\alpha}$ \\
Convergence: Not attained \\
Order reduction: Observed}}} \vspace{1mm} \\
\bottomrule \\
\end{tabular} \label{Table_NatureOfConvergence_OnlySpace_OnlyTime}
\end{table}

\subsection{Local Truncation Error of a Linear Inhomogeneous Variable-Coefficient Advection Equation} \label{LTE_Linear_Inhomogeneous_Variable_Coefficient_Advection_Equation}

We consider the linear variable-coefficient one-dimensional inhomogeneous advection equation
\begin{equation}
u_t + p(x) u + (q(x) u)_x = f(x,t),
\label{LinearAdvection1D_1}
\end{equation}
which can be expressed as 
\begin{equation}
u_t + F_x = s \equiv -p(x) u + f(x,t),
\label{LinearAdvection1D}
\end{equation}
where $F = q(x) u$ is the flux with $q(x) > 0$, and the source term $s(u,x,t) = -p(x) u + f(x,t)$ consists of two parts: a linear variable-coefficient function of the dependent variable $-p(x) u$ and a function of the independent variables $f(x,t)$. The term $-p(x) u$ is motivated by the Coriolis acceleration $f \hat{k} \times \vec{u}$ appearing in the horizontal momentum equations of geophysical flows. Here $f$ denotes the Coriolis parameter, which may be a constant on the idealized f-plane, or linear in latitude on the beta-plane, which is an example of the variable-coefficient case in equation \eqref{LinearAdvection1D}. 

Leveraging the computational power of SymPy, a symbolic package of Python, we calculate the first few relevant terms containing ${\Delta t}^l {\Delta x}^k$ for $l = 1,2,\ldots$ and $k=0,1,\ldots$ of the local truncation error of~\eqref{LinearAdvection1D} for various spatial and temporal discretizations. Tables~B.7--B.11b list these terms for the first-order upwind finite difference spatial discretization, and the five explicit time-stepping methods of List~\ref{myListOfExplicitTimeSteppingMethodsForAnalysis}. The supplementary text file `LocalTruncationError\_Output.rtf' contains these results for second- and third-order upwind finite difference spatial discretizations. Determining symbolic representations of the local truncation error of~\eqref{LinearAdvection1D} with SymPy consists of a few steps. We start by using the spatial and temporal discertizations to find expressions of $\hat{u}_j^{n+1}$ as functions of quantities defined at spatial locations adjacent to and including $x_j$ and temporal locations adjacent to and including~$t^n$. Next, Taylor expansions are used to expand every term about $(x_j,t^n)$. The third step requires expressing the temporal and mixed derivatives of $u$ as functions of quantities at the current time level $t^n$, which are assumed to be known a priori. For the advection equation~\eqref{LinearAdvection1D}, these derivatives are in Tables~B.3--B.6. In the final step, we compute the difference between the exact solution $u_j^{n+1}$ and its numerical counterpart~$\hat{u}_j^{n+1}$ to arrive at the final form of the local truncation error $\hat{\tau}_j^{n+1}$. For all spatial and temporal discretizations we consider, $\hat{\tau}_j^{n+1}$ can be expressed as
\begin{equation}
\hat{\tau}_j^{n+1} = \sum \limits_{k=1}^{\infty} \frac{\Delta t^k}{k!} \left(c_k + \mathcal{O} \left(\Delta x^{\alpha}\right)\right)_j^n, \label{LocalTruncationErrorNumericalSolutionFinalForm_Compact}
\end{equation}
where $\alpha$ and $\beta$ represent the orders of the spatial and temporal discretizations, and $c_k = 0$ for $k=1,2,\ldots,\beta$, and this expression is a compact form of~\eqref{LocalTruncationErrorNumericalSolutionFinalForm}.

Now we consider the special situation when $p(x)$ is constant, $p(x) = p_0$, and $q(x)$ is linear, $q(x) = q_0 + q_1 x$, so that $q_x(x) = q_1$, and $u(x,t)$ and $f(x,t)$ are functions of only $t$. Then,
\begin{enumerate}[label=(\alph*),noitemsep] 
\item the linear advection equation~\eqref{LinearAdvection1D} reduces to the linear ODE~\eqref{ODE1DParticularChoice};
\item the coefficients of $\Delta x^k$ for $k \ge \alpha$ in the $\bigO\left(\Delta x^{\alpha}\right)$ terms reduce to zero;
\item the local truncation error of~\eqref{LinearAdvection1D}, assuming form~\eqref{LocalTruncationErrorNumericalSolutionFinalForm_Compact}, reduces to that of~\eqref{ODE1DParticularChoice}, assuming form~\eqref{ODE1DTruncationErrorNumericalSolution_2} given by
\begin{equation}
\hat{\tau}^{n+1} = \sum \limits_{k=\beta+1}^{\infty} \frac{c_k^n}{k!} \Delta t^k = \frac{c_{\beta+1}}{(\beta+1)!} \Delta t^{\beta+1} + \mathcal{O} \left(\Delta t^{\beta+2}\right).
\end{equation}
\end{enumerate}

Tables B.12--B.16 contain the coefficients of $\Delta t^l \Delta x^k$ for $l=1,2,\ldots,\beta+1$ and $k = 0,1,2,3$ in the local truncation error of the linear advection equation~\eqref{LinearAdvection1D} discretized in space with the first-order upwind finite difference scheme ($\alpha=1$), and advanced in time with the five explicit time-stepping methods of List~\ref{myListOfExplicitTimeSteppingMethodsForAnalysis} i.e.~the first-order Forward Euler method ($\beta=1$), the second-order explicit midpoint method ($\beta=2$), the second-order Adams-Bashforth method ($\beta=2$), Williamson's low-storage third-order Runge-Kutta method~\citep{williamson1980low} ($\beta=3$), and the third-order Adams-Bashforth method ($\beta=3$). Since $c_{\beta+1}/(\beta+1)!$ is the coefficient of $\Delta t^{\beta+1} \Delta x^0$ in the local truncation error, its explicit expression is present in the row with $l=\beta+1$ and $k=0$. With the assumption
\begin{equation}
p_x(x) = q_{xx}(x) = u_x(x,t) = f_x(x,t) = 0,
\end{equation} 
the linear advection equation~\eqref{LinearAdvection1D} reduces to the linear ODE~\eqref{ODE1DParticularChoice}, and the above-mentioned coefficient of $\Delta t^{\beta+1} \Delta x^0$ reduces to that of $\Delta t^{\beta+1}$ in the local truncation error of~\eqref{ODE1DParticularChoice} advanced with the same time-stepping method and listed in Table~B.2. One can also verify that the local truncation error of the generic hyperbolic PDE~\eqref{AdvectionEquation1DFunctionalForm} of Section~\ref{sec:pdes} using any of the time-stepping methods of List~\ref{myListOfExplicitTimeSteppingMethodsForAnalysis} reduces to that of the generic ODE~\eqref{ODE1D} of Section~\ref{sec:odes} with approximations~\eqref{myApproximationsForPDEsUsedToModelPhysicalPhenomena_1} and~\eqref{myApproximationsForPDEsUsedToModelPhysicalPhenomena_3}, advanced with the same time-stepping method. Finally, expressing the particular ODE~\eqref{ODE1DParticularChoice} as 
\begin{equation}
u_t = \mathcal{F}(u,t) \equiv -\left(p_0 + q_1\right) u + f(t), \label{ODE1DParticularChoice_ExpressedInGenericForm}
\end{equation} 
and the particular PDE~\eqref{LinearAdvection1D} as
\begin{equation}
u_t = \mathcal{F} \left(u,u_x,x,t\right) \equiv -p(x) u - (q(x) u)_x + f(x,t), \label{LinearAdvection1D_ExpressedInGenericForm}
\end{equation} 
the local truncation error of the generic ODE~\eqref{ODE1D} and the generic hyperbolic PDE~\eqref{AdvectionEquation1DFunctionalForm} advanced with any of the time-stepping methods of List~\ref{myListOfExplicitTimeSteppingMethodsForAnalysis} reduce to that of the particular ODE~\eqref{ODE1DParticularChoice} and the particular advection equation~\eqref{LinearAdvection1D} with the specific formulation of $\mathcal{F}$ given by~\eqref{ODE1DParticularChoice_ExpressedInGenericForm} and~\eqref{LinearAdvection1D_ExpressedInGenericForm}, respectively.

\subsection{Local Truncation Error of a Non-Linear Inhomogeneous Advection Equation}

We conclude our analysis by considering the non-linear advection equation
\begin{equation}
\left(\bar{u} + u\right)_t + \left(\bar{u} + u\right) \left(\bar{u} + u\right)_x + p_0 \left(\bar{u} + u\right) = \hat{f}(x,t) \label{NonLinearAdvection1D_1},
\end{equation} 
which can be expressed in conservative form as
\begin{equation}
u_t + \left(\bar{u} u + \frac{u^2}{2}\right)_x + p_0 u = f(x,t), \label{NonLinearAdvection1D}
\end{equation} 
where $p_0$ is a constant and $f(x,t) = \hat{f}(x,t) - \bar{u} \bar{u}_x - p_0 \bar{u}$. Motivated by applications in fluid dynamics, $u$ has been decomposed into a constant mean component $\bar{u}$, and a perturbation term $u$, which is a function of space and time. If $\bar{u}$, $p_0$ and $f(x,t)$ are reduced to zero,~\eqref{NonLinearAdvection1D} reduces to the inviscid Burgers' equation. The supplementary text file `LocalTruncationError\_Output.rtf' contains the leading order terms of the local truncation error of~\eqref{NonLinearAdvection1D} discretized in space with the first-, second-, and third-order upwind finite difference schemes and advanced in time with the five explicit time-stepping methods of List~\ref{myListOfExplicitTimeSteppingMethodsForAnalysis}.

Similar to our reasoning in Section \ref{LTE_First_Order_Linear_ODE}, we cannot employ implicit time-stepping methods to advance~\eqref{LinearAdvection1D_1} and~\eqref{NonLinearAdvection1D_1}, since it would require knowledge of the functional forms of $p(x)$, $q(x)$, and $f(x,t)$ for~\eqref{LinearAdvection1D_1}, and $\bar{u}(x)$ and $f(x,t)$ for~\eqref{NonLinearAdvection1D_1}.

\section{Numerical Results} \label{sec:numerical_results}

In this section, we numerically verify our theoretical findings for the spatial and temporal order of convergence of hyperbolic PDEs. We perform convergence studies on the linear variable-coefficient inhomogeneous advection equation
\begin{equation}
u_t + x u_x + 2 u \equiv u_t + (x u)_x + u = s, \quad x \in [0,1], \: t>0, \label{LinearAdvection1D_NumericalExperiment}
\end{equation}
and the nonlinear inhomogeneous advection equation
\begin{equation}
(1 + u)_t + (1 + u) (1 + u)_x + (1 + u) \equiv u_t + \left(u + \frac{1}{2} u^2\right)_x + 1 = s, \quad x \in [0,1], \: t>0, \label{NonLinearAdvection1D_NumericalExperiment}
\end{equation}
with periodic boundary conditions. The linear advection equation~\eqref{LinearAdvection1D_NumericalExperiment} is a special case of~\eqref{LinearAdvection1D} with $p(x) = 1$ and $q(x) = x$, while the non-linear advection equation~\eqref{NonLinearAdvection1D_NumericalExperiment} is a special case of~\eqref{NonLinearAdvection1D} with $\bar{u}=p_0=1$. The exact solution is chosen to be
\begin{equation}
u_{\text{exact}}(x,t) = \hat{u} \sin (kx - \omega t) + 2 \hat{u} \cos (2kx - \omega t), \label{ExactSolution_NumericalExperiments}
\end{equation} 
which is a superposition of two sinusoidal wave modes, with $\hat{u}$, $k$, and $\omega$ representing the amplitude, wavenumber, and angular velocity of the first wave mode. The second wave mode has twice the amplitude, half the wavelength, and half the phase speed as the first one, and leads in phase by 90 degrees. If $c$ denotes the phase speed of the first wave mode, we can write $\omega = ck = \frac{c}{2} (2k)$, so that the angular velocity remains the same for both wave modes. We specify $k=2\pi$, $c=1$, and $\hat{u} = 1$ for the linear advection equation and $\hat{u} = 0.01$ for the non-linear advection equation. By substituting $t=0$ in \eqref{ExactSolution_NumericalExperiments}, we obtain the initial condition
\begin{equation}
u_{\text{exact}}(x,0) = \hat{u} \sin (kx) + 2 \hat{u} \cos (2kx). \label{Initialondition_NumericalExperiments}
\end{equation}
The motivation behind the choice \eqref{ExactSolution_NumericalExperiments} for the exact solution is to eliminate artificially high rates of convergence that sometimes occurs with `nice' test problems such as a single sinusoidal wave mode whose leading order terms in the local truncation error can be zero. By substituting the exact solution~\eqref{ExactSolution_NumericalExperiments} into the left-hand side of~\eqref{LinearAdvection1D_NumericalExperiment} and~\eqref{NonLinearAdvection1D_NumericalExperiment}, we obtain the corresponding source terms on the right-hand side. We employ finite difference and finite volume methods for spatial discretization and the following set of predictor-corrector and multistep time-stepping methods, ranging from first- to fourth-order, to advance our numerical solution in time:
\begin{mylist}\mbox{Time-stepping methods for numerical experiments:}
\begin{enumerate}[noitemsep] 
\label{myListOfTimeSteppingMethodsForNumericalExperiments}
\item[FE1:] first-order Forward Euler method
\item[RK2:] explicit midpoint method, belonging to the second-order Runge-Kutta family
\item[RK3:] low-storage third-order Runge-Kutta method of~\citet{williamson1980low}
\item[RK4:] low-storage five-stage fourth-order Runge-Kutta method of~\citet{carpenter1994fourth}
\item[AB2:] second-order Adams-Bashforth method
\item[AB3:] third-order Adams-Bashforth method
\item[AB4:] fourth-order Adams-Bashforth method
\end{enumerate}
\end{mylist}

\subsection{Spatial Discretization} \label{SpatialDiscretization}
We consider two spatial discretization methods: a standard first-order finite difference upwind scheme, and a piecewise parabolic reconstruction (PPR) finite volume scheme that is equivalent to the spatial discretization part of the piecewise parabolic method (PPM) of~\citet{colella1984piecewise}. The PPR involves fitting a parabolic profile within each cell. The three constants needed to uniquely define this parabola are determined by solving a linear system for the cell-averaged solution and the left and right edge estimates. Finalizing the values of these edge estimates consists of a few steps. Starting with the mean solution within each cell, PPR first interpolates the solution to the edges. This interpolation is fourth-order accurate on a uniform mesh. Then PPR applies the monotonized-central slope limiter and adjusts the edge estimates to flatten any local maximum or minimum with the cell. A close variant of the original PPR scheme appears in~\citet{engwirda2016weno}, where the application of the slope limiter and adjustments of the edge estimates are performed in the reverse order. We have used both versions of PPR, and obtained similar results. In this paper, however, we present our results with the original version. 

The application of the slope limiter and the edge adjustments guarantee that the parabolic profile within each cell is oscillation-free, monotonicity-preserving, and total variation diminishing. Even though these monotone slope-limiting strategies ensure numerical stability, it comes at the cost of accuracy. More specifically, the combined effect of the slope limiter and the flattening of new local extrema manifests as spurious numerical dissipation and results in discontinuities at cell edges, thereby compromising both the spatial and temporal orders of approximation. The numerical flux at every edge is typically a function of these edge estimates among other parameters, and is inevitably different from the flux computed with the edge estimate from the first interpolation step, which would be fourth-order accurate on a uniform mesh. In the original PPM~\cite{colella1984piecewise}, the numerical flux is computed as the integral average of the flux passing through each edge from the current time step to the next one. Using this time-centered approximation of the numerical fluxes, the time-centered tendencies are computed and used to determine the solution at the next time step. Instead of adopting this approach to advance our numerical solution in time, we combine PPR with the seven time-stepping methods of List~\ref{myListOfTimeSteppingMethodsForNumericalExperiments}. However, we still need to compute numerical fluxes at a previous time level for a multistep time-stepping method, or at a fraction of a time step for a predictor-corrector time-stepping method. This is where we employ the slightly dissipative local Lax-Friedrichs Riemann solver of~\citet{rusanov1961calculation}, which determines the dot product of the numerical flux vector in the direction of the flow at every edge with the outward unit normal vector as
\begin{equation}
\mathbb{F}^*\left(u^{\text{int}},u^{\text{ext}};\hat{n}\right) = \vec{\mathbb{F}} \cdot \hat{n} = \frac{1}{2} \left\{\left(F^{\text{int}} + F^{\text{ext}}\right) \text{sign}\left(\hat{n}\right) - \left|\lambda\right|_{\max} \left(u^{\text{ext}} - u^{\text{int}}\right)\right\}.
\end{equation}
Here $F^{\text{int}}$ and $F^{\text{ext}}$ are the fluxes of $u^{\text{int}}$ and $u^{\text{ext}}$, the internal and external states at the edge of the cell, where the solution tendency is being computed, and $\hat{n}$ is the unit normal vector directed from the internal state to the external one. For example, to determine the numerical flux at the right edge of the cell $[x_{j-\frac{1}{2}},x_{j+\frac{1}{2}}]$, we specify $u^{\text{int}} = u^R_j$, $u^{\text{ext}} = u^L_{j+1}$, and $\hat{n} = \hat{x}$, where the superscripts $L$ and $R$ represent the left and right edge estimates. The term $\left|\lambda\right|_{\max} = \max(F^{\text{int}}_u,F^{\text{ext}}_u)$ is the larger of the magnitudes of the two wave speeds, the first one computed as a function of the internal state $u^{\text{int}}$, and the second one as a function of the external state $u^{\text{ext}}$. We use this formulation of the numerical flux in our experiments modeling both linear and non-linear advection. We evaluate the gradient of the flux within the tendency of the cell-averaged solution $\bar{u}_j$ as 
\begin{equation}
\left[F_x\right]_j = \frac{1}{2 \Delta x}\left\{\left.\mathbb{F}^*\left(u^{\text{int}},u^{\text{ext}};\hat{n}\right)\right|_{j-\frac{1}{2}} + \left.\mathbb{F}^*\left(u^{\text{int}},u^{\text{ext}};\hat{n}\right)\right|_{j+\frac{1}{2}}\right\}.
\end{equation}
Equipped with the flux gradient and the source terms, we compute the solution tendency and advance the solution to the next time step. 

\subsection{Computing the Error Norm}

In addition to the slope limiter, the monotonicity-preserving strategies, and the (dissipative) Riemann solver, the order of convergence of a hyperbolic PDE with refinement in both space and time also depends on whether or not the prognostic variable is chosen to be cell-integrated or cell-averaged for a finite volume method, and whether the numerical solution (or the error) is interpolated to the coarsest mesh for a finite difference method.

For refinement only in time, we use a mesh with the same spatial resolution. When we refine only in space, we need to interpolate the numerical solution (or the error) to the coarsest mesh so that we can compute the difference between the numerical solution (or the error) for successive pairs of spatial resolutions. However, when we perform a refinement in space and time simultaneously, it is not immediately clear if we need to perform the above-mentioned interpolation to the coarsest mesh.

\subsubsection{Interpolation to the Coarsest Mesh for a Finite Difference Method}

If $\alpha$ and $\beta$ denote the spatial and temporal orders of accuracy of a PDE, we know that the coefficients of $\Delta x^k$ for $k = \alpha, \alpha+1, \ldots$ in the coefficients of $\Delta t^l$ for $l = 1,2,\ldots$ within the local truncation error of a PDE are functions of the spatial gradients of the dependent variable, the coefficients of the PDE, and the source terms at the current time. We now consider the error of the numerical solution computed at a set of grid points, as in a finite difference method. If we compute the error norm over the entire spatial interval at a certain time horizon using the magnitudes of the error at every grid point, we may obtain a higher error norm for a fine mesh than a coarse one. The reason for this discrepancy is that the decrease in the magnitude of the local truncation error at a set of points on the fine mesh due to the reduction in the magnitude of $\Delta x$ may be off-set by the increase in the magnitude of the above-mentioned spatial gradients at a subset of these points, which may not even exist on the coarse mesh. This is more pronounced in convergence studies for higher-dimensional problems involving unstructured meshes where a fine mesh is not necessarily embedded within a coarse one. With the global solution error being approximated by the global truncation error, which in turn is one order of $\Delta t$ less than the local truncation error, we may not even obtain numerical convergence by computing the error norm based on the magnitudes of the error obtained at the native set of points within each mesh. Even if we achieve convergence, the order may be less than the expected one. Therefore, for refinement in both space and time using a finite difference method, it is advisable to interpolate the error to the set of coarsest mesh points and then determine the error norm.

\subsubsection{Cell-Integrated vs.~Cell-Averaged Quantity as Prognostic Variable for a Finite Volume Method}

We consider the formulation of any finite volume method for solving the one-dimensional inhomogeneous advection equation, expressed in conservative form
\begin{equation}
u_t + F_x = s, \label{FiniteVolume1DConservativeForm}
\end{equation}
where $u$ is the scalar quantity being advected, and the flux $F$ can be a linear or non-linear function of $u$ with constant or spatially dependent coefficients. The source term $s$ can be a function of $x$, $t$ and $u$, but not of $u_x$. Integrating~\eqref{FiniteVolume1DConservativeForm} with respect to $x$ over the cell~$j$ with,
\begin{equation}
U_t = F_{j-\frac{1}{2}} - F_{j+\frac{1}{2}} + S, \label{FiniteVolume1DConservativeFormIntegratedOverCell}
\end{equation}
where
\begin{equation}
U = \int \limits_{x_{j-\frac{1}{2}}}^{x_{j+\frac{1}{2}}} u dx, \quad \text{and} \quad S = \int \limits_{x_{j-\frac{1}{2}}}^{x_{j+\frac{1}{2}}} s dx.
\end{equation}
Therefore the spatial order of accuracy of the right-hand side of~\eqref{FiniteVolume1DConservativeFormIntegratedOverCell} depends on how the flux terms are constructed. In the PPR scheme, before the application of the slope limiter and the monotonicity-preserving strategies, this order of accuracy is 4. 

If $m$ is the interpolant's order of accuracy, and we define the prognostic variable to be the cell-integrated solution $U$, then the spatial order of accuracy of the finite volume method is $m$. If, however, we define our prognostic variable to be the cell-averaged solution 
\begin{align}
\bar{u} = \frac{U}{\Delta x} = \frac{1}{\Delta x} \int \limits_{x_{j-\frac{1}{2}}}^{x_{j+\frac{1}{2}}} u dx,
\end{align}
which is the standard practice with finite volume methods, our prognostic equation becomes
\begin{equation}
\bar{u}_t = \frac{F_{j-\frac{1}{2}} - F_{j+\frac{1}{2}}}{\Delta x} + \bar{s}, \quad
\text{where} \quad
\bar{s} = \frac{S}{\Delta x} = \frac{1}{\Delta x} \int \limits_{x_{j-\frac{1}{2}}}^{x_{j+\frac{1}{2}}} s dx,
\label{FiniteVolume1DConservativeFormAveragedOverCell}
\end{equation}
and the spatial order of accuracy of the finite volume method drops to $m-1$. Based on Theorem~\ref{Theorem_1}, the order of convergence of $U$ and $\bar{u}$ at constant ratio of time step to cell width will be $\min(m,n)$ and $\min(m-1,n)$, respectively, where $n$ is the order of the time-stepping method. 

\subsubsection{Interpolation to the Coarsest Mesh for a Finite Volume Method}

Since a finite volume method advances the cell-integrated or the cell-averaged solution in time, it considers the entire variation of the solution within the cells. Therefore it may not be necessary to interpolate the finite volume solution (or its error) to the coarsest mesh for simultaneous refinement in space and time, unlike a finite difference solution (or its error) defined only at a set of grid points. However, if the finer meshes are not embedded within the coarsest mesh, or if we are refining only in space, we need to perform this interpolation.

\begin{figure}[!htp]
\centering
\includegraphics[scale=.305]{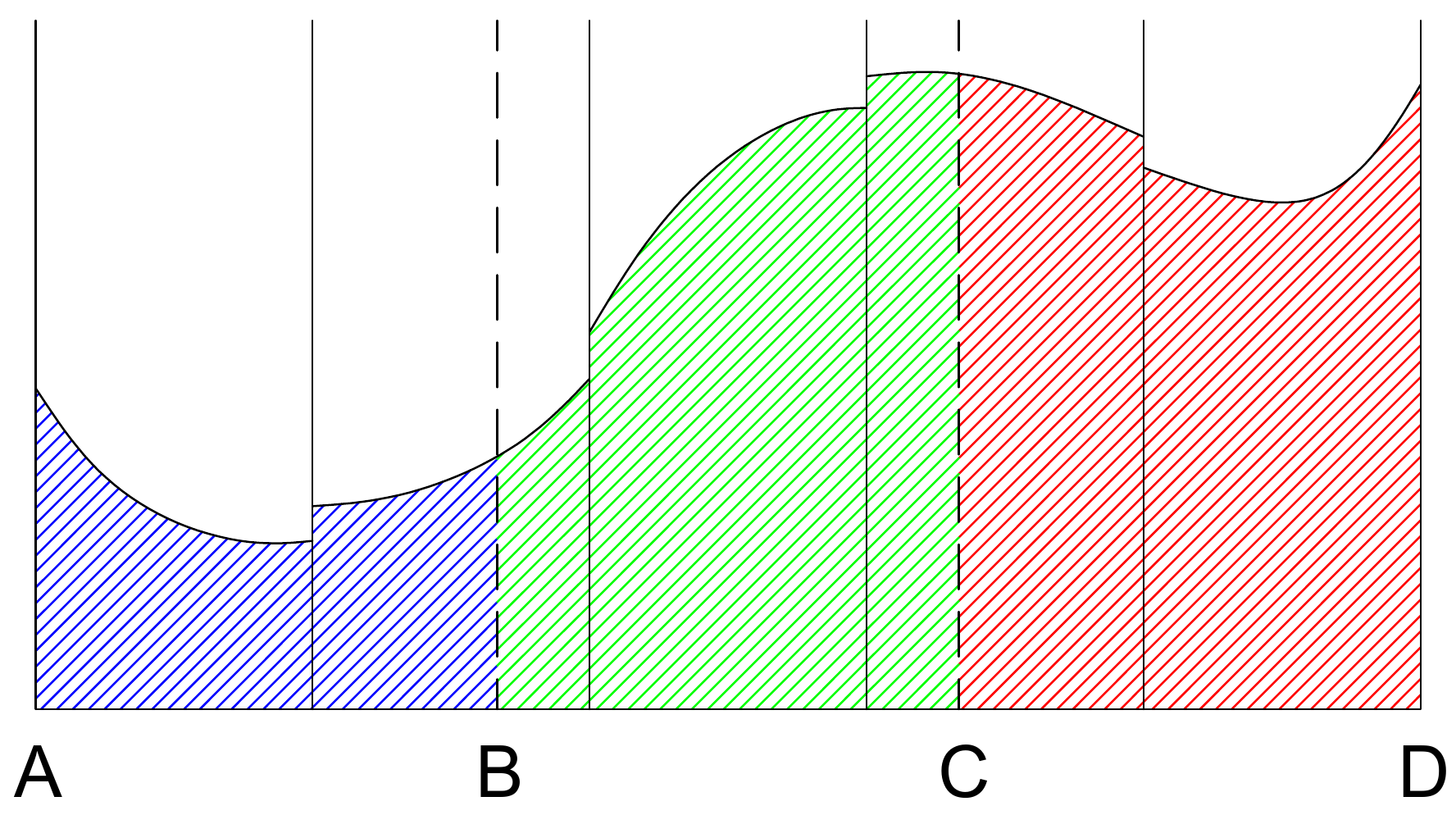} 
\caption{A schematic of the cell-integrated solution of a fine mesh consisting of five cells (between solid lines) being interpolated to the coarsest mesh consisting of three cells (shown by fill colors), for a finite volume method.}
\label{Figure_PPR_InterpolationToCoarsestMesh}
\end{figure}

By integrating the parabolic profiles of the solution within the cells of a fine mesh which are either partially or entirely contained within the cells of the coarsest mesh, we can interpolate the cell-integrated solution to the cells of the coarsest mesh. Let $U^{f \to c}$ denote the cell-integrated solution and $E^{f \to c}$ the cell-integrated error interpolated to the coarsest mesh. Figure \ref{Figure_PPR_InterpolationToCoarsestMesh} illustrates an example where 5 cells of the fine mesh are contained within 3 cells of the coarsest mesh. The cell-integrated solution of the fine mesh interpolated to cells $AB$, $BC$, and $CD$ of the coarsest mesh are the blue, green, and red shaded areas, respectively. Now, the cell-averaged solution and error interpolated to the coarsest mesh are $\bar{u}^{f \to c} = U^{f \to c}/\Delta x_c$ and $\bar{e}^{f \to c} = E^{f \to c}/\Delta x_c$, respectively, where $\Delta x_c$ is the cell-width of the coarsest mesh. Since $\Delta x_c$ is constant, a log-log plot of $U^{f \to c}$ $\left(\text{or }E^{f \to c}\right)$ versus $\Delta x_f$ will have the same slope but a different intercept as that of $\bar{u}^{f \to c}$ or $\left(\text{or }\bar{e}^{f \to c}\right)$ versus $\Delta x_f$, where $\Delta x_f$ represents the cell width of the fine meshes. 




\subsection{Convergence Plots} \label{sec:convergence_plots}

\begin{figure}[!htp]
\centering
\hspace{0.165cm}
\includegraphics[scale=.3125]{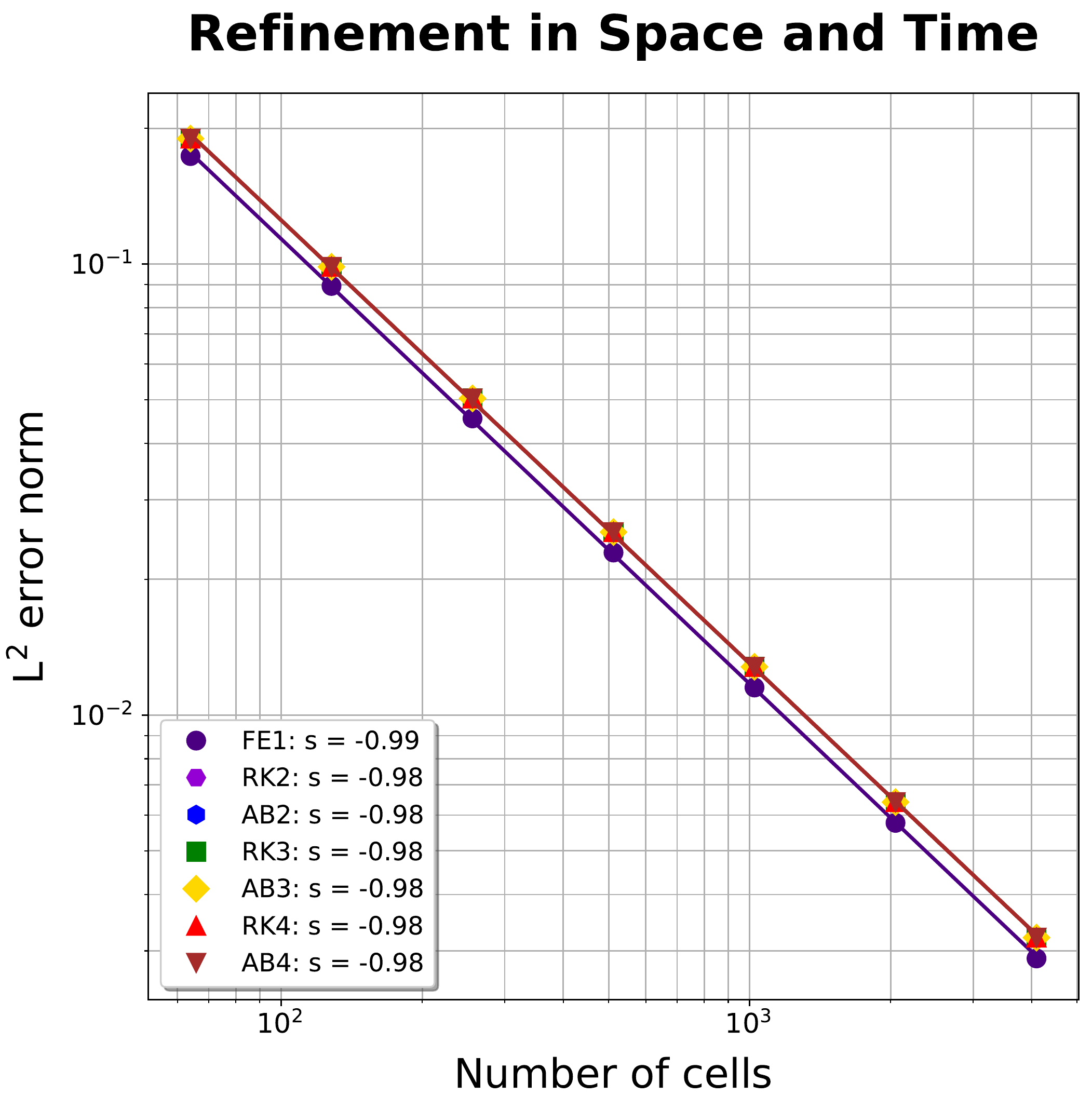} \hspace{0.785cm}
\includegraphics[scale=.3125]{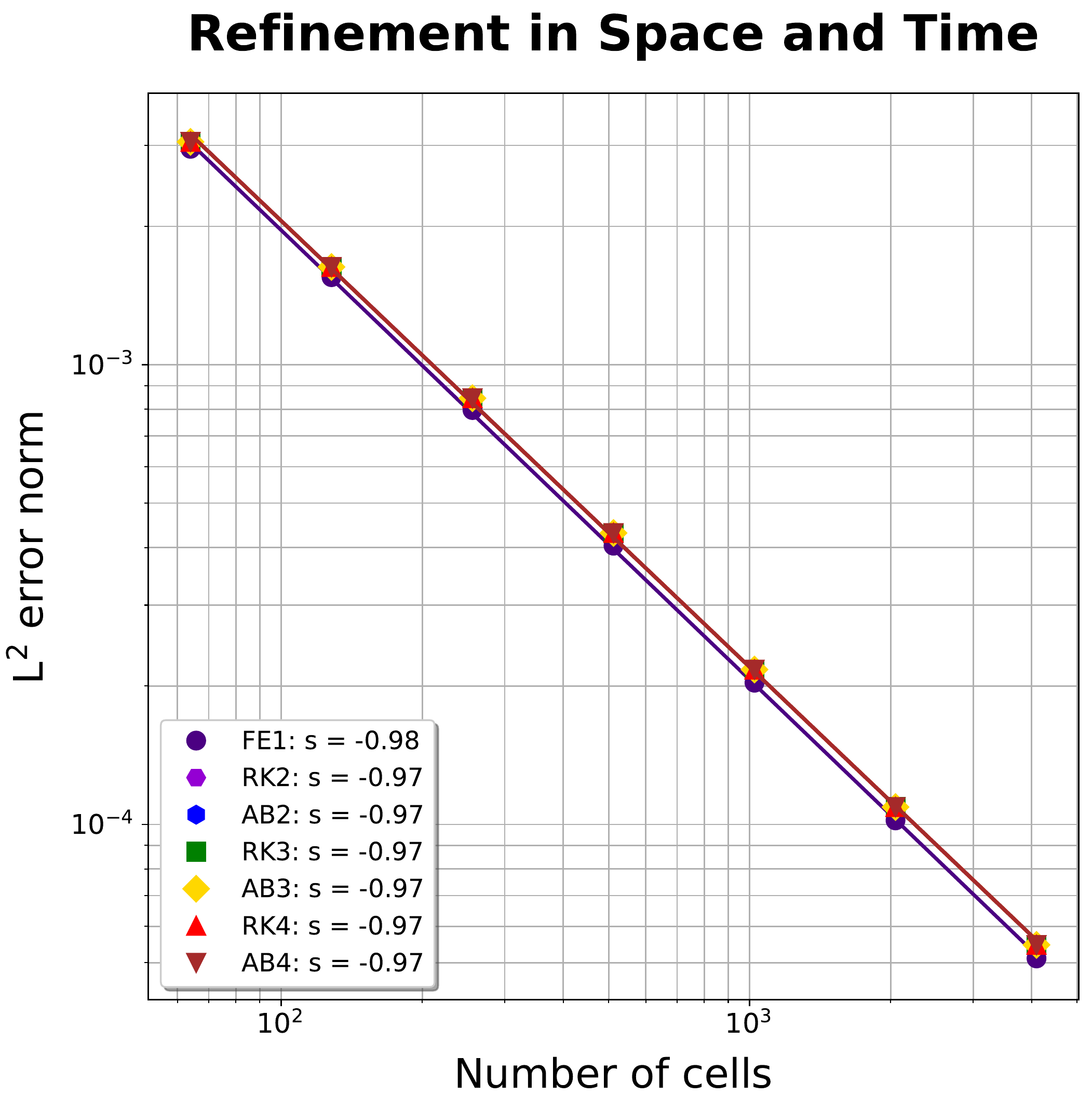}
\includegraphics[scale=.3125]{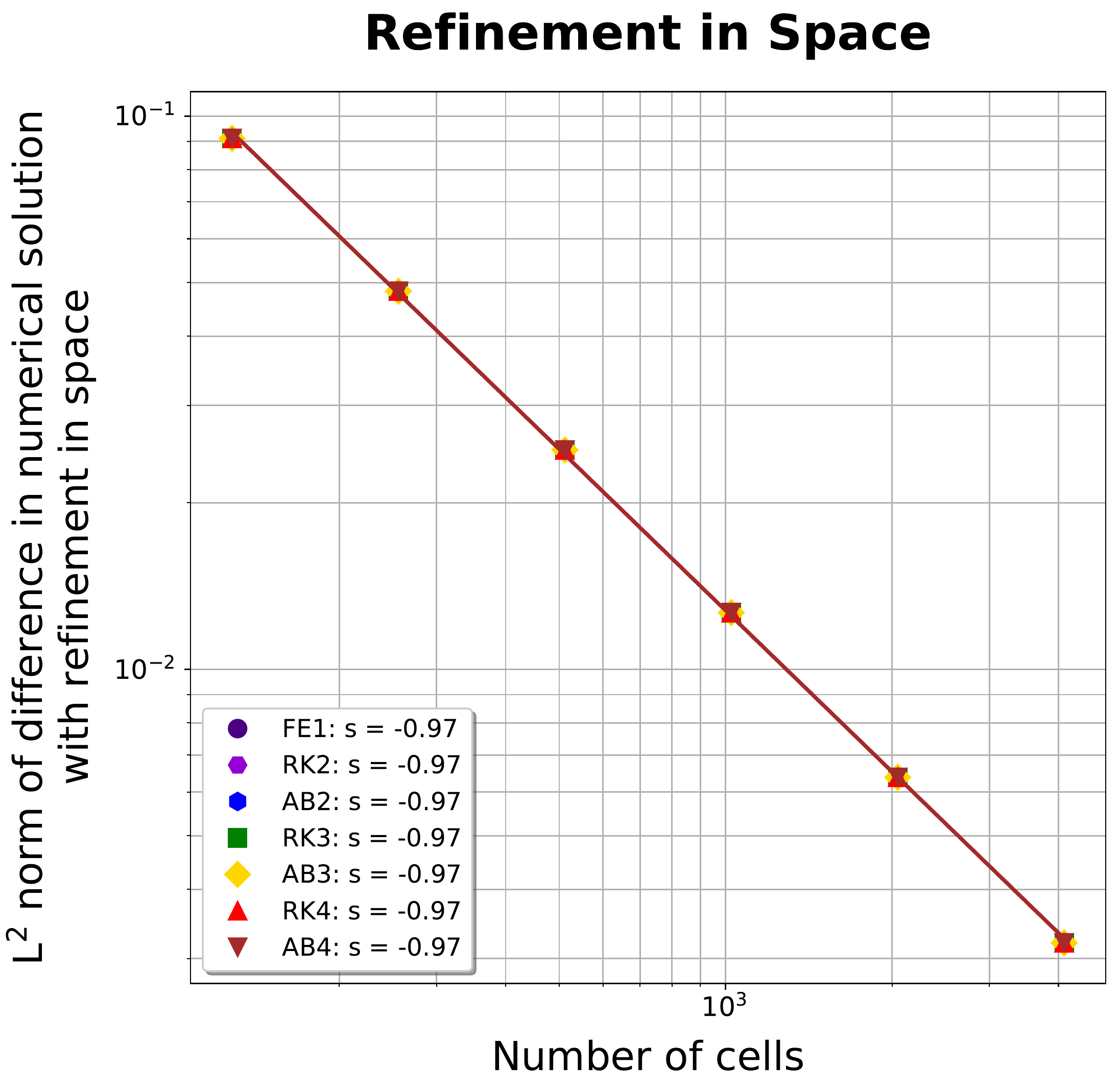} \hspace{0.5cm}
\includegraphics[scale=.3125]{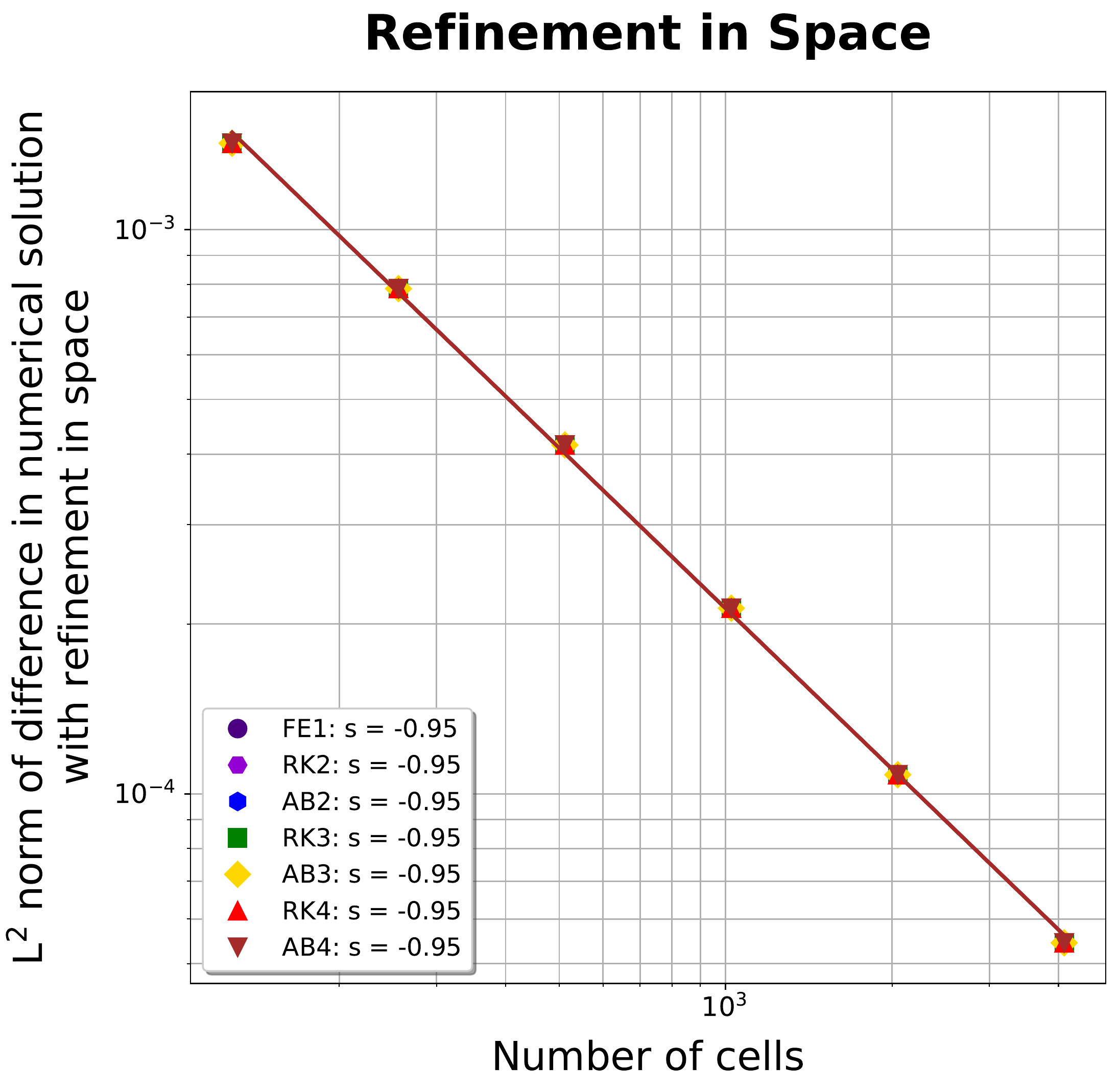}
\includegraphics[scale=.3125]{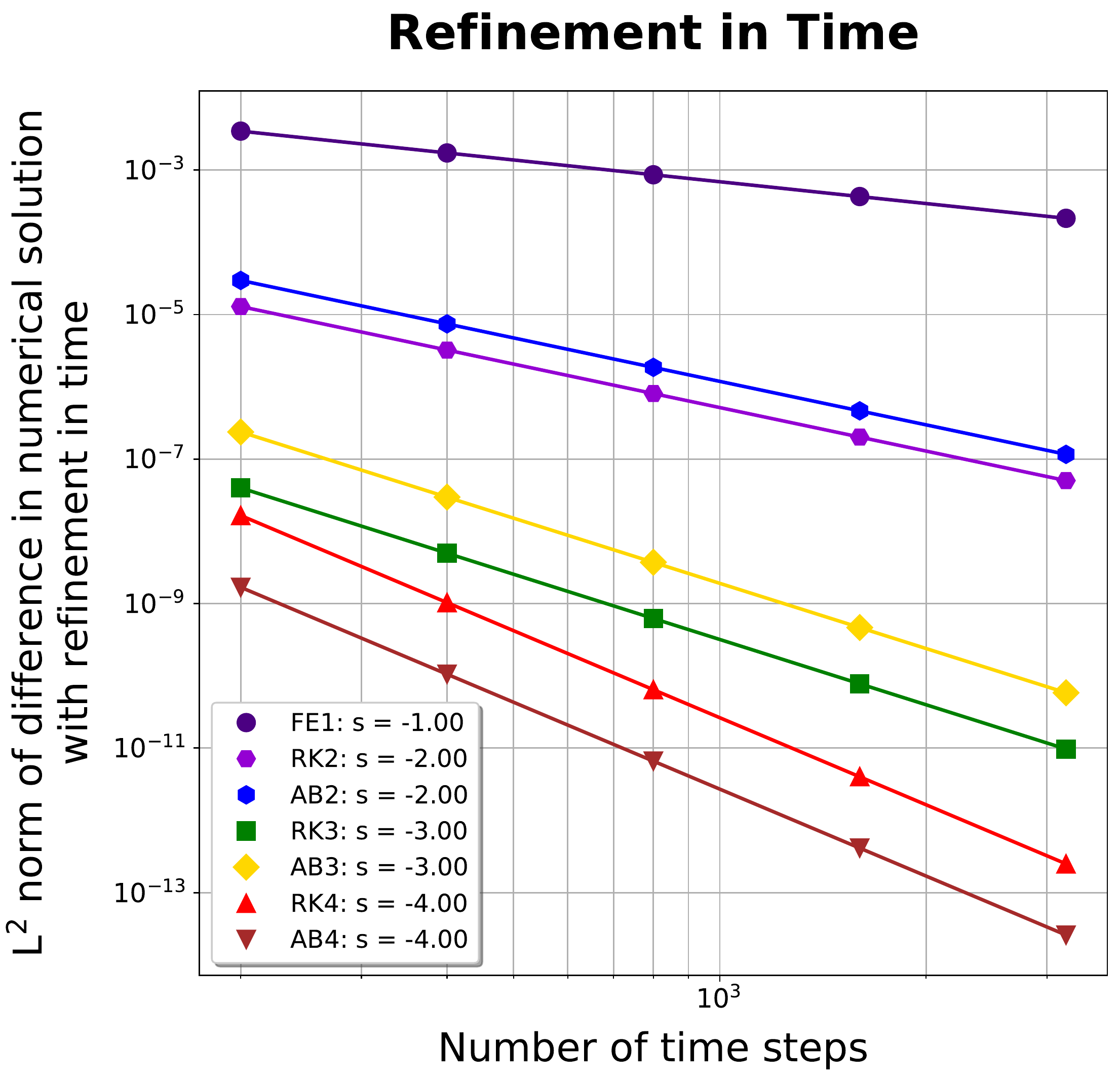} \hspace{0.5cm}
\includegraphics[scale=.3125]{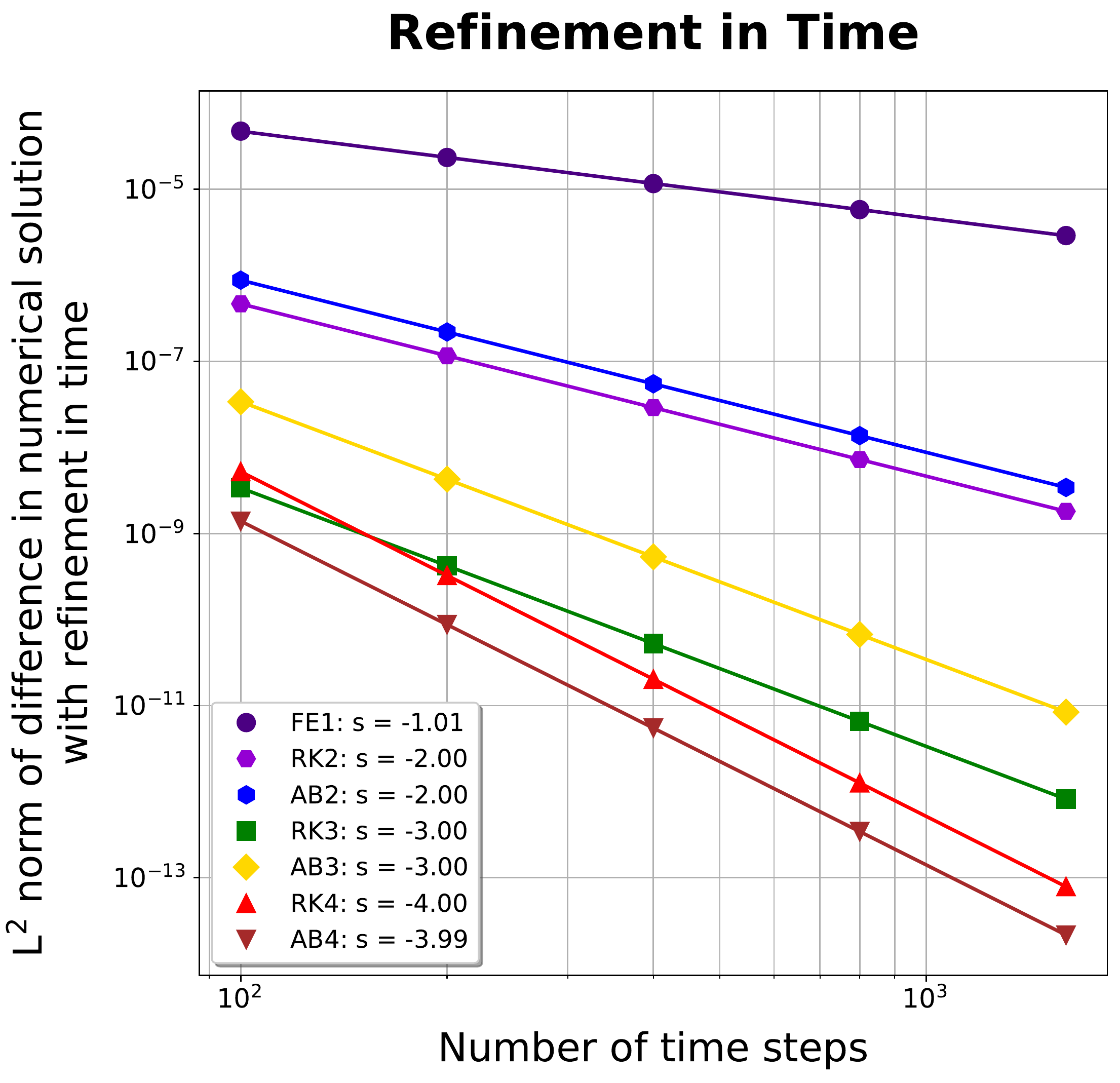}
\caption{Convergence plots of the linear variable-coefficient advection equation \eqref{LinearAdvection1D_NumericalExperiment} (left column) and the non-linear advection equation~\eqref{NonLinearAdvection1D_NumericalExperiment} (right column) using first-order upwind in space, for refinement in both space and time (first row), refinement only in space (second row), and refinement only in time (third row). Abbreviations for the seven time-stepping methods (legends) are given in List 3.} \label{Figure_Numerical_Convergence_Plots_1}
\end{figure}

\begin{figure}[!htp]
\centering
\hspace{0.165cm}
\includegraphics[scale=.3175]{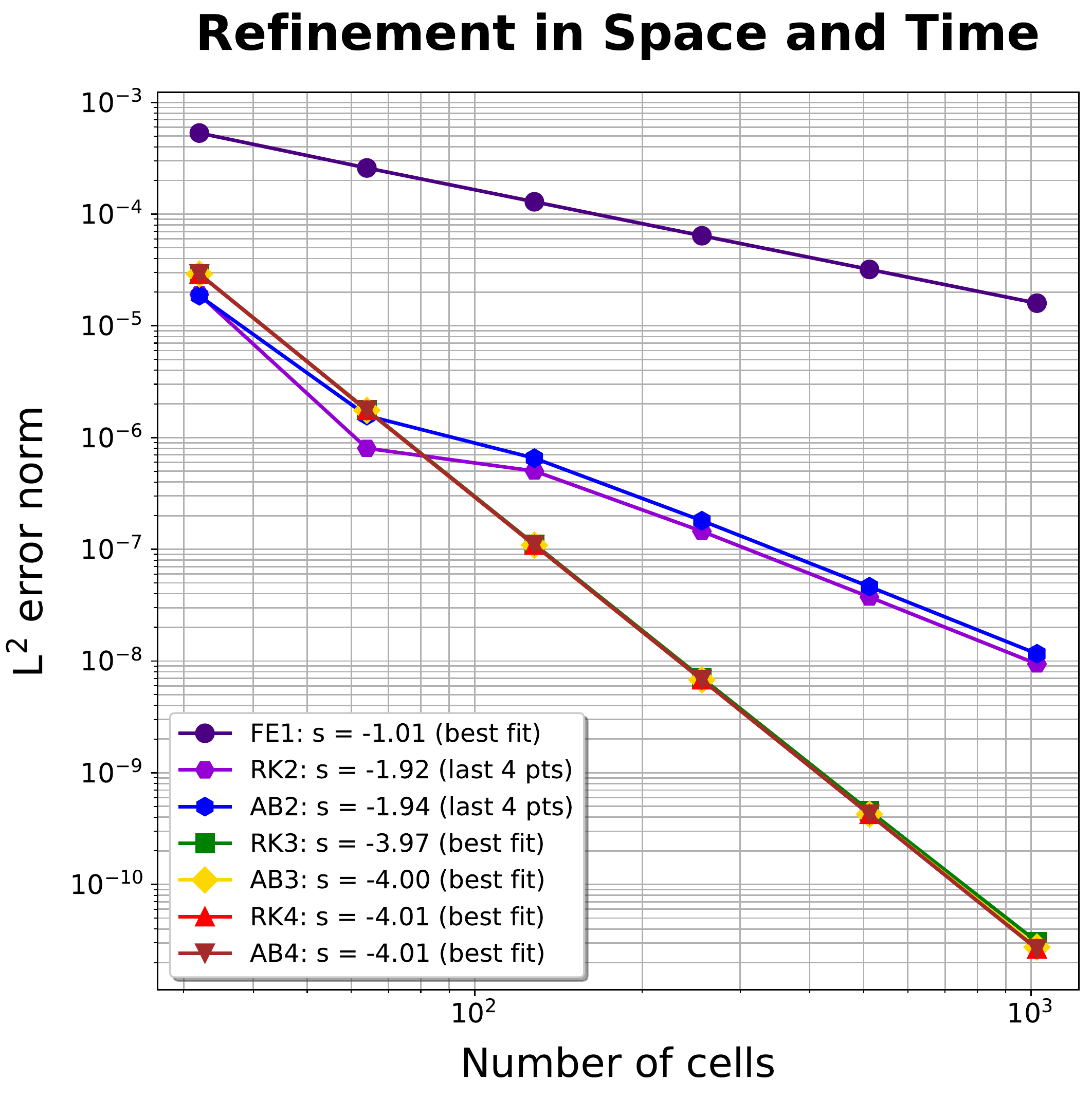} \hspace{0.785cm}
\includegraphics[scale=.3175]{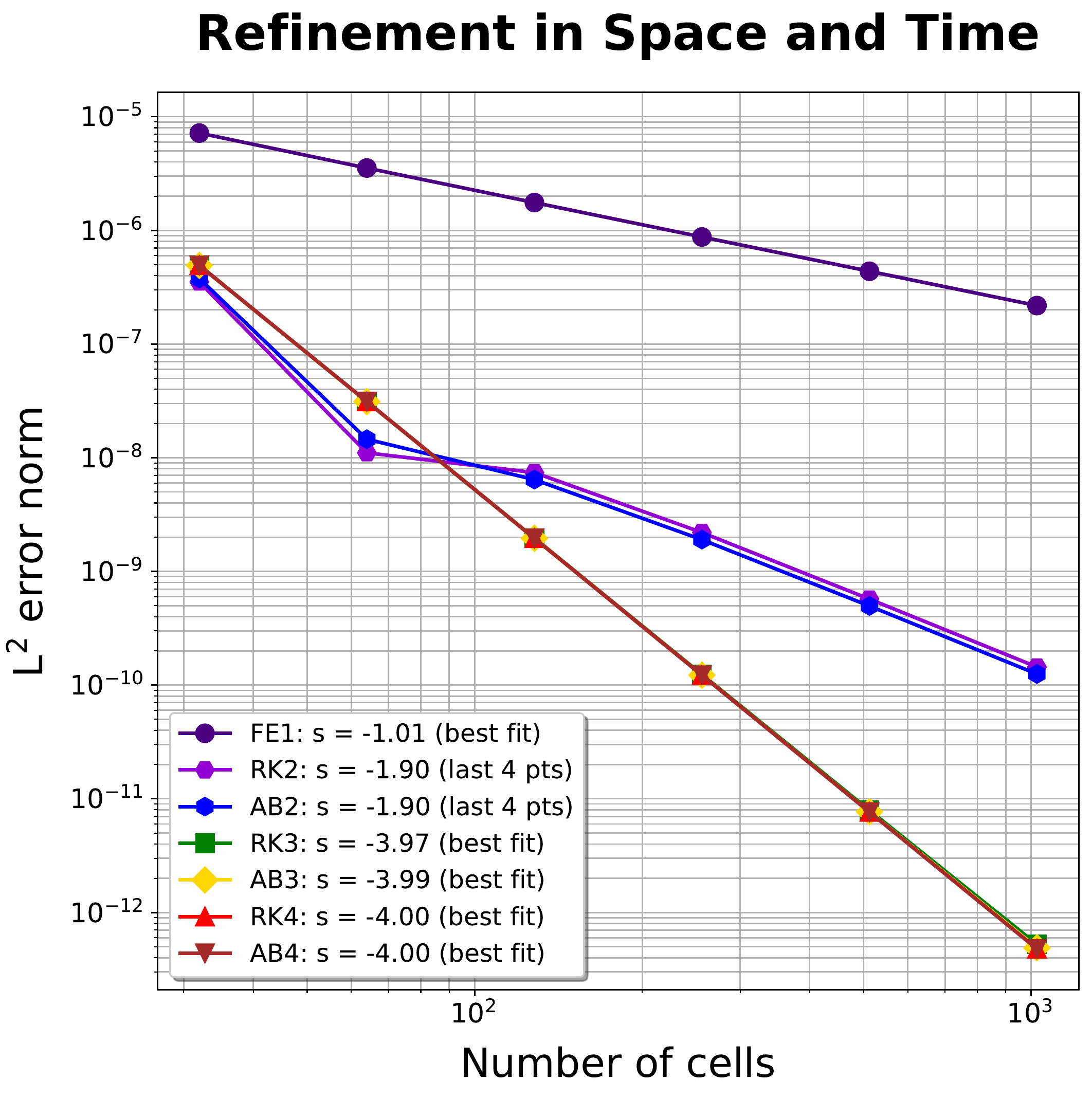}
\includegraphics[scale=.3175]{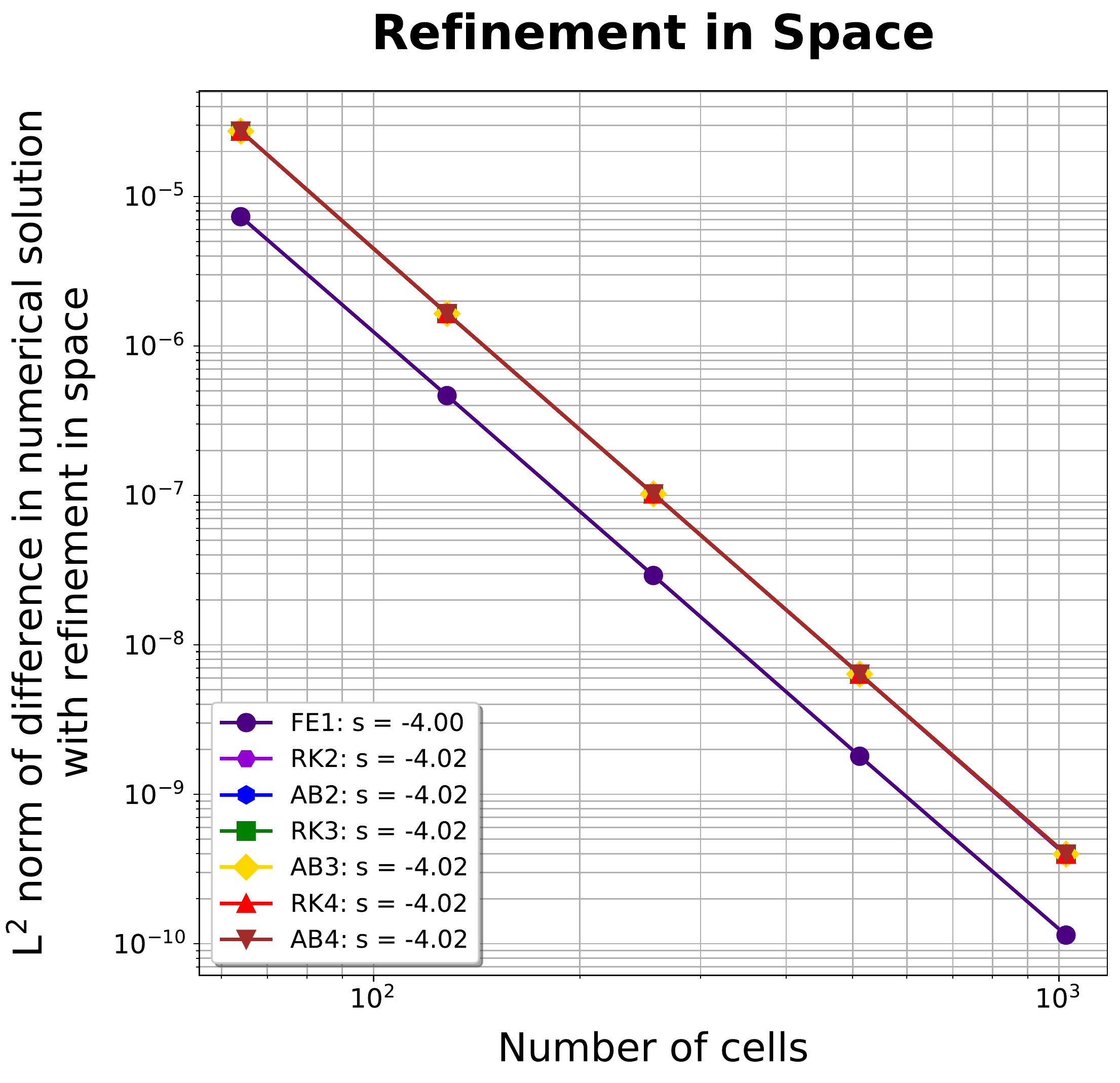} \hspace{0.5cm}
\includegraphics[scale=.3175]{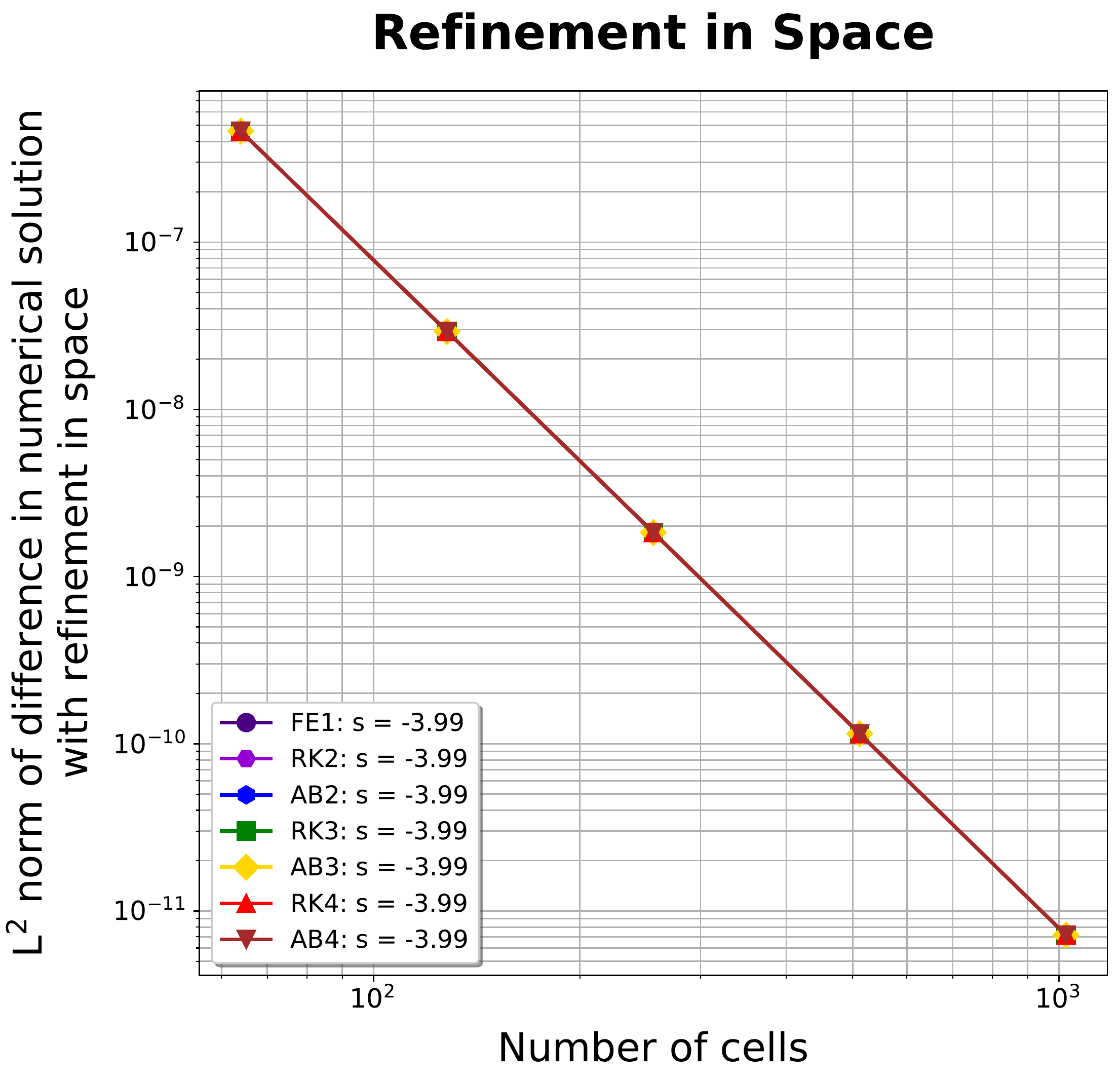}
\includegraphics[scale=.3175]{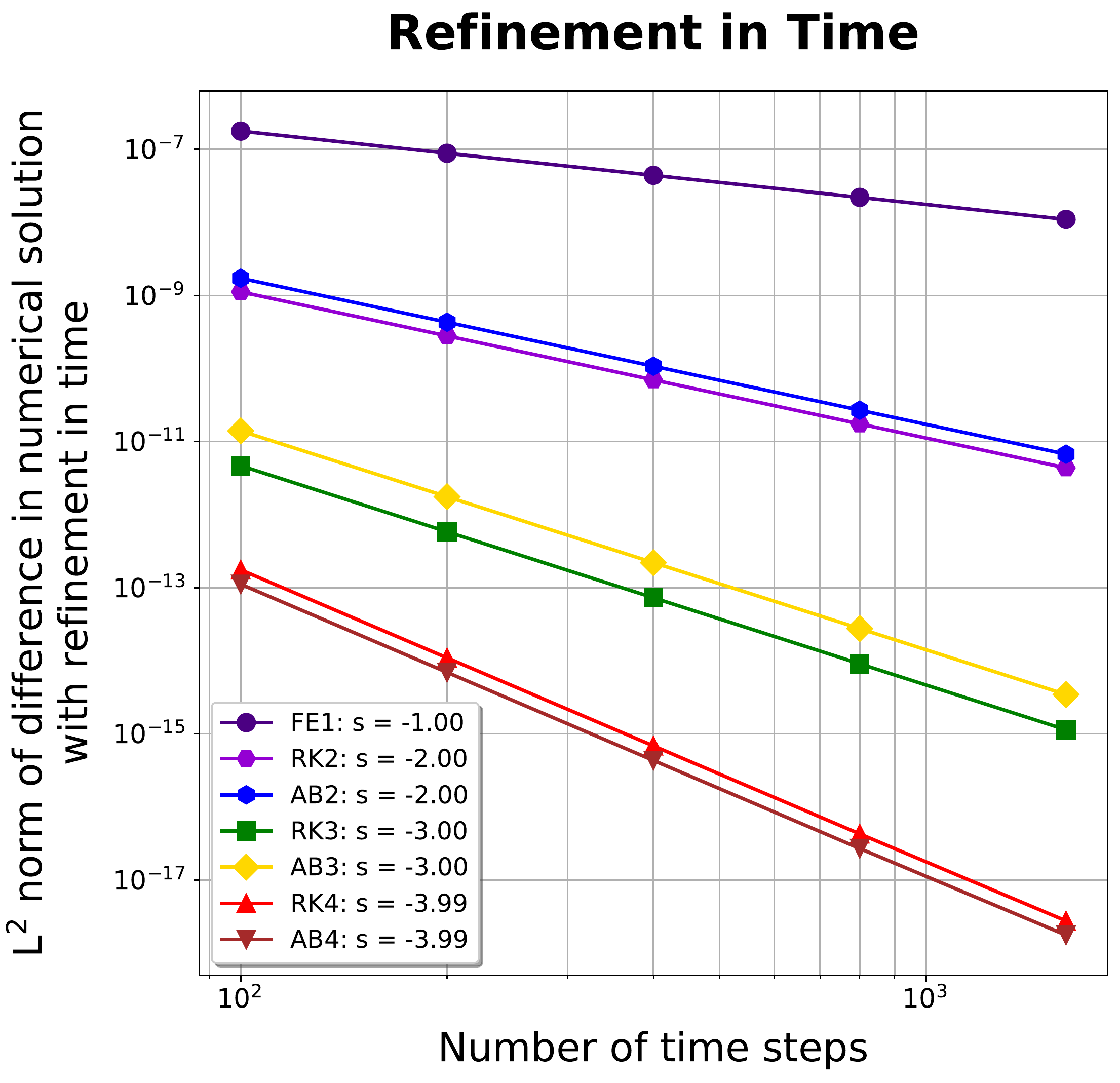} \hspace{0.5cm}
\includegraphics[scale=.3175]{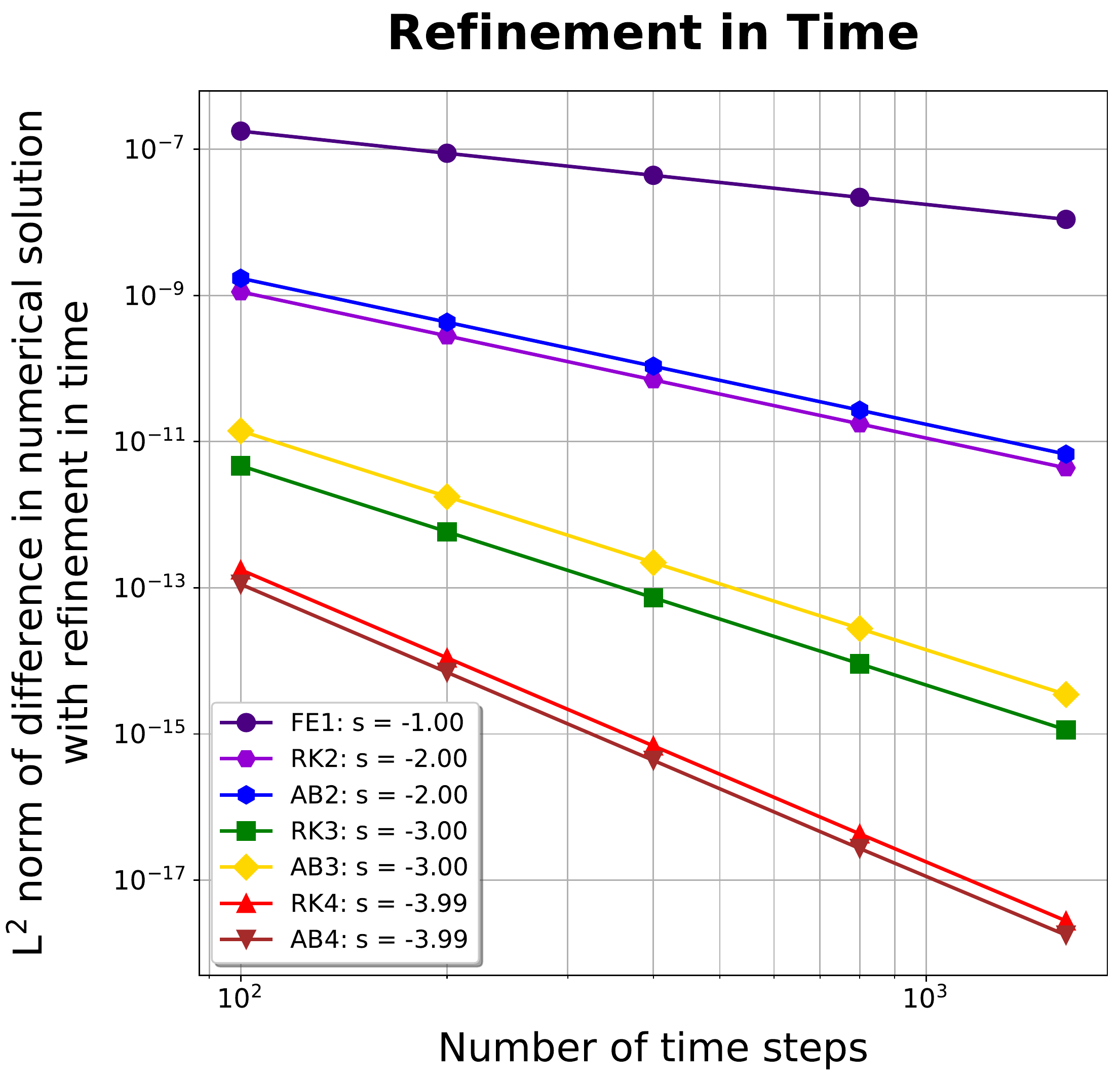}
\caption{Similar to Figure \ref{Figure_Numerical_Convergence_Plots_2} but using fourth-order accurate PPR in space, without the application of the slope limiter and monotonicity-preserving strategies.} \label{Figure_Numerical_Convergence_Plots_2}
\end{figure}

\begin{figure}[!htp]
\centering
\hspace{0.165cm}
\includegraphics[scale=.3175]{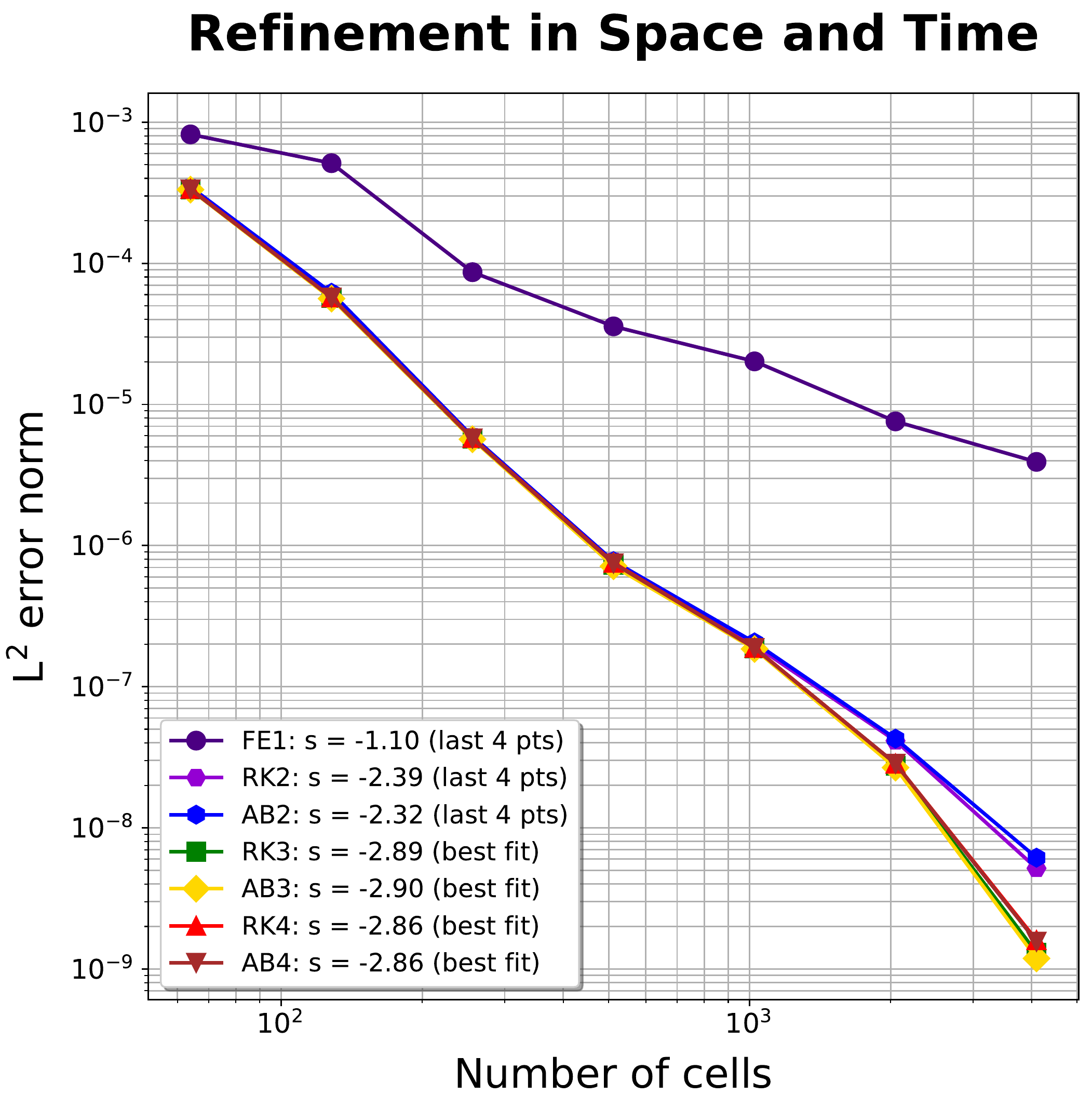} \hspace{0.785cm}
\includegraphics[scale=.3175]{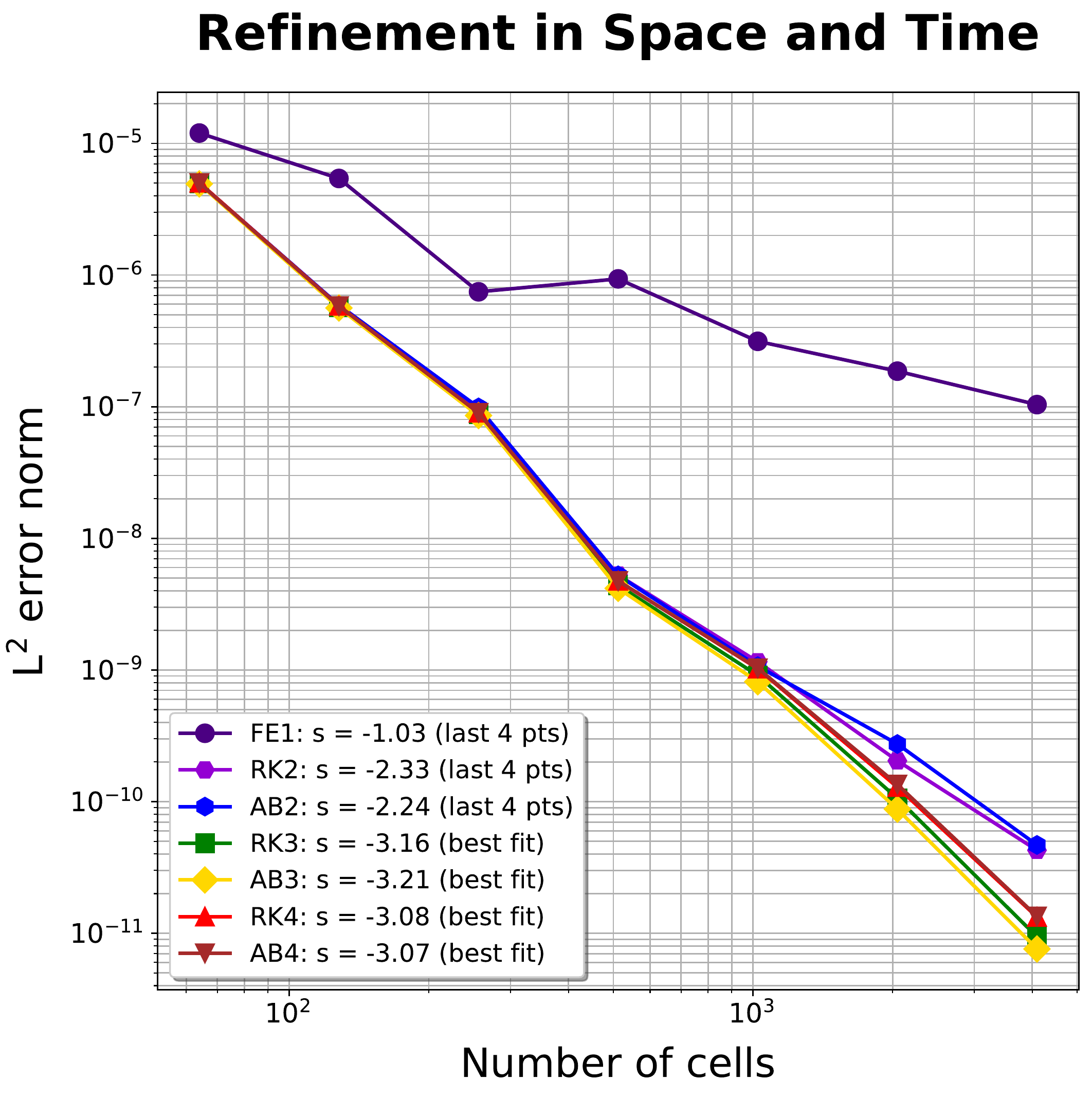}
\includegraphics[scale=.3175]{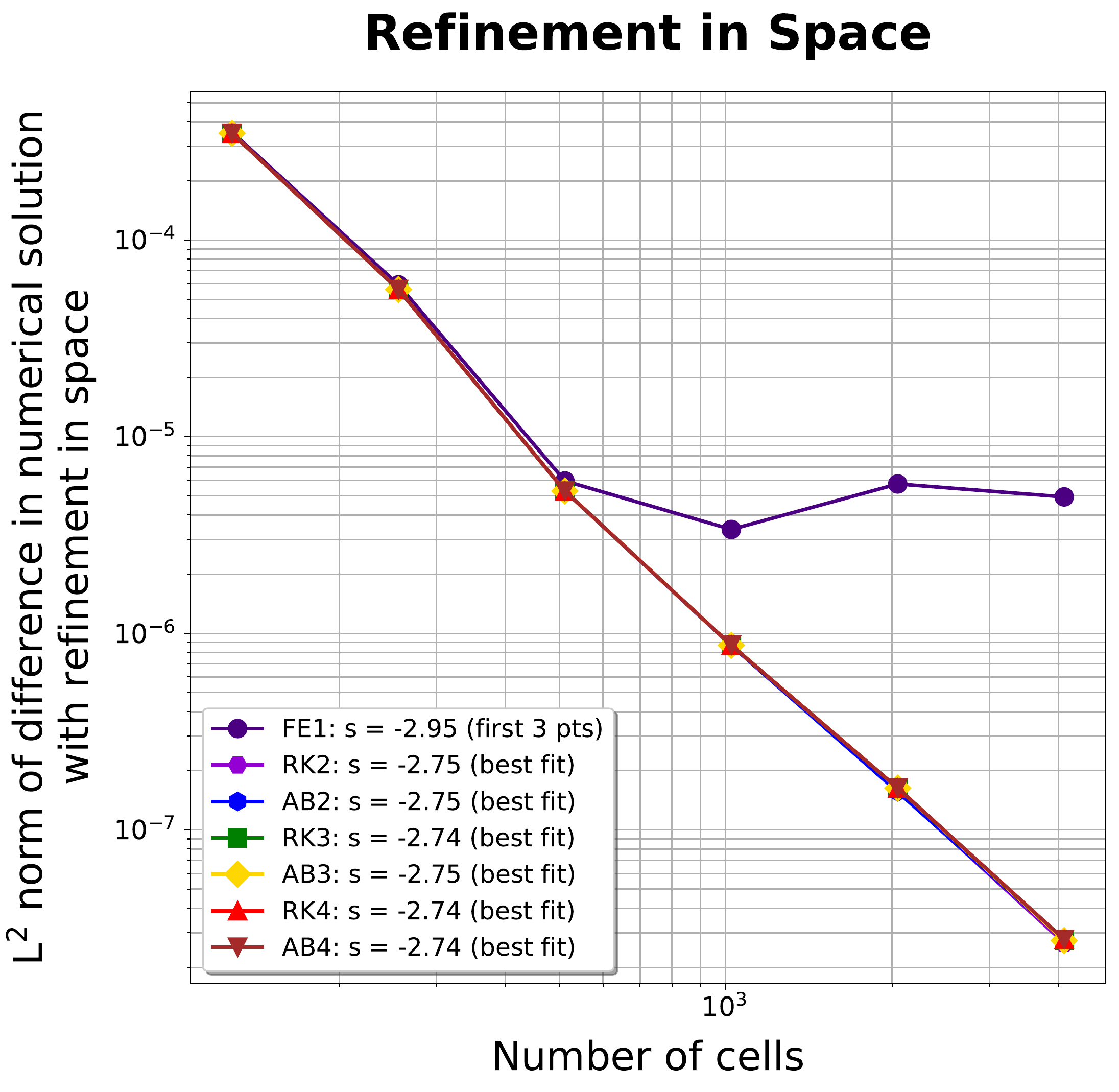} \hspace{0.5cm}
\includegraphics[scale=.3175]{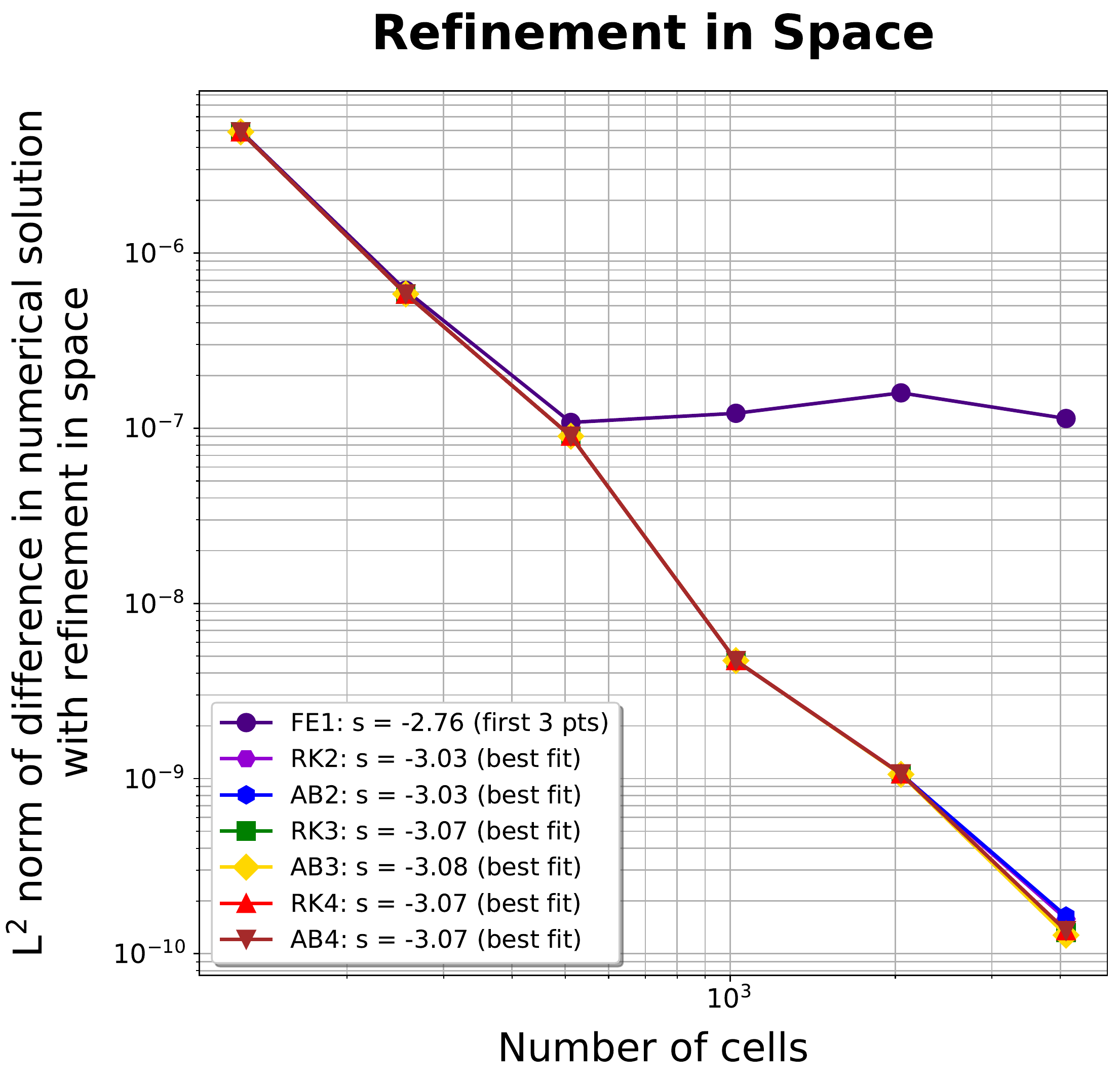}
\includegraphics[scale=.3175]{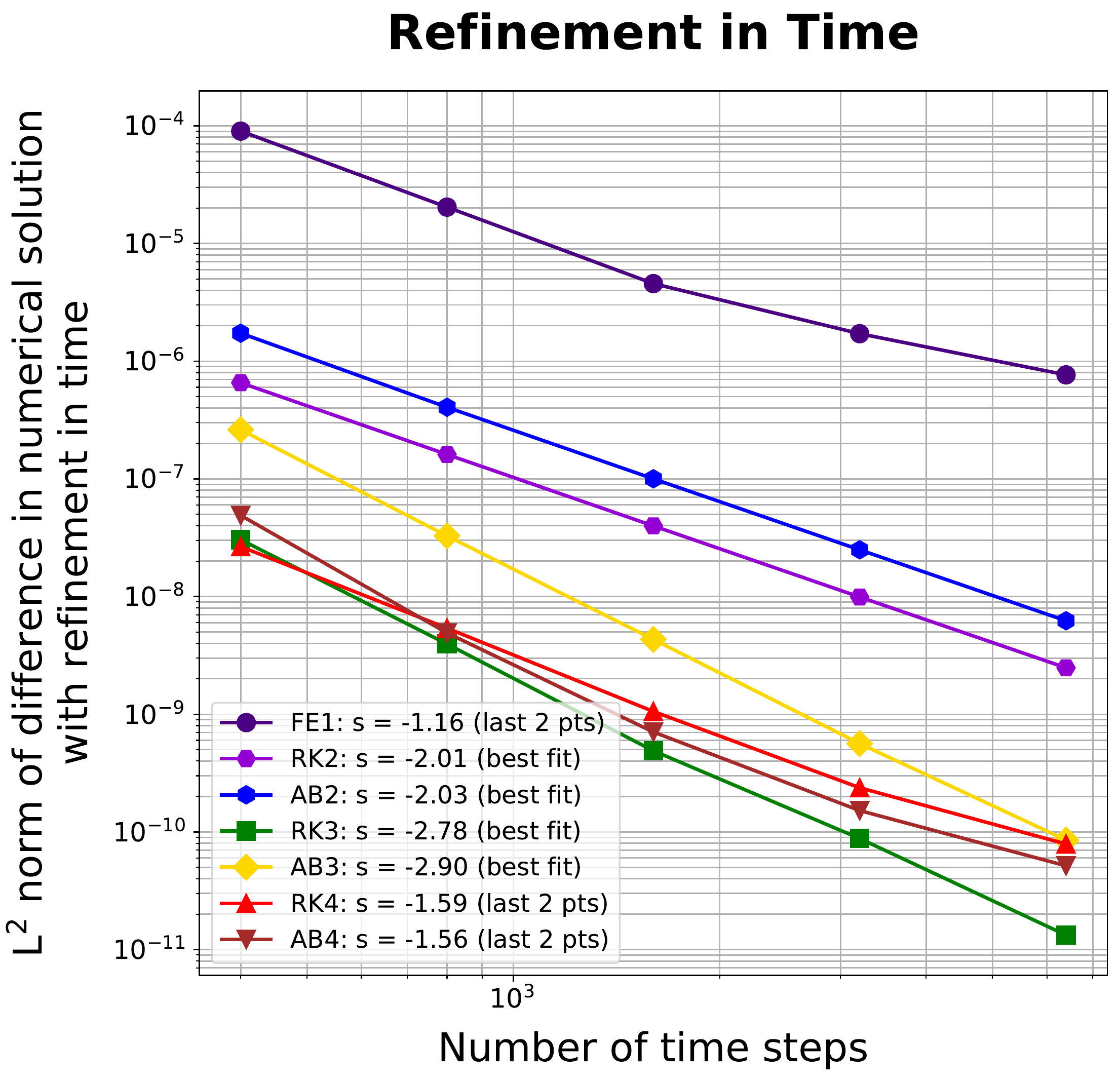} \hspace{0.5cm}
\includegraphics[scale=.3175]{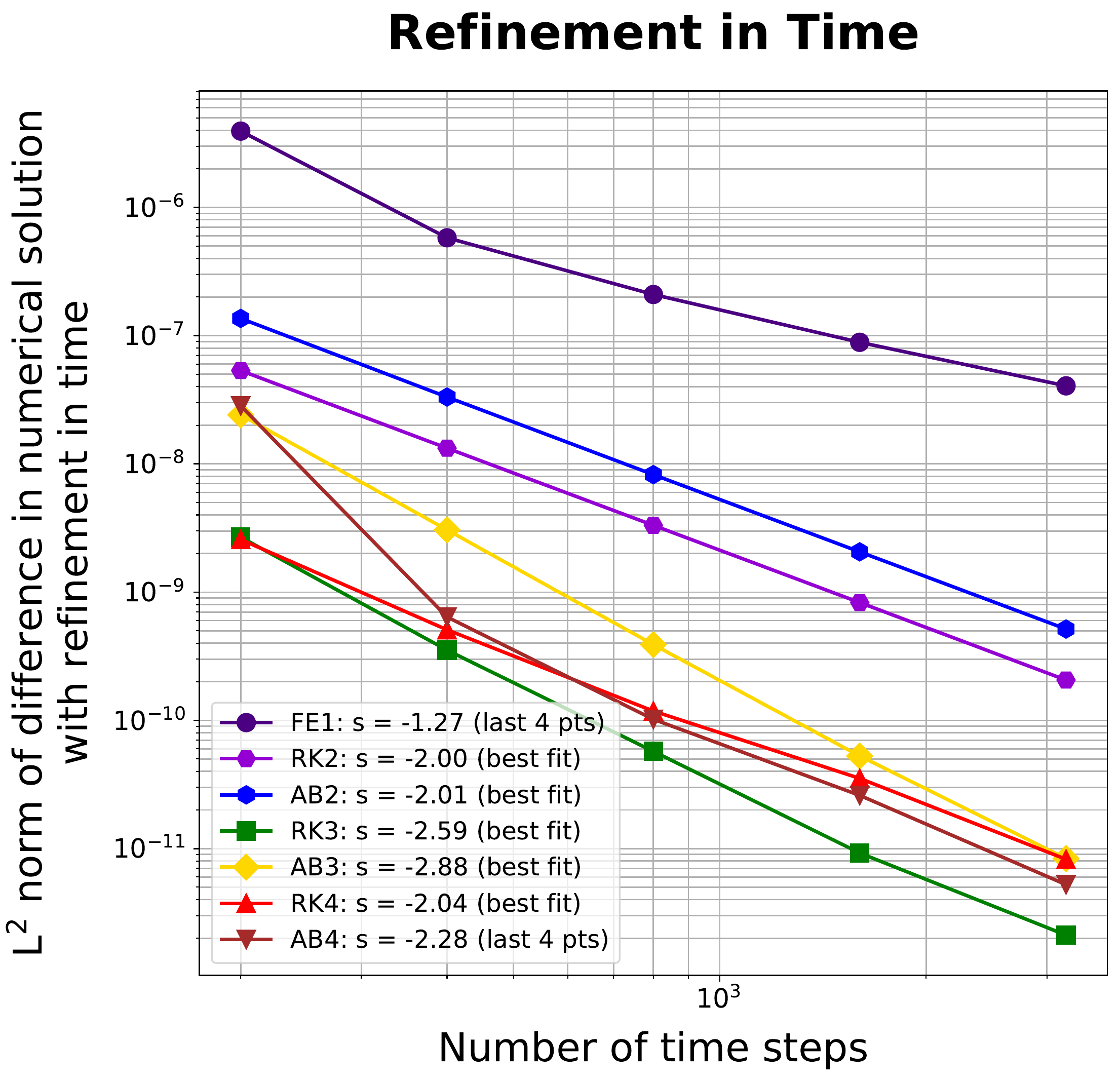}
\caption{Similar to Figure \ref{Figure_Numerical_Convergence_Plots_2} but using fourth-order accurate PPR in space, with the application of the slope limiter and monotonicity-preserving strategies.} \label{Figure_Numerical_Convergence_Plots_3}
\end{figure}

Figures~\ref{Figure_Numerical_Convergence_Plots_1}--\ref{Figure_Numerical_Convergence_Plots_3} show the convergence plots of the linear variable-coefficient advection equation~\eqref{LinearAdvection1D_NumericalExperiment} and the non-linear advection equation~\eqref{NonLinearAdvection1D_NumericalExperiment} using first-order upwind spatial discretization (Figure~\ref{Figure_Numerical_Convergence_Plots_1}), and fourth-order accurate PPR in space, without (Figure~\ref{Figure_Numerical_Convergence_Plots_2}) and with (Figure~\ref{Figure_Numerical_Convergence_Plots_3}) the application of the slope limiter and monotonicity-preserving strategies. From here onward, we refer to these two finite volume methods as non-monotone and monotone, respectively. The seven time-stepping methods of List \ref{myListOfTimeSteppingMethodsForNumericalExperiments} are applied. Refinement is performed in both space and time (first row), only in space (second row), and only in time (third row). Since our domain size is one, the cell width is $\Delta x = 1/N_{\text{cells}}$, where $N_{\text{cells}}$ is the number of cells. Table~\ref{Table_Numerical_Experiments_Parameters} lists the values of $N_{\text{cells}}$ for the various combinations of spatial and temporal discretizations, and refinement types. Letting $\eta = \Delta t/\Delta x$ denote the ratio of the time step to the cell width, we define two more parameters $\hat{\eta}_{\text{space}}$ and $\hat{\eta}_{\text{time}}$, and assign to them the values listed in Table~\ref{Table_Numerical_Experiments_Parameters}. For refinement in both space and time, we specify $\eta = \hat{\eta}_{\text{space}}$ and $\Delta t = \hat{\eta}_{\text{space}} \Delta x$. For refinement only in space, we specify $\eta = \hat{\eta}_{\text{space}}$ for the smallest value of $\Delta x$, say $\Delta x_{\text{smallest}}$, and use $\Delta t = \hat{\eta}_{\text{space}} \Delta x_{\text{smallest}}$ throughout the study. For refinement only in time, we specify the largest $\Delta t$, say $\Delta t_{\text{largest}}$, as $\Delta t_{\text{largest}} = \hat{\eta}_{\text{time}} \Delta x$, and vary $\Delta t$ from $\Delta t_{\text{largest}}$ to $\Delta t_{\text{largest}}/2^5$ by factors of $1/2$. We know that the width of the absolute stability region of the Runge-Kutta methods, used to advance the characteristic ODE, either remains the same or slightly increases with the order, whereas absolute stability region of the Adams-Bashforth methods decreases by almost half with an increase in order. This is what motivated us to specify smaller values of $\hat{\eta}_{\text{space}}$ and $\hat{\eta}_{\text{time}}$ for the Adams-Bashforth methods than for the Runge-Kutta methods. Table~\ref{Table_Numerical_Experiments_Parameters} also lists the values of the time horizon $T_{\text{horizon}}$ as a fractions of the time period of the first wave mode, $T_1$. Since stability is not guaranteed for the non-monotone finite volume method, we are sometimes compelled to use smaller values of $N_{\text{cells}}$, $\hat{\eta}_{\text{space}}$, $\hat{\eta}_{\text{time}}$, and $T_{\text{horizon}}$, so that the numerical solutions remain bounded. Even though it is advisable to ensure monotonicity in practice, for the purpose of numerical verification of our theory, the non-monotone method has proven to be a helpful exercise. 

\begin{table}[!htp]
\centering
\caption{Choice of parameters $N_{\text{cells}}$, $\hat{\eta}_{\text{space}}$, $\hat{\eta}_{\text{time}}$, and $T_{\text{horizon}}$ for the linear variable-coefficient advection equation \eqref{LinearAdvection1D_NumericalExperiment} and the non-linear advection equation \eqref{NonLinearAdvection1D_NumericalExperiment} using first-order upwind and fourth-order PPR in space, with and without the application of the slope limiter and monotonicity-preserving strategies, for refinement in both space and time, refinement only in space, and refinement only in time. The parameter $\hat{\eta}_{\text{space}}$ is defined for refinement in both space and time, and for refinement only in space, whereas the parameter $\hat{\eta}_{\text{time}}$ is defined for refinement only in time. $T_1$ represents the time period of the first wave mode.}
\setlength{\tabcolsep}{0.6em}
\renewcommand{\arraystretch}{1.125}
\begin{tabular}{cccccc}
\toprule
\multirow{2}{*}{Parameter} & Advection & Refinement & Spatial & Time & \multirow{2}{*}{Values} \\
& Type & Type & Discretization & Integrators & \\
\midrule
\multirow{10}{*}{$N_{\text{cells}}$} & \multirow{2}{*}{Both} & Space-Time, & First-Order Upwind, & \multirow{2}{*}{All} & \multirow{2}{*}{$2^6,\ldots,2^{12}$} \\
& & Space & Fourth-Order PPR (Monotone) & & \\
& \multirow{2}{*}{Both} & Space-Time, & \multirow{2}{*}{Fourth-Order PPR (Non-Monotone)} & \multirow{2}{*}{All} & \multirow{2}{*}{$2^5,\ldots,2^{10}$} \\
& & Space & & & \\
& \multirow{2}{*}{Linear} & \multirow{2}{*}{Time} & First-Order Upwind, & \multirow{2}{*}{All} & \multirow{2}{*}{$2^7$} \\
& & & Fourth-Order PPR (Monotone) & & \\
& \multirow{2}{*}{Non-Linear} & \multirow{2}{*}{Time} & First-Order Upwind, & \multirow{2}{*}{All} & \multirow{2}{*}{$2^6$} \\
& & & Fourth-Order PPR (Monotone) & & \\
& Both & Time & Fourth-Order PPR (Non-Monotone) & All & $2^7$ \\
\midrule
\multirow{10}{*}{$\hat{\eta}_{\text{space}}$} & \multirow{2}{*}{Both} & \multirow{2}{*}{Both} & \multirow{2}{*}{First-Order Upwind} & AB4 & 0.125 \\
& & & & Rest & 0.25 \\
& Both & Both & Fourth-Order PPR (Monotone) & FE1 & 0.2 \\
& Linear & \multirow{2}{*}{Space-Time} & \multirow{2}{*}{Fourth-Order PPR (Non-Monotone)} & \multirow{2}{*}{FE1} & 0.15 \\
& Non-Linear & & & & 0.2 \\
& Linear & \multirow{2}{*}{Space} & \multirow{2}{*}{Fourth-Order PPR (Non-Monotone)} & \multirow{2}{*}{FE1} & 0.0125 \\
& Non-Linear & & & & 0.1 \\
& \multirow{3}{*}{Both} & \multirow{3}{*}{Both} & \multirow{3}{*}{Fourth-Order PPR (Both)} & AB2, AB3, AB4 & 0.15 \\
& & & & RK2 & 0.2 \\
& & & & RK3, RK4 & 0.25 \\
\midrule
\multirow{3}{*}{$\hat{\eta}_{\text{time}}$} & \multirow{2}{*}{Both} & \multirow{2}{*}{Time} & \multirow{2}{*}{First-Order Upwind} & RK4 & 0.32 \\
& & & & Rest & 0.16 \\
& Both & Time & Fourth-Order PPR (Both) & All & 0.16 \\
\midrule
\multirow{7}{*}{$T_{\text{horizon}}$} & Both & All & First-Order Upwind & All & $0.25T_1$ \\
& Both & All & Fourth-Order PPR (Monotone) & All & $0.25T_1$ \\
& \multirow{2}{*}{Linear} & \multirow{2}{*}{Space} & \multirow{2}{*}{Fourth-Order PPR (Non-Monotone)} & FE1 & $0.03125T_1$ \\
& & & & Rest & $0.125T_1$ \\
& \multirow{2}{*}{Linear} & Space-Time, & \multirow{2}{*}{Fourth-Order PPR (Non-Monotone)} & \multirow{2}{*}{All} & \multirow{2}{*}{$0.125 T_1$} \\
& & Time & & & \\
& Non-Linear & All & Fourth-Order PPR (Non-Monotone) & All & $0.125 T_1$ \\
\bottomrule 
\end{tabular} \label{Table_Numerical_Experiments_Parameters}
\end{table}

\begin{table}[!htp]
\centering
\caption{Spatial and temporal discretizations for which the error of the linear variable-coefficient advection equation \eqref{LinearAdvection1D_NumericalExperiment} and the non-linear advection equation \eqref{NonLinearAdvection1D_NumericalExperiment} decreases with increase in $\Delta \xi$, for at least some values of $\Delta \xi$, where $\Delta \xi = \Delta x$ for refinement only in space, and $\Delta \xi = \Delta t$ for refinement only in time.}
\setlength{\tabcolsep}{0.6em}
\renewcommand{\arraystretch}{1.125}
\begin{tabular}{cccc}
\toprule
Advection & Refinement & Spatial & Time \\
Type & Type & Discretization & Integrators \\
\midrule
Linear & Time & First-Order Upwind & FE1, RK2, AB2, RK4 \\
Non-Linear & Time & First-Order Upwind & FE1, AB3, RK4, AB4 \\
Linear & Time & Fourth-Order PPR (Non-Monotone) & RK2, AB2, AB3 \\
Non-Linear & Time & Fourth-Order PPR (Non-Monotone) & RK2, AB2, AB3 \\
Linear & Space & Fourth-Order PPR (Monotone) & FE1 \\
Non-Linear & Space & Fourth-Order PPR (Monotone) & FE1 \\
Linear & Time & Fourth-Order PPR (Monotone) & RK3, AB3, AB4 \\
Non-Linear & Time & Fourth-Order PPR (Monotone) & RK3, AB3 \\
\bottomrule 
\end{tabular} \label{Table_Numerical_Experiments_Reduction_in_Error_With_Refinement}
\end{table}

Before discussing the nature of the convergence plots, we point out that an increase in the error with refinement only in space or only in time is a more common phenomenon than one might think. By studying the behavior of the actual error norm with only spatial or temporal refinement, we have noted all such occurrences in our numerical experiments and listed them in Table~\ref{Table_Numerical_Experiments_Reduction_in_Error_With_Refinement}.

\pagebreak

\begin{figure}[!htp]
\centering
\includegraphics[scale=.3175]{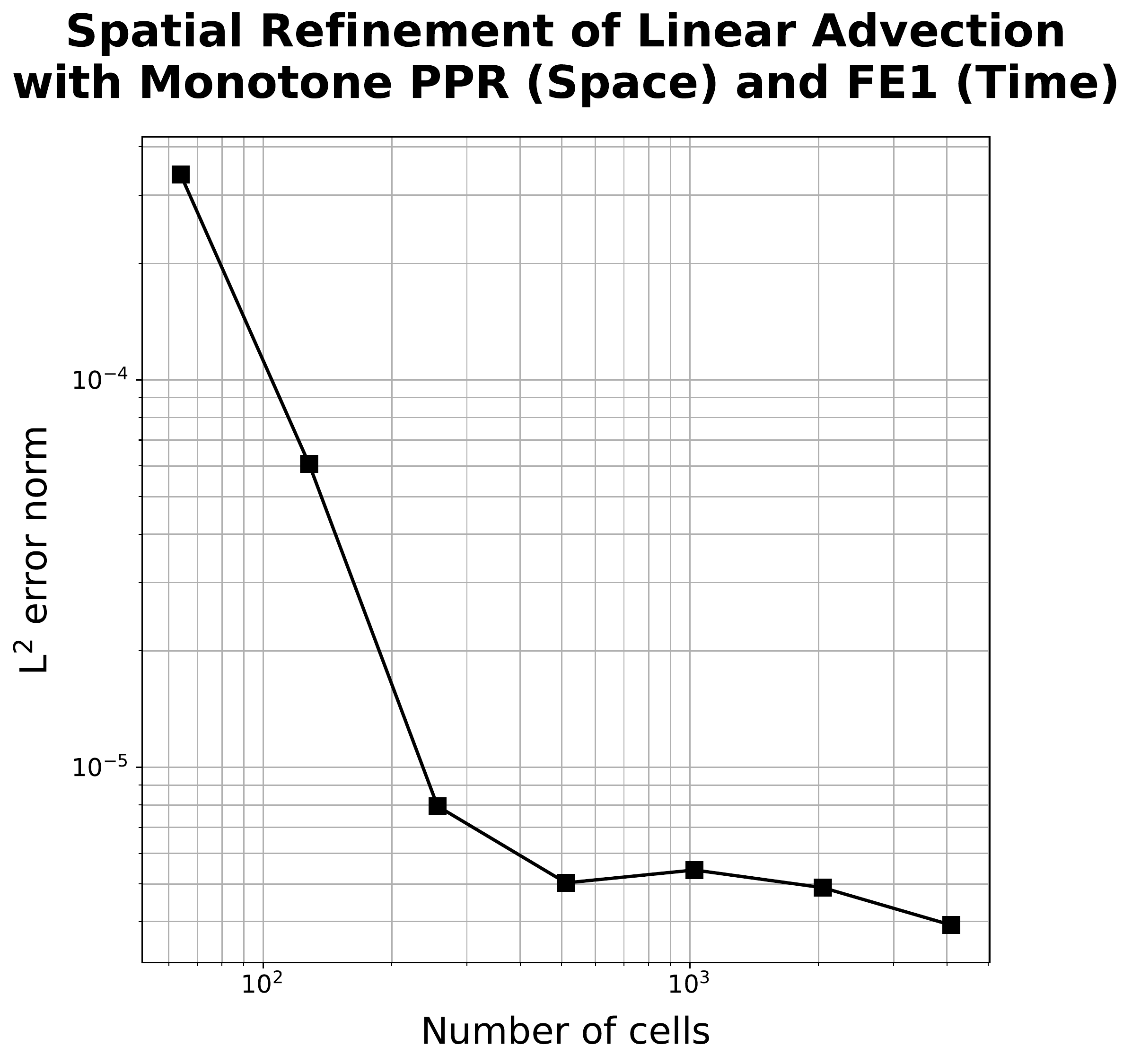} \hspace{0.15cm}
\includegraphics[scale=.3175]{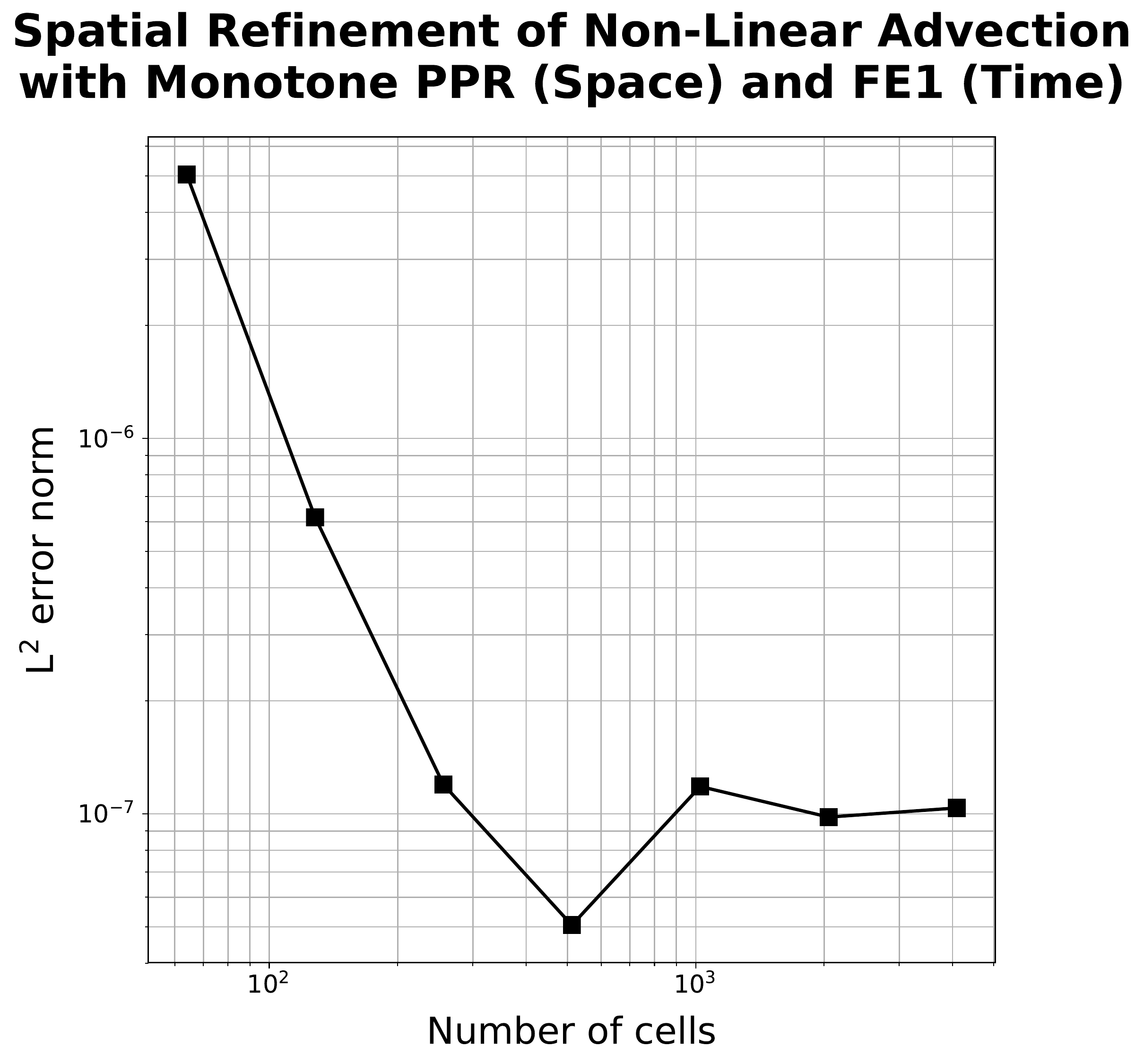}
\caption{Variation of the error of the linear variable-coefficient advection equation \eqref{LinearAdvection1D_NumericalExperiment} (left) and the non-linear advection equation \eqref{NonLinearAdvection1D_NumericalExperiment} (right) employing the monotone finite volume method and advanced with Forward Euler, with refinement only in space.} \label{NonConvergentError}
\end{figure}

The convergence curves of the linear and non-linear advection equations are similar in nature, and the following explanations are applicable to both. With the finite difference method, the spatial and temporal resolutions have reached the asymptotic regime, and the order of convergence at $\eta = \hat{\eta}_{\text{space}}$ is limited by the first-order accuracy of the spatial approximation. By plotting differences in the numerical solution or the error between successive pairs of spatial (or temporal) resolutions, we capture the true order of the spatial (or temporal) discretization. With the finite volume method, we know that the spatial approximation of the cell-integrated solution, as determined by the flux approximation of the PPR, is fourth-order accurate on a uniform mesh. With the non-monotone finite volume method, the resolutions have reached the asymptotic regime for all but the convergence study in both space and time using third-order Runge-Kutta and Adams-Bashforth methods. For these methods, the slope of the convergence curves is~4 instead of 3, which can happen for one of two reasons. The first reason is that the coefficient of $\Delta t^3$ in the global truncation error is actually zero. However, if this were the case, we would not have attained third-order convergence by plotting the difference in the norm of the numerical solution for successive temporal resolutions with refinement only in time. The second reason, and the only plausible explanation is that $\bigO\left(\Delta x^4\right) + \bigO\left(\Delta t^4\right) \gg \bigO\left(\Delta t^3\right)$ for the range of values of $\Delta x$ and $\Delta t$ used in the convergence study in both space and time. As a result, when $\Delta t \propto \Delta x$, the convergence slope is 4, and not 3. With further refinement in $\Delta x$ and $\Delta t$, we expect to reach the asymptotic regime and obtain a convergence slope of 3. However, as mentioned before, we were unable to do so and keep the solution stable, without applying the slope limiter and the monotonicity-preserving strategies. With the explicit midpoint method and the second-order Adams-Bashforth method, we observe order reduction, as we reach the asymptotic regime, and the slope of the convergence curves drop to 3 (as expected) after the first 3 points. Just like the finite difference method, plotting the differences in the numerical solution or the error between successive pairs of spatial (or temporal) resolutions reveals the true order of the spatial (or temporal) discretization. For the monotone finite volume method, the combined effect of (a) the slope limiter, (b) the monotonicity-preserving strategies, and (c) the dissipative Riemann solver modify the expression of the truncation error and reduce the spatial and temporal orders of accuracy. For spatial and temporal discretizations of order $\alpha$ and $\beta$, we expect the global truncation error to assume the form
\begin{equation}
\hat{\tau}_G = \bigO\left({\Delta x}^{\alpha}\right) + \Delta t \bigO\left({\Delta x}^{\alpha}\right) + {\Delta t}^2 \bigO\left({\Delta x}^{\alpha}\right) + \cdots + {\Delta t}^{\beta-1} \bigO\left({\Delta x}^{\alpha}\right) + \bigO\left({\Delta t}^{\beta}\right). \label{GlobalTruncationErrorNumericalSolutionFinalForm_Chapter3Summary}
\end{equation}
However, operations (a)--(c) can modify the global truncation error to
\begin{equation}
\left[\hat{\tau}_G\right]_{\text{modified}} = \bigO\left({\Delta x}^{\alpha_0}\right) + \Delta t \bigO\left({\Delta x}^{\alpha_1}\right) + {\Delta t}^2 \bigO\left({\Delta x}^{\alpha_2}\right) + \cdots + {\Delta t}^{\beta-1} \bigO\left({\Delta x}^{\alpha_{\beta-1}}\right) + \bigO\left({\Delta t}^{\beta}\right), \label{ModifiedGlobalTruncationErrorNumericalSolutionFinalForm_Chapter3Summary}
\end{equation}
where $\alpha_i \le \alpha$ for $i=0,1,\ldots,\beta-1$. This effectively reduces the spatial order of accuracy to
\begin{equation}
[\alpha]_{\text{modified}} = \min\left(\alpha_0,\alpha_1,\ldots,\alpha_{\beta-1}\right),
\end{equation}
and the temporal order of accuracy to
\begin{equation}
[\beta]_{\text{modified}} = r,
\end{equation}
where $r \in [1, \beta]$ signifies the first occurrence of $\alpha_r = 0$ and an $\bigO(1)$ term as the coefficient of $\Delta t^r$. The continued application of operations (a)--(c) can keep introducing an $\bigO(1)$ term to every coefficient of $\Delta t^i$ for $i=0,1,\ldots$, and as a result, we cannot expect to attain convergence at all in the asymptotic regime. This is indeed what we observe in the case of refinement only in space with the Forward Euler method. We will discuss the nature of this particular plot in more detail later. In many practical applications, we do not reach the asymptotic regime with the specified values of $\Delta x$ and $\Delta t$, in which case the dominant coefficients come to our rescue. We observe this phenomenon for every other convergence plot with the monotone finite volume method. With $[\alpha]_{\text{modified}} \approx 3.0$, we do not attain a convergence slope larger than 3 even with the fourth-order Runge-Kutta or Adams-Bashforth methods. We observe order reduction with the Forward Euler, the explicit midpoint and the second-order Adams-Bashforth methods, as the slope of the convergence curve reduces to almost 1 for Forward Euler and close to 2 for the other two methods after the first 4 points. For refinement only in space, we obtain a convergence slope of almost 3 for all but the Forward Euler method. For refinement only in time, we obtain expected orders of convergence for the second- and third-order methods. With Forward Euler, we observe order reduction as the convergence slope drops to 1 (as expected) after the first 3 points. However, unlike the previous experiments (employing finite difference and non-monotone PPR in space), with the fourth-order Runge-Kutta and Adams-Bashforth methods, the convergence slope drops to 2 from 4. So, the only reasonable explanation is that operations (a)--(c) reduce the value of $r$ and $[\beta]_{\text{modified}}$ to at most 2, and the coefficient of $\Delta t^2$ turns out to be the dominant one.

\begin{table}[!htp]
\centering
\caption{Coefficients $\left[\zeta_{\gamma},\zeta_{\gamma+1}\right]$ of the best fit polynomial $\zeta_{\gamma} \Delta \xi^{\gamma} + \zeta_{\gamma+1} \Delta \xi^{\gamma+1}$ to (a) the error of the linear variable-coefficient advection equation~\eqref{LinearAdvection1D_NumericalExperiment} with refinement in both space and time, and (b) the norm of the difference in the numerical solution of~\eqref{LinearAdvection1D_NumericalExperiment} for successive resolutions with refinement only in space or only in time. For refinement in both space and time and refinement only in space, $\Delta \xi = \Delta x$, and for refinement only in time, $\Delta \xi = \Delta t$. For the finite difference and the non-monotone finite volume methods, the resolutions have already reached the asymptotic regime, and $\gamma = \alpha$ for refinement only in space, $\gamma = \beta$ for refinement only in time, and $\gamma = \min(\alpha,\beta)$ for refinement in both space and time, with $\alpha$ and~$\beta$ representing the spatial and temporal orders of accuracy. For the monotone finite volume method, the same definitions hold, except $\alpha$ and $\beta$ are replaced by their reduced equivalents. Since the norm of the difference in the numerical solution with only spatial refinement is non-convergent for the monotone finite volume method advanced with Forward Euler after the first 3 points, we fit the polynomial to only these points.}
\setlength{\tabcolsep}{0.05em}
\renewcommand{\arraystretch}{1.5}
\begin{tabular}{cccccc}
\toprule
Refinement & Time & \multicolumn{3}{c}{$\left[\zeta_{\gamma},\zeta_{\gamma+1}\right]$} \\
Type & Integrator & First-Order & Fourth-Order PPR & Fourth-Order PPR \\
& & Upwind & (Non-Monotone) & (Monotone) \\
\midrule
Space-Time & FE1 & $\left[+1.18 \times 10^{+01}, -4.46 \times 10^{+01}\right]$ & $\left[+1.62 \times 10^{-02}, +2.52 \times 10^{-02}\right]$ & $\left[+1.94 \times 10^{-02}, -5.03 \times 10^{-01}\right]$ \\ 
Space-Time & RK2 & $\left[+1.31 \times 10^{+01}, -6.28 \times 10^{+01}\right]$ & $\left[+1.07 \times 10^{-02}, -3.17 \times 10^{-01}\right]$ & $\left[+2.03 \times 10^{-01}, -3.01 \times 10^{-01}\right]$ \\
Space-Time & AB2 & $\left[+1.31 \times 10^{+01}, -6.31 \times 10^{+01}\right]$ & $\left[+1.29 \times 10^{-02}, -2.77 \times 10^{-01}\right]$ & $\left[+2.20 \times 10^{-01}, -8.45 \times 10^{+00}\right]$ \\
Space-Time & RK3 & $\left[+1.31 \times 10^{+01}, -6.29 \times 10^{+01}\right]$ & $\left[-2.72 \times 10^{-02}, +3.14 \times 10^{+01}\right]$ & $\left[+1.49 \times 10^{+02}, -3.91 \times 10^{+03}\right]$ \\
Space-Time & AB3 & $\left[+1.31 \times 10^{+01}, -6.28 \times 10^{+01}\right]$ & $\left[-3.31 \times 10^{-02}, +3.15 \times 10^{+01}\right]$ & $\left[+1.47 \times 10^{+02}, -3.83 \times 10^{+03}\right]$ \\
Space-Time & RK4 & $\left[+1.31 \times 10^{+01}, -6.29 \times 10^{+01}\right]$ & $\left[+2.84 \times 10^{+01}, +6.79 \times 10^{+01}\right]$ & $\left[+1.51 \times 10^{+02}, -3.99 \times 10^{+03}\right]$ \\
Space-Time & AB4 & $\left[+1.31 \times 10^{+01}, -6.29 \times 10^{+01}\right]$ & $\left[+2.84 \times 10^{+01}, +6.79 \times 10^{+01}\right]$ & $\left[+1.51 \times 10^{+02}, -4.03 \times 10^{+03}\right]$ \\
Space & FE1 & $\left[+1.31 \times 10^{+01}, -1.82 \times 10^{+02}\right]$ & $\left[+1.27 \times 10^{+02}, -2.41 \times 10^{+02}\right]$ & $\left[+1.23 \times 10^{+03}, -6.25 \times 10^{+04}\right]$ \\
Space & Rest & $\left[+1.31 \times 10^{+01}, -1.82 \times 10^{+02}\right]$ & $\left[+1.27 \times 10^{+02}, -2.43 \times 10^{+02}\right]$ & $\left[+1.13 \times 10^{+03}, -5.13 \times 10^{+04}\right]$ \\
Time & FE1 & $\left[+2.75 \times 10^{+00}, +1.36 \times 10^{+01}\right]$ & $\left[+1.37 \times 10^{-02}, +3.76 \times 10^{-02}\right]$ & $\left[+1.53 \times 10^{-02}, +8.88 \times 10^{+01}\right]$ \\
Time & RK2 & $\left[+8.25 \times 10^{+00}, +1.12 \times 10^{+01}\right]$ & $\left[+4.89 \times 10^{-02}, +7.43 \times 10^{-02}\right]$ & $\left[+1.62 \times 10^{+00}, +8.89 \times 10^{+01}\right]$ \\
Time & AB2 & $\left[+1.90 \times 10^{+01}, -4.48 \times 10^{+01}\right]$ & $\left[+1.03 \times 10^{-01}, -4.42 \times 10^{-01}\right]$ & $\left[+3.90 \times 10^{+00}, +8.41 \times 10^{+02}\right]$ \\
Time & RK3 & $\left[+2.04 \times 10^{+01}, +5.09 \times 10^{+01}\right]$ & $\left[+1.27 \times 10^{-01}, +8.75 \times 10^{-02}\right]$ & $\left[+1.34 \times 10^{+02}, -1.56 \times 10^{+04}\right]$ \\
Time & AB3 & $\left[+1.22 \times 10^{+02}, -6.83 \times 10^{+02}\right]$ & $\left[+7.11 \times 10^{-01}, -1.00 \times 10^{+01}\right]$ & $\left[+1.08 \times 10^{+03}, -1.56 \times 10^{+04}\right]$ \\
Time & RK4 & $\left[+4.21 \times 10^{+02}, +7.98 \times 10^{+02}\right]$ & $\left[+5.24 \times 10^{+00}, +1.45 \times 10^{+01}\right]$ & $\left[+4.13 \times 10^{-02}, +4.26 \times 10^{+01}\right]$ \\
Time & AB4 & $\left[+6.96 \times 10^{+02}, -7.58 \times 10^{+03}\right]$ & $\left[+4.24 \times 10^{+00}, -6.05 \times 10^{+01}\right]$ & $\left[-2.17 \times 10^{-02}, +2.35 \times 10^{+02}\right]$ \\
\bottomrule 
\end{tabular} \label{Table_Linear_Numerical_Advection_Error_Coefficients}
\end{table}

\begin{table}[!htp]
\centering
\caption{Similar to Table \ref{Table_Linear_Numerical_Advection_Error_Coefficients} but for the non-linear advection equation~\eqref{NonLinearAdvection1D_NumericalExperiment}.}
\setlength{\tabcolsep}{0.05em}
\renewcommand{\arraystretch}{1.5}
\begin{tabular}{cccccc}
\toprule
Refinement & Time & \multicolumn{3}{c}{$\left[\zeta_{\gamma},\zeta_{\gamma+1}\right]$} \\
Type & Integrator & First-Order & Fourth-Order PPR & Fourth-Order PPR \\
& & Upwind & (Non-Monotone) & (Monotone) \\
\midrule
Space-Time & FE1 & $\left[+2.09 \times 10^{-01}, -1.28 \times 10^{+00}\right]$ & $\left[+2.23 \times 10^{-04}, +2.17 \times 10^{-04}\right]$ & $\left[+2.53 \times 10^{-04}, +1.13 \times 10^{-01}\right]$ \\
Space-Time & RK2 & $\left[+2.23 \times 10^{-01}, -1.75 \times 10^{+00}\right]$ & $\left[+1.66 \times 10^{-04}, -5.69 \times 10^{-03}\right]$ & $\left[+9.97 \times 10^{-04}, +1.95 \times 10^{-01}\right]$ \\
Space-Time & AB2 & $\left[+2.23 \times 10^{-01}, -1.75 \times 10^{+00}\right]$ & $\left[+1.43 \times 10^{-04}, -4.90 \times 10^{-03}\right]$ & $\left[+8.94 \times 10^{-04}, +2.44 \times 10^{-01}\right]$ \\
Space-Time & RK3 & $\left[+2.23 \times 10^{-01}, -1.76 \times 10^{+00}\right]$ & $\left[+2.45 \times 10^{-04}, +5.08 \times 10^{-01}\right]$ & $\left[+1.11 \times 10^{+00}, +1.23 \times 10^{+01}\right]$ \\
Space-Time & AB3 & $\left[+2.23 \times 10^{-01}, -1.76 \times 10^{+00}\right]$ & $\left[+2.20 \times 10^{-04}, +5.09 \times 10^{-01}\right]$ & $\left[+1.08 \times 10^{+00}, +1.33 \times 10^{+01}\right]$ \\
Space-Time & RK4 & $\left[+2.23 \times 10^{-01}, -1.76 \times 10^{+00}\right]$ & $\left[+5.31 \times 10^{-01}, -4.54 \times 10^{-01}\right]$ & $\left[+1.14 \times 10^{+00}, +1.12 \times 10^{+01}\right]$ \\
Space-Time & AB4 & $\left[+2.23 \times 10^{-01}, -1.76 \times 10^{+00}\right]$ & $\left[+5.30 \times 10^{-01}, -4.54 \times 10^{-01}\right]$ & $\left[+1.14 \times 10^{+00}, +1.06 \times 10^{+01}\right]$ \\
Space & FE1 & $\left[+2.22 \times 10^{-01}, -5.07 \times 10^{+00}\right]$ & $\left[+7.96 \times 10^{+00}, -1.43 \times 10^{+01}\right]$ & $\left[+1.03 \times 10^{+01}, +1.23 \times 10^{+01}\right]$ \\
Space & Rest & $\left[+2.22 \times 10^{-01}, -5.06 \times 10^{+00}\right]$ & $\left[+7.95 \times 10^{+00}, -1.43 \times 10^{+01}\right]$ & $\left[+9.39 \times 10^{+00}, +1.18 \times 10^{+02}\right]$ \\
Time & FE1 & $\left[+1.85 \times 10^{-02}, +1.80 \times 10^{-01}\right]$ & $\left[+1.40 \times 10^{-04}, +9.36 \times 10^{-04}\right]$ & $\left[+4.67 \times 10^{-04}, +6.50 \times 10^{-01}\right]$ \\ 
Time & RK2 & $\left[+7.42 \times 10^{-02}, +1.74 \times 10^{-01}\right]$ & $\left[+7.15 \times 10^{-04}, +4.26 \times 10^{-04}\right]$ & $\left[+3.37 \times 10^{-02}, +3.75 \times 10^{-01}\right]$ \\
Time & AB2 & $\left[+1.40 \times 10^{-01}, +8.77 \times 10^{-02}\right]$ & $\left[+1.10 \times 10^{-03}, -3.35 \times 10^{-03}\right]$ & $\left[+8.30 \times 10^{-02}, +3.37 \times 10^{+00}\right]$ \\ 
Time & RK3 & $\left[+2.15 \times 10^{-01}, +7.57 \times 10^{-01}\right]$ & $\left[+2.39 \times 10^{-03}, +1.35 \times 10^{-03}\right]$ & $\left[+1.54 \times 10^{+00}, -1.32 \times 10^{+02}\right]$ \\ 
Time & AB3 & $\left[+2.21 \times 10^{+00}, -1.21 \times 10^{+01}\right]$ & $\left[+7.27 \times 10^{-03}, -7.97 \times 10^{-02}\right]$ & $\left[+1.27 \times 10^{+01}, -3.04 \times 10^{+02}\right]$ \\ 
Time & RK4 & $\left[+8.26 \times 10^{+00}, +2.32 \times 10^{+01}\right]$ & $\left[+7.20 \times 10^{-02}, +6.95 \times 10^{-02}\right]$ & $\left[+9.75 \times 10^{-04}, +5.37 \times 10^{-01}\right]$ \\ 
Time & AB4 & $\left[+3.62 \times 10^{+01}, -2.62 \times 10^{+02}\right]$ & $\left[+4.55 \times 10^{-02}, +2.45 \times 10^{-01}\right]$ & $\left[+5.04 \times 10^{-04}, +1.81 \times 10^{+00}\right]$ \\ 
\bottomrule 
\end{tabular} \label{Table_NonLinear_Numerical_Advection_Error_Coefficients}
\end{table}

Now, the nature of the convergence curve of the monotone finite volume method with refinement only in space and advanced with the Forward Euler method deserves a standalone explanation. Figure \ref{NonConvergentError} shows the variation of the actual error of the linear variable-coefficient advection equation \eqref{LinearAdvection1D_NumericalExperiment} (left) and the non-linear advection equation \eqref{NonLinearAdvection1D_NumericalExperiment} (right) with refinement only in space. The error decreases for the first 4 points at third-order (as verified from Figure \ref{Figure_Numerical_Convergence_Plots_3}) after which it becomes non-convergent due to the application of the slope limiter, the monotonicity-preserving strategies, and the dissipative local Lax-Friedrichs Riemann solver. More specifically, the error passes through phases of local maxima and minima. As a result, when we plot the norm of the difference in the error (or the numerical solution) for successive spatial resolutions, the curve is non-convergent after the first 3 points and follows a similar pattern, as observed in Figure \ref{Figure_Numerical_Convergence_Plots_3}.


To understand the role of the coefficients of the leading order terms of the error, we have fit the polynomial $\zeta_{\gamma} \Delta \xi^{\gamma} + \zeta_{\gamma+1} \Delta \xi^{\gamma+1}$ to (a) the error of the linear variable-coefficient advection equation~\eqref{LinearAdvection1D_NumericalExperiment} and the non-linear advection equation~\eqref{NonLinearAdvection1D_NumericalExperiment} with refinement in both space and time, and (b) the norm of the difference in the numerical solution of~\eqref{LinearAdvection1D_NumericalExperiment} and~\eqref{NonLinearAdvection1D_NumericalExperiment} with refinement only in space or only in time. For the advection equations~\eqref{LinearAdvection1D_NumericalExperiment} and~\eqref{NonLinearAdvection1D_NumericalExperiment}, Tables~\ref{Table_Linear_Numerical_Advection_Error_Coefficients} and~\ref{Table_NonLinear_Numerical_Advection_Error_Coefficients} list the coefficients $\left[\zeta_{\gamma},\zeta_{\gamma+1}\right]$. For refinement in both space and time and refinement only in space, $\Delta \xi = \Delta x$, and for refinement only in time, $\Delta \xi = \Delta t$. For the finite difference and the non-monotone finite volume methods, the resolutions have reached the asymptotic regime, and $\gamma = \alpha$ for refinement only in space, $\gamma = \beta$ for refinement only in time, and $\gamma = \min(\alpha,\beta)$ for refinement in both space and time, with $\alpha$ and~$\beta$ representing the spatial and temporal orders of accuracy. The convergence curve of the advection equations discretized with the non-monotone finite volume method and advanced with the explicit midpoint method and the second-order Adams-Bashforth methods, reach the asymptotic regime after the first three points for refinement in both space and time. So, we do not include these points for determining the best fit polynomial. Even though $\zeta_{\gamma+1} \Delta \xi^{\gamma+1} \ll \zeta_{\gamma} \Delta \xi^{\gamma}$ in the asymptotic regime, if the slope of the convergence curve is $\gamma+1$ instead of $\gamma$, it is immediately clear that we have not reached the aymptotic regime, and we would expect $\zeta_{\gamma+1}$ to be at least a few orders of magnitude larger than $\zeta_{\gamma}$. This is what we observe for the convergence of the advection equations discretized with the non-monotone finite volume method and advanced with the third-order Runge-Kutta or Adams-Bashforth methods, with refinement in both space and time. For these cases, $\zeta_{4}$ is three orders of magnitude larger than $\zeta_{3}$. If, however, the resolutions have reached the asymptotic regime, and the slope of the convergence curve is $\gamma$, we observe $\zeta_{\gamma+1}$ to be (a) less than $\zeta_{\gamma}$, or (b) of the same order of magnitude as $\zeta_{\gamma}$, or (c) at most one order of magnitude larger than $\zeta_{\gamma}$. For the monotone finite volume method, $\alpha$ and $\beta$ are replaced by their reduced equivalents, as observed in the convergence plots of Figure \ref{Figure_Numerical_Convergence_Plots_2}. If we observe order reduction in the convergence curves obtained with these monotone finite volume methods, and we use all points to obtain the best fit polynomial, $\zeta_{\gamma+1}$ can be a few orders of magnitude larger than $\zeta_{\gamma}$. Finally, the norm of the difference in the numerical solution with only spatial refinement is non-convergent for the monotone finite volume method advanced with Forward Euler after the first 3 points. So, we fit the polynomial to only these points.

Summarizing, we obtain the expected orders of convergence with the finite difference method, but not with the finite volume method for some of the time-stepping methods. This is because of two main reasons. First, the resolutions have not reached the asymptotic regime for some of these time-stepping methods, for example, with the advection equations discretized with the non-monotone finite volume method and advanced with the third-order Runge-Kutta or Adams-Bashforth method, for refinement in both space and time. As a result, the coefficients of the leading order terms in the truncation error predominate over the powers of $\Delta x$ and $\Delta t$. This trend is expected to reverse with spatial or temporal refinement as we approach the asymptotic regime. However, we were unable to do so with the non-monotone finite volume method without the solution becoming unstable. Now, the extent of refinement of the discretization parameters $\Delta x$ and $\Delta t$ required to reach the asymptotic regime depends on the problem being solved. In our effort to reach the asymptotic regime, we may keep refining $\Delta x$ and $\Delta t$, and eventually the numerical error with the reduced values of $\Delta x$ and $\Delta t$ can drop below machine precision. In practice, we may not reach the asymptotic regime, and the magnitude of the error is dictated by the coefficients of the leading order terms in the truncation error, rather than the powers of $\Delta x$ and $\Delta t$. If, however, the coefficients of the leading order terms in the truncation error dominate for the largest values of $\Delta x$ and $\Delta t$, and we approach (or reach) the asymptotic regime before the machine precision error dominates, we expect to obtain reduction in the convergence slope. This is observed with the refinement of the advection equations in both space and time, when they are (a) discretized in space with the non-monotone finite volume method, and advanced in time with the explicit midpoint and the second-order Adams-Bashforth methods, and (b) discretized in space with the monotone finite volume method, and advanced in time with Forward Euler, the explicit midpoint and the second-order Adams-Bashforth methods, with simultaneous refinement in space and time. The second reason for some of the convergence slopes obtained with the finite volume method not matching the theoretical predictions is the following one. The slope limiter, the monotonicity-preserving strategies, and the disspative Riemann solver, all of which enable us to employ a high order finite volume method and increase the resolution while ensuring numerical stability, drop the order of accuracy. So, for verification purposes, we should adhere to the non-monotone finite volume methods. In other words, we should refrain from adopting any strategy which can reduce the order of accuracy. But despite the order reduction due to these strategies adopted in the monotone finite volume methods, we observe that the optimum order of convergence at constant ratio of time step to cell width is obtained by a time-stepping method of at least the same order of accuracy as that of the spatial discretization.


\section{Conclusion} \label{sec:conclusion}

We have derived expressions for the local truncation error of generic and specific hyperbolic PDEs, consisting of linear and non-linear advection equations, advanced in time with a variety of time-stepping methods, belonging to the Method of Lines, e.g.~Forward Euler, predictor-corrector methods like Runge-Kutta and multistep methods like Adams-Bashforth. We used first-, second-, and third-order upwind spatial discretization on a uniform mesh. The local truncation error assumes the form 
\begin{equation}
\hat{\tau} = \Delta t \bigO\left({\Delta x}^{\alpha}\right) + {\Delta t}^2 \bigO\left({\Delta x}^{\alpha}\right) + {\Delta t}^3 \bigO\left({\Delta x}^{\alpha}\right) + \cdots + {\Delta t}^{\beta} \bigO\left({\Delta x}^{\alpha}\right) + \bigO\left({\Delta t}^{\beta+1}\right). \label{LocalTruncationErrorNumericalSolutionFinalForm_Conclusion}
\end{equation}
The form of the local truncation error does not depend on whether the advection equation is linear or non-linear, constant- or variable-coefficient, homogeneous or inhomogeneous, and whether the time-stepping method is explicit or implicit, predictor-corrector or multistep. The leading order terms of the local truncation error only depend on the orders of the spatial and temporal discretizations. If the PDE is reduced to an ODE by specifying all spatial gradients to be zero, the local truncation error reduces to that of the ODE, thereby attesting to the robustness of our theory. At a time horizon, the global truncation error is one order of $\Delta t$ less than its local counterpart, and assumes the form 
\begin{equation}
\hat{\tau}_G = \bigO\left({\Delta x}^{\alpha}\right) + \Delta t \bigO\left({\Delta x}^{\alpha}\right) + {\Delta t}^2 \bigO\left({\Delta x}^{\alpha}\right) + \cdots + {\Delta t}^{\beta-1} \bigO\left({\Delta x}^{\alpha}\right) + \bigO\left({\Delta t}^{\beta}\right) \approx \bigO\left({\Delta x}^{\alpha}\right) + \bigO\left({\Delta t}^{\beta}\right) \text{ for $\Delta t \ll 1$}. \label{GlobalTruncationErrorNumericalSolutionFinalForm_Conclusion}
\end{equation}
When performing convergence tests of a hyperbolic PDE with the assumptions of
\begin{enumerate}[label=(\alph*),noitemsep]
\item a stable numerical scheme,
\item the global solution error being of the same order of accuracy as the global truncation error,
\item having reached the asymptotic regime, where the truncation error is dominated by the powers of the cell width and the time step rather than their coefficients, before machine precision error dominates,
\end{enumerate}
the following hold:
\begin{enumerate}[label=(\roman*),noitemsep]
\item The order of convergence at constant ratio of time step to cell width is determined by the minimum of the orders of the spatial and temporal discretizations. So, for a spatial discretization of a given order, a time-stepping method of at least the same order should be employed to attain the optimum order of convergence, with the most computationally efficient choice being the order of the spatial discretization.
\item Convergence of the error norm cannot be guaranteed under only spatial or temporal refinement.
\item By plotting the difference in the numerical solution or the error between two successive pairs of spatial (or temporal) resolutions, the convergence rates of the spatial (or temporal) discretization can be determined.
\end{enumerate}
We have conducted numerical experiments with linear and non-linear advection equations to demonstrate and underline our theoretical findings. We have employed finite difference and finite volume spatial discretizations, and a variety of time-stepping methods including Forward Euler, Runge-Kutta, and Adams-Bashforth methods from second up to fourth order. With the finite difference and the majority of the non-monotone finite volume methods, the spatial and temporal resolutions have reached the asymptotic regime and the convergence rates match our theoretical predictions. However, for the finite volume method with some of the time-stepping methods, the resolutions do not reach the asymptotic regime before machine precision error takes over. Under such circumstances, the coefficients of the leading order terms in the truncation error cannot be ignored, and consequently (i)--(iii) may not necessarily hold. Moreover, the slope-limiting monotonicity-preserving strategies and the dissipation provided by the Riemann solver in the monotone finite volume method drops the spatial and temporal orders of accuracy. These are the practical aspects of numerical models we need to take into consideration, which are not covered by our theory. However, for a spatial discretization of a given order, we still observe that the optimum order of convergence is attained by a time-stepping method of at least the same order.

The points presented in this paper on how spatial and temporal convergence interact, specifically equations \eqref{LocalTruncationErrorNumericalSolutionFinalForm_Conclusion} and \eqref{GlobalTruncationErrorNumericalSolutionFinalForm_Conclusion}, are both straightforward and fundamental. Ongoing and future work includes extending our theory to parabolic PDEs and higher-order and spectral discretizations in space and time.


\section{Acknowledgements} \label{sec:acknowledgements}

Siddhartha Bishnu has been supported by the Scientific Discovery through Advanced Computing (SciDAC) projects LEAP (Launching an Exascale ACME Prototype) and CANGA (Coupling Approaches for Next Generation Architectures) under the U.S.~Department of Energy (DOE), Office of Science, Office of Biological and Environmental Research (BER). Mark Petersen was supported as part of the Energy Exascale Earth System Model (E3SM) project, also funded by the DOE BER. This research used resources provided by the Los Alamos National Laboratory Institutional Computing Program, which is supported by the U.S. Department of Energy National Nuclear Security Administration under Contract No.~89233218CNA000001. The authors thank Tomek Plewa, Darren Engwirda, Giacomo Capodaglio, and Pedro da Silva Peixoto for helpful discussions.








\end{document}